%% file: wvf4.3.tex
\documentclass[12pt]{report}
\usepackage{a4}
\usepackage{amsfonts,amsthm}
\usepackage{mathrsfs}
\newtheorem*{proposition}{Proposition}
\newtheorem*{lemma}{Lemma}
\newtheorem*{question}{Question/Conjecture}

\newtheorem{corollarynumber}{Corollary}
\theoremstyle{remark}
\newtheorem*{remark}{Remark}
\newtheorem*{example}{Example}
\newtheorem{examplenumber}{Example}
\newtheorem*{sketcht}{Sketch of proof}
\newtheorem{tablef}{Table}

\newenvironment{itemf}{\vspace{-4pt}\begin{itemize}\setlength{\parskip}{-4pt}}{\end{itemize}\vspace{-4pt}}

\newenvironment{sketche}{\begin{sketcht}}{\hspace{\stretch{1}}$\square$\end{sketcht}}

\newcommand{\m}{{}^-\!}
\newcommand{\p}{\,\grave{}}

\renewcommand{\mathcal}{\mathscr}
\renewcommand{\phi}{\varphi}

\begin{document}

\input{lunadiagrams}

\input{figures}

\title{An introduction to wonderful varieties\\ with many examples of type $\mathsf F_4$}
\author{P.\ Bravi and D.\ Luna}
\date{March 9, 2009}
\maketitle

\thispagestyle{empty}
\section*{}
\setcounter{page}{2}
\newpage
\setcounter{page}{1}
\section*{Introduction}

Let $G$ be a complex algebraic group. A smooth connected projective $G$-variety $X$ is called \textit{wonderful} of rank $r$, if $G$ stabilizes exactly $r$ irreducible divisors in $X$ having the following properties:
\begin{itemf}
\item[-] they are smooth and have a nonempty transversal intersection;
\item[-] two points of $X$ are on the same $G$-orbit if (and only if) they are contained in the same $G$-stable divisors. 
\end{itemf}

In this paper, unless otherwise stated, we will assume the group $G$ to be semisimple and simply connected, and the action of the center of $G$ on wonderful varieties to be trivial. Wonderful varieties of rank 0 are then just the (generalized) flag varieties. Other examples are the wonderful completions of (adjoint) homogeneous symmetric spaces (\cite{DP83}). 

In the recent years a general theory of wonderful varieties has been developed. Our aim in this paper is to give a somewhat unusual introduction to this theory, by presenting (all) examples when $G$ is simple of type $\mathsf F_4$, and studying many of them more in detail. This $G$, although of small rank, seems complicated enough to illustrate most of the phenomena present in the general case.

In Chapter~1, we recall some basic facts about wonderful varieties and their combinatorial invariants called \textit{spherical systems} (see \cite{Lu01} and \cite{Bri07a}). We very much emphasize \textit{spherical diagrams}: these diagrams provide a convenient way of visualizing spherical systems (just as Dynkin diagrams do for root systems), and will be omnipresent throughout the rest of the paper. We give the list of all 266 spherical diagrams of type $\mathsf F_4$:
\begin{itemf}
\item[-] 16 of rank 0 (the generalized flag varieties),
\item[-] 41 of rank 1,
\item[-] 61 of rank 2,
\item[-] 77 of rank 3,
\item[-] 71 of rank 4.
\end{itemf}

In Chapter~2, we present some further notions and results about wonderful varieties. Since our main intention here is only to prepare for the study of examples of type $\mathsf F_4$, this chapter is somewhat experimental and sketchy: some proofs are only outlined and several questions are left open (as already in Chapter~1). An exception is Section~\ref{sssphericalclosure}, where we reexamine and clarify the notion of spherical closure (introduced in \cite{Lu01}, Section~6), and obtain a general classification of spherical orbits in all simple projective $G$-spaces, which seems to be new.

In Chapter~3, we reach the main object of our study: in ten small sections, we examine more closely several of the 266 wonderful varieties of type $\mathsf F_4$, those we believe to be the most interesting. In a last section, we discuss also some examples which are not of type $\mathsf F_4$.

\bigskip

Let us say some words about the place of wonderful varieties in mathematics, to explain the usefulness of a general theory of wonderful varieties. 

A good general notion of algebraic varieties ``having a big algebraic symmetry group'' is that of spherical varieties, i.e.\ of normal algebraic $G$-varieties, under a reductive connected group $G$, which are almost homogeneous under a Borel subgroup of $G$ (see \cite{Bri97} or \cite{T06} for an introduction). Spherical varieties are strongly related to (are the algebraic analog of) real symplectic manifolds equipped with a multiplicity-free Hamiltonian action of a compact Lie group (see \cite{Wo98}).

Wonderful varieties are spherical and, moreover, play a central role inside the theory of spherical varieties: indeed, to every spherical $G$-variety $Y$, one can associate, in a functorial way, a wonderful variety $\mathbf Y$, and $Y$ is then determined by the spherical system of $\mathbf Y$ together with some additional combinatorial data (see \cite{Lu01}, Theorem~3). This means that many properties of spherical varieties can be understood in geometrical terms of wonderful varieties. So it is not surprising that wonderful varieties play an important role in the classification of model homogeneous spaces (\cite{GZ84,AHV98,Lu07}), or in the multiplicity-free case of invariant Hilbert schemes (\cite{AB05,BC08}), or in the study of total coordinate rings of spherical varieties (\cite{Bri07a}).

Furthermore, wonderful varieties have some features of symmetric spaces: for instance, they naturally come with an interesting ``little Weyl group'' (\cite{Bri90,Kn95}). Lastly, let us mention that wonderful varieties are special cases of log homogeneous varieties (\cite{Bri07b,Bri08}).

%\clearpage
\chapter{Wonderful varieties and\\ spherical systems}

In the first section, we will recall some notions which can be naturally attached to each wonderful variety (spherical roots, colors, Cartan pairing, \ldots ), with appropriate notations, and explain how they are related to the combinatorial invariants of $G$ (root system and Cartan matrix). This will prepare the more formal approach to the combinatorial invariants of wonderful varieties (i.e.\ spherical systems) given in the second section. In the third and last section, we will recall some facts about the simple group of type $\mathsf F_4$, and give the list of all 266 spherical systems of type $\mathsf F_4$.

In what follows, we fix two opposite Borel subgroups $B$ and $\m B$ of $G$, so that $T=B\cap \m B$ is a maximal torus of $G$. If $K$ is any affine algebraic group, we write $K^r$ (resp.\ $K^u$) for the radical (resp.\ the unipotent radical) of $K$, and $\Xi(K)$ for the group of characters of $K$. We identify $\Xi(B)$ and $\Xi(T)$, and denote $R\subset\Xi(T)$ the root system of $G$ and $S\subset R$ the basis of $R$ corresponding to $B$. The elements of $S$ are traditionally called \textit{simple roots}. For all $\alpha\in S$, we denote by $\alpha^\vee$ the co-root of $\alpha$, and by $P_\alpha\supset B$ and $\m P_\alpha\supset \m B$ the corresponding minimal parabolic subgroups of $G$. Remember that the integers $\langle\alpha^\vee,\beta\rangle$, for $(\alpha,\beta)\in S\times S$, are the coefficients of the Cartan matrix of $G$. We denote by $\Omega$ the set of fundamental weights of $G$; so $\mathbb N\Omega$ can be considered as the set of dominant weights.

Remember that an algebraic $G$-variety $Y$ is called \textit{spherical}, if it is normal and if $B$ has an open (dense) orbit in $Y$. It is known that $B$ has only finitely many orbits in spherical varieties. A \textit{color} of a spherical variety $Y$ is by definition an irreducible $B$-stable but not $G$-stable divisor in $Y$. Let us denote by $\Delta_Y$ the set of colors of $Y$ and by $\Xi(Y)$ the free abelian group of weights of $B$ in the field of rational functions on $Y$.

It is known that wonderful varieties are spherical.

\vspace{6ex}\section{Basic definitions and properties}\label{ssbasic}

In what follows $X$ will always denote a wonderful $G$-variety.

\subsection{Spherical roots}

Let $z$ be the unique point in $X$ fixed by $\m B$. \textit{Spherical roots} of $X$ are by definition weights of $T$ in $T_zX/T_z(G.z)$ (the normal space at $z$ of the closed orbit $G.z$ in $X$). Let us denote by $\Sigma_X$ the set of spherical roots of $X$. Since the action of the center of $G$ is assumed to be trivial, one has $\Sigma_X\subset\mathbb NS$. The smoothness of $X$ implies that $\Sigma_X$ is a basis of $\Xi(X)$.

For every $\sigma\in\Sigma_X$, let us denote by $D^\sigma$ the unique $G$-stable divisor such that $\sigma$ is the weight of $T$ in $T_zX/T_zD^\sigma$. One obtains in this way a natural bijection between $\Sigma_X$ and the set of irreducible $G$-stable divisors of $X$. In particular, the cardinality of $\Sigma_X$ equals the rank of $X$. For all $\Sigma'\subset\Sigma_X$, we set $X_{\Sigma'}=\cap_{\sigma\in\Sigma_X\setminus\Sigma'}D^\sigma$; each $X_{\Sigma'}$ is a wonderful $G$-subvariety of $X$ whose set of spherical roots is $\Sigma'$. One obtains in this way a bijection between the set of subsets of $\Sigma_X$ and the set of closed irreducible $G$-stable subsets of $X$.

The spherical root of any wonderful $G$-variety of rank 1 is also called a \textit{spherical root of $G$}. Let us denote by $\Sigma(G)$ the (finite) set of spherical roots of $G$. Since the spherical roots of $X$ are exactly those of the wonderful varieties $X_{\{\sigma\}}$ for all $\sigma\in\Sigma_X$, one has $\Sigma_X\subset\Sigma(G)$.

\subsection{Colors and simple roots}\label{ssscolors}

Let us now explain the relations between the set of colors $\Delta_X$ and the set of simple roots $S$.

For all $\alpha\in S$, let us set $\Delta(\alpha)=\{D\in\Delta_X: P_\alpha.D\neq D\}$. One has $\Delta_X=\cup_{\alpha\in S}\Delta(\alpha)$ and $0\leq\mathrm{card}(\Delta(\alpha))\leq2$.

Define $S^p_X$ as the set of simple roots $\alpha$ such that $\Delta(\alpha)=\emptyset$. The parabolic subgroup of $G$ generated by the $P_\alpha$'s, $\alpha\in S^p_X$, is sometimes denoted by $P=P_X$; it can also be defined as the stabilizer of the open orbit of $B$ in $X$. One has the useful formula: $\dim X=\mathrm{card}\, \Sigma_X+\dim P_X^u$. The parabolic subgroup $\m P$ of $G$ (opposite with respect to $T$) is equal to the isotropy group at $z$.

Set $S^a_X=S\cap\Sigma_X$. One has $\mathrm{card}(\Delta(\alpha))=2$ if and only if $\alpha\in S^a_X$. Set $S^{2a}_X=\{\alpha\in S: 2\alpha\in\Sigma_X\}$ and  $S^b_X=S\setminus(S^p_X\cup S^a_X\cup S^{2a}_X)$.

Furthermore, let us set 
\begin{itemf}
\item[-] $\Delta^a_X=\cup \Delta(\alpha)$, $\alpha\in S^a_X$; 
\item[-] $\Delta^{2a}_X=\cup \Delta(\alpha)$, $\alpha\in S^{2a}_X$; 
\item[-] $\Delta^b_X=\cup \Delta(\alpha)$, $\alpha\in S^b_X$. 
\end{itemf}
The union $\Delta_X=\Delta^a_X\cup\Delta^{2a}_X\cup\Delta^b_X$ is disjoint.

\subsection{Cartan pairing}\label{ssscartan}

The classes $[D]$, $D\in\Delta_X$, form a basis of $\mathrm{Pic}(X)$. The pairing $c\colon\Delta_X\times\Sigma_X\to\mathbb Z$ given by $[D^\sigma]=\sum_{D\in\Delta_X} c(D,\sigma)[D]$ is called the \textit{Cartan pairing} of $X$.

The natural map $G/\m B\to G.z\subset X$ induces a group homomorphism $\omega$ from $\mathbb Z\Delta_X\cong\mathrm{Pic}(X)$ to $\mathrm{Pic}(G/\m B)\cong\mathbb Z\Omega$, the weight lattice of $G$. The way we labeled the $G$-stable divisors of $X$ implies $\omega([D^\sigma])=\sigma$.

The Cartan pairing of $X$ is related to the Cartan matrix of $G$ as follows:
\begin{itemf}
\item[-] if $\alpha\in S^a_X$ and if we denote by $D^+_\alpha$, $D^-_\alpha$ the two colors in $\Delta(\alpha)$, 
\item[] then $c(D^+_\alpha,\sigma)+c(D^-_\alpha,\sigma)=\langle\alpha^\vee,\sigma\rangle$, for all $\sigma\in\Sigma_X$;
\item[-] if $\alpha\in S^{2a}_X$ and if we denote by $D_{2\alpha}$ the unique color in $\Delta(\alpha)$, 
\item[] then $c(D_{2\alpha},\sigma)=\langle\alpha^\vee,\sigma\rangle/2$, for all $\sigma\in\Sigma_X$;
\item[-] if $\alpha\in S^b_X$ and if we denote by $D_\alpha$ the unique color in $\Delta(\alpha)$, 
\item[] then $c(D_\alpha,\sigma)=\langle\alpha^\vee,\sigma\rangle$, for all $\sigma\in\Sigma_X$. 
\end{itemf}

Moreover, if $\alpha, \beta\in S^b_X$, then $D_\alpha=D_\beta$ if and only if $\alpha\perp\beta$ and $\alpha+\beta\in\Sigma_X$.

\subsection{Localization (at a parabolic subgroup $Q$)}\label{localization}

Let $Q$ be a parabolic subgroup of $G$ containing $B$. We will denote by $\m Q$ the opposite parabolic with respect to $T$, and by $M=Q\cap \m Q$ the common Levi subgroup of $Q$ and $\m Q$.

Let $X$ be a wonderful $G$-variety. Let us denote by $Z$ the set of points of $X$ fixed by $\m Q^r$ (the radical of $\m Q$). This variety $Z$ is known to be a wonderful $\m Q/ \m Q^r$-variety called the wonderful variety obtained by \textit{localization} of $X$ at $Q$. 

\subsection{Parabolic induction}\label{ssparabolic}

Conversely, from a wonderful $\m Q/ \m Q^r$-variety $Z$, one gets a
wonderful $G$-variety by $X=G\ast_{\m Q}Z$, the quotient of $G\times
Z$ under the action of $ \m Q$: $q.(g,y)=(gq^{-1},q.y)$, $q\in \m Q$,
$g\in G$, $y\in Z$. The localization of $G\ast_{\m Q}Z$ at $Q$ is clearly isomorphic to $Z$. This construction is called \textit{parabolic induction}. A wonderful variety which cannot be obtained by proper parabolic induction is sometimes called \textit{cuspidal}.

\subsection{Wonderful varieties of rank 1}\label{ssswonderfulvarietiesofrank1}

All wonderful varieties of rank 1 are well known (\cite{Ah83}). They are obtained by parabolic induction from a list of cuspidal wonderful varieties of rank 1. In particular, for every $G$, the set $\Sigma(G)$ (the set of spherical roots of $G$) is known. 

A spherical root $\sigma$ has one of the following shapes. If its support is connected (except one case, it is always connected) then we write $\sigma=n_1\alpha_1+\ldots+n_r\alpha_r$, labeling the simple roots $\alpha_1,\ldots,\alpha_r$ in $\mathrm{supp}(\sigma)$ as in Bourbaki. We group them according to their type of support.
\bigskip

\begin{center}
\begin{tabular}{ll}
type of $\mathrm{supp}(\sigma)$&shape of $\sigma$\\
\hline
$\mathsf A_1$&$\alpha_1$\\
&$2\alpha_1$\vspace{0.1cm}\\
$\mathsf A_1\times\mathsf A_1$&$\alpha_1+\alpha_1'$\vspace{0.1cm}\\
$\mathsf A_r$, $r\geq2$&$\alpha_1+\ldots+\alpha_r$\vspace{0.1cm}\\
$\mathsf A_3$&$\alpha_1+2\alpha_2+\alpha_3$\vspace{0.1cm}\\
$\mathsf B_r$, $r\geq2$&$\alpha_1+\ldots+\alpha_r$\\
&$2\alpha_1+\ldots+2\alpha_r$\vspace{0.1cm}\\
$\mathsf B_3$&$\alpha_1+2\alpha_2+3\alpha_3$\vspace{0.1cm}\\
$\mathsf C_r$, $r\geq3$&$\alpha_1+2\alpha_2+\ldots+2\alpha_{r-1}+\alpha_r$\vspace{0.1cm}\\
$\mathsf D_r$, $r\geq4$&$2\alpha_1+\ldots+2\alpha_{r-2}+\alpha_{r-1}+\alpha_r$\vspace{0.1cm}\\
$\mathsf F_4$&$\alpha_1+2\alpha_2+3\alpha_3+2\alpha_4$\vspace{0.1cm}\\
$\mathsf G_2$&$\alpha_1+\alpha_2$\\
&$2\alpha_1+\alpha_2$\\
&$4\alpha_1+2\alpha_2$
\end{tabular}
\end{center}
\bigskip

For every wonderful $G$-variety $X$ of rank 1, the couple ($\sigma_X$, $S^p_X$) (where $\sigma_X$ is the spherical root of $X$ and $S^p_X$ is as defined in \ref{ssscolors}) determines $X$ up to $G$-isomorphism.

We call a couple ($\sigma$, $S^p$) (where $\sigma\in\Sigma(G)$ and $S^p\subset S$) \textit{compatible}, if it is obtained as above from a $G$-wonderful variety of rank 1. 

A couple ($\sigma$, $S^p$) is compatible, if and only if $S^{pp}(\sigma)\subset S^p\subset S^p(\sigma)$, where $S^p(\sigma)$ denotes the set of simple roots orthogonal to $\sigma$, and where $S^{pp}(\sigma)$ is equal to
\begin{itemf}
\item[-] $S^p(\sigma)\cap\mathrm{supp}(\sigma)\setminus\{\alpha_r\}$ if $\sigma=\alpha_1+\ldots+\alpha_r$ with support of type $\mathsf B_r$, or
\item[-] $S^p(\sigma)\cap\mathrm{supp}(\sigma)\setminus\{\alpha_1\}$ if $\sigma$ has support of type $\mathsf C_r$, or 
\item[-] $S^p(\sigma)\cap\mathrm{supp}(\sigma)$ otherwise.
\end{itemf}

\vspace{6ex}\section{Spherical systems}

In this section we introduce the notion of \textit{spherical system}, as a purely abstract (axiomatic) version of the invariants introduced in Section~\ref{ssbasic}, and state then a conjecture on the combinatorial classification of wonderful varieties.

Spherical roots are already, by definition, combinatorial invariants (related to the root system of $G$), but colors a priori are not. In order to give a combinatorial meaning to colors, we retain only their relations with spherical roots via the Cartan pairing. Moreover, we have seen in Section~\ref{ssbasic} that, if one knows the set of spherical roots and the isomorphism class of the closed orbit (i.e.\ the set $S^p$), then one knows already automatically the \textit{abstract} colors of types $2a$ and $b$. These remarks lead naturally to the following definition.

\subsection{The definition of spherical systems}

We call \textit{spherical system} of $G$ any triplet $\mathcal S=(\Sigma, S^p, \mathbf A)$, where
\begin{itemf}
\item[-] $\Sigma$ is a subset (without proportional elements) of $\Sigma(G)$,
\item[-] $S^p$ is a subset of $S$,
\item[-] $\mathbf A$ is a finite multi-subset of $(\mathbb Z\Sigma)^\ast=\mathrm{Hom}_\mathbb Z(\mathbb Z\Sigma,\mathbb Z)$ (i.e.\ $\mathbf A$ is a finite abstract set, together with a map $c\colon \mathbf A\to(\mathbb Z\Sigma)^\ast$, which we will consider also as a pairing $c\colon \mathbf A\times\Sigma\to\mathbb Z$) 
\end{itemf}
satisfying the following properties (axioms):
\begin{itemf}
\item[(S)] $S^p$ is compatible with all $\sigma\in\Sigma$;
\item[(A1)] for all $\delta\in\mathbf A$ and $\sigma\in\Sigma$, one has $c(\delta,\sigma)\leq1$, and $c(\delta,\sigma)=1$ implies $\sigma\in S\cap\Sigma$;
\item[(A2)] for all $\alpha\in S\cap\Sigma$, the set $\mathbf A(\alpha)=\{\delta\in\mathbf A : c(\delta,\sigma)=1\}$ contains exactly two elements, and if $\mathbf A(\alpha)=\{\delta_\alpha^+,\delta_\alpha^-\}$ then $c(\delta_\alpha^+,\sigma)+c(\delta_\alpha^-,\sigma)=\langle\alpha^\vee,\sigma\rangle$, for all $\sigma\in\Sigma$;
\item[(A3)] $\mathbf A$ is the union of the $\mathbf A(\alpha)$'s, for $\alpha\in S\cap\Sigma$;
\item[($\Sigma$1)] if $2\alpha\in\Sigma\cap2S$, then $\langle\alpha^\vee,\sigma\rangle/2$ is a non-positive integer, for all $\sigma\in\Sigma\setminus\{2\alpha\}$;
\item[($\Sigma$2)] if $\alpha+\beta\in\Sigma$, with $\alpha,\beta\in S$ and $\alpha\perp\beta$, then $\langle\alpha^\vee,\sigma\rangle=\langle\beta^\vee,\sigma\rangle$ for all $\sigma\in\Sigma$.
\end{itemf}

If $X$ is any wonderful $G$-variety, then $\mathcal S_X=(\Sigma_X, S^p_X, \Delta^a_X)$ is a spherical system.

Let us explain here how all the colors and the full Cartan pairing can be recovered from the spherical system, since this will play a role in what follows. Let $\mathcal S=(\Sigma,S^p,\mathbf A)$ be a spherical system. Define $S^{2a}=\{\alpha\in S: 2\alpha\in\Sigma\cap2S\}$, $S^b=S\setminus(S^p\cup(S\cap\Sigma)\cup S^{2a})$ and $\Delta$ as the disjoint union $\mathbf A\cup S^{2a}\cup(S^b/\sim)$, where $\alpha,\beta$ are identified in $S^b$ if they are orthogonal and $\alpha+\beta\in\Sigma$. It is then clear how to define $c\colon\Delta\to(\mathbb Z\Sigma)^\ast$, that is, by the map $c$ already given on $\mathbf A$ and by the corresponding half-co-roots and co-roots on $S^{2a}$ and $S^b/\sim$. In this way one clearly recovers also the decomposition $\Delta=\cup_{\alpha\in S}\Delta(\alpha)$.

\subsection{The main question}\label{question}

Some years ago, the following question has been formulated:

\begin{question}[\cite{Lu01}]
Are wonderful varieties classified by their spherical systems?
\end{question} 

This question has been answered positively in many cases
(\cite{W96,Lu01,BP05,Bra07,BC08,Cu08}). The ``uniqueness part'' (i.e.\
that two wonderful varieties having same spherical system are
$G$-isomorphic) has been proved last year in general (this follows
from results of I.V.~Losev in \cite{Lo07}). Although many spherical
systems have been geometrically realized (i.e.\ as spherical systems
of wonderful varieties), at present no general proof for this
``existence part'' exists in the literature (even for the case of type
$\mathsf F_4$).

\bigskip
\subsection{Spherical systems and localization}

Let $\mathcal S=(\Sigma, S^p, \mathbf A)$ be a spherical system.

For every $\Sigma'\subset\Sigma$, we define another spherical system by $\mathcal S_{\Sigma'}=(\Sigma',S^p,\mathbf A')$, where $\mathbf A'$ is the union of the $\mathbf A(\alpha)$'s, $\alpha\in S\cap\Sigma'$, and where the Cartan pairing $c\colon\mathbf A'\times\Sigma'\to\mathbb Z$ is obtained by restriction from $c\colon\mathbf A\times\Sigma\to\mathbb Z$. We will say that $\mathcal S_{\Sigma'}$ is obtained by \textit{localization} of $\mathcal S$ with respect to $\Sigma'$.

For every $S'\subset S$, we define still another spherical system $\mathcal S_{S'}=(\Sigma',(S^p)',\mathbf A')$, where
\begin{itemf}
\item[-] $\Sigma'=\Sigma\cap\mathbb NS'$;
\item[-] $(S^p)'=S^p\cap S'$;
\item[-] $\mathbf A'$ is the union of the $\mathbf A(\alpha)$'s, $\alpha\in S'\cap\Sigma$ and the Cartan pairing\\ $c\colon\mathbf A'\times\Sigma'\to\mathbb Z$ is obtained by restriction from $c\colon\mathbf A\times\Sigma\to\mathbb Z$.
\end{itemf}
We will say that $\mathcal S_{S'}$ is obtained by \textit{localization} of $\mathcal S$ with respect to $S'$. Notice that $\mathcal S_{S'}$ is a spherical system of the root system $R\cap\mathbb ZS'$.

Let $X$ be a wonderful $G$-variety having spherical system $\mathcal S=(\Sigma,S^p,\mathbf A)$. Let $Q$ be a parabolic subgroup of $G$ containing $B$, and let $S'$ be the set of $\alpha\in S$ such that $P_\alpha\subset Q$. Let us denote by $X_{S'}=Z$, the localization of $X$ at $Q$ introduced in \ref{localization}.

%\newpage
\begin{proposition}
\begin{itemf}
\item[]
\item[1)] The spherical system of $X_{\Sigma'}$ is $\mathcal S_{\Sigma'}$.
\item[2)] The spherical system of $X_{S'}$ is $\mathcal S_{S'}$.
\end{itemf}
\end{proposition}

\begin{sketche}
One uses the following interpretation of the Cartan pairing.

(*) Let $\alpha\in S\cap\Sigma$. Then $X_{\{\alpha\}}$, the localization of $X$ with respect to $P_\alpha$,
is isomorphic to $\mathbb P^1 \times \mathbb P^1$, where $\m P_\alpha$ acts via the natural morphism
$\m P_\alpha\to \mathrm{PGL}(2)$. From this follows that $T$ has four fixed points in $X_{\{\alpha\}}$:
$z$, $s_\alpha.z$, and two other points $z_\alpha^+$, $z_\alpha^-$ exchanged by $s_\alpha$ (the involution
of $\mathrm{N}_G(T)/T$ associated to $\alpha$). For any $\sigma\in\Sigma\setminus\{\alpha\}$, the two points $z_\alpha^\pm$
are contained in $D^\sigma$, and one can show that the weight of $T$ in the
normal bundles of $D^\sigma$ in $X$ at $z_\alpha^\pm$ is given by $\sigma - c(D_\alpha^\pm , \sigma)\alpha$.

1) All assertions on the spherical system of $X_{\Sigma'}$ are easy, except those
concerning the Cartan pairing, which follow from (*).

2) One can identify the set of roots of $\m Q/\m Q^r$ to $R\cap\mathbb ZS'$, which is also the
set of $\beta\in R$ trivial on $M^r$. The variety $Z$ can also be characterized as the
connected component of the set of points of $X$ fixed by $M^r$, containing $z$
(the unique point of $X$ fixed by $\m B$). From this follows easily that the
spherical roots of the wonderful $\m Q/\m Q^r$-variety $Z$ are those of $\Sigma\cap \mathbb NS'$.
The other assertions are either easy, or follow from (*).
\end{sketche}

\bigskip
\textit{Remarks}

1) From the proposition follows: if a spherical system is geometrically realizable, then so are all its localizations.

2) The conjecture in \ref{question} can be reformulated in more geometrical terms. 
Let us define the \textit{essential skeleton} of a wonderful $G$-variety, 
as the union
of all wonderful $G$-subvarieties which either have rank 1, or have rank 2
and contain at least one simple root as spherical root. Part 1) of the proposition implies that the
essential skeleton determines the spherical system; conversely, since
wonderful varieties of rank 1 and 2 are known (\cite{W96}), by gluing some of
these together, one can associate to each spherical system its essential
skeleton. Then one can ask: does every essential skeleton come from a (unique) wonderful
variety? 

\bigskip
\subsection{Diagrams}\label{sssdiagrams}

We will now introduce (spherical) \textit{diagrams}, which are a way to visualize
spherical systems (exactly as Dynkin diagrams allow to visualize root
systems). They are obtained by adding information to the Dynkin diagram
of the root system $R$.

Here is our way to represent spherical roots on the Dynkin diagram (see also \ref{ssswonderfulvarietiesofrank1}):

\begin{center}
\begin{tabular}{cl}
diagram&spherical root\\
\hline
\begin{picture}(900,2100)\put(0,0){\usebox{\aone}}\end{picture}&\begin{picture}(6000,1800)\put(0,600){$\alpha_1$}\end{picture}\\
\begin{picture}(900,1800)\put(0,0){\usebox{\aprime}}\end{picture}&\begin{picture}(6000,1800)\put(0,600){$2\alpha_1$}\end{picture}\\
\begin{picture}(3000,1800)\put(300,0){\line(1,0){2400}}\multiput(0,0)(2400,0){2}{\put(300,0){\line(0,1){600}}\put(300,900){\circle{600}}\put(300,900){\circle*{150}}}\end{picture}&\begin{picture}(6000,1800)\put(0,600){$\alpha_1+\alpha_1'$}\end{picture}\\
\begin{picture}(6000,1800)\put(0,600){\usebox{\mediumam}}\end{picture}&\begin{picture}(6000,1800)\put(0,600){$\alpha_1+\ldots+\alpha_r$}\end{picture}\\
\begin{picture}(3600,1800)\put(0,600){\usebox{\dthree}}\end{picture}&\begin{picture}(6000,1800)\put(0,600){$\alpha_1+2\alpha_2+\alpha_3$}\end{picture}\\
\begin{picture}(7500,1800)\put(0,600){\usebox{\shortbm}}\end{picture}&\begin{picture}(6000,1800)\put(0,600){$\alpha_1+\ldots+\alpha_r$}\end{picture}\\
\begin{picture}(7500,1800)\put(0,600){\usebox{\shortbprimem}}\end{picture}&\begin{picture}(6000,1800)\put(0,600){$2\alpha_1+\ldots+2\alpha_r$}\end{picture}\\
\begin{picture}(3900,1800)\put(0,600){\usebox{\bthirdthree}}\end{picture}&\begin{picture}(6000,1800)\put(0,600){$\alpha_1+2\alpha_2+3\alpha_3$}\end{picture}\\
\begin{picture}(9000,1800)\put(0,600){\usebox{\shortcm}}\end{picture}&\begin{picture}(6000,1800)\put(0,600){$\alpha_1+2\alpha_2+\ldots+2\alpha_{r-1}+\alpha_r$}\end{picture}\\
\begin{picture}(6900,2400)\put(0,0){\usebox{\shortdm}}\end{picture}&\begin{picture}(14400,2400)\put(0,900){$2\alpha_1+\ldots+2\alpha_{r-2}+\alpha_{r-1}+\alpha_r$}\end{picture}\\
\begin{picture}(5700,1800)\put(0,600){\usebox{\ffour}}\end{picture}&\begin{picture}(6000,1800)\put(0,600){$\alpha_1+2\alpha_2+3\alpha_3+2\alpha_4$}\end{picture}\\
\begin{picture}(2400,1800)\put(0,600){\usebox{\gsecondtwo}}\end{picture}&\begin{picture}(6000,1800)\put(0,600){$\alpha_1+\alpha_2$}\end{picture}\\
\begin{picture}(2400,1800)\put(0,600){\usebox{\gtwo}}\end{picture}&\begin{picture}(6000,1800)\put(0,600){$2\alpha_1+\alpha_2$}\end{picture}\\
\begin{picture}(2400,1800)\put(0,600){\usebox{\gprimetwo}}\end{picture}&\begin{picture}(6000,1800)\put(0,600){$4\alpha_1+2\alpha_2$}\end{picture}\\
%\hline
\end{tabular}
\end{center}
\bigskip

Let us now explain how to represent a spherical system $\mathcal S = (\Sigma , S^p , \mathbf A)$. One begins by representing all its spherical roots. Then one represents $S^p$ by adding some (not shadowed) circles around vertices in such a way that $S^p$ becomes exactly the set of vertices having no circles around, below or above. If $S\cap\Sigma = \emptyset$, we are done: the result is what we call a diagram, which allows one to visualize spherical systems of the form $\mathcal S = (\Sigma , S^p , \emptyset)$.

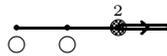
\bigskip
\begin{examplenumber}
\[\begin{picture}(5700,1200)\multiput(300,900)(1800,0){2}{\usebox{\segm}}\put(3900,900){\usebox{\rightbisegm}}\multiput(0,0)(1800,0){2}{\usebox{\aprime}}\put(3600,600){\usebox{\GreyCircleTwo}}\end{picture}\]
Here $G$ is of type $\mathsf B_4$. There are 3 spherical roots, $\sigma_1=2\alpha_1$, $\sigma_2=2\alpha_2$ and $\sigma_3=2\alpha_3+2\alpha_4$. Notice that, for $\alpha_1$ and $\alpha_2$, Axiom~$\Sigma$1 is satisfied. Moreover, $S^p=\{\alpha_4\}$. 
\end{examplenumber}

\bigskip
\begin{examplenumber}
\[\begin{picture}(6000,1200)\multiput(300,900)(1800,0){2}{\usebox{\segm}}\put(3900,900){\usebox{\leftbisegm}}\multiput(300,900)(5400,0){2}{\circle{600}}\multiput(300,0)(5400,0){2}{\line(0,1){600}}\put(300,0){\line(1,0){5400}}\put(1800,600){\usebox{\atwo}}\end{picture}\]
Here $G$ is of type $\mathsf C_4$. There are 2 spherical roots, $\sigma_1=\alpha_1+\alpha_4$ and $\sigma_2=\alpha_2+\alpha_3$. Notice that, for $\alpha_1$ and $\alpha_4$, Axiom~$\Sigma$2 is satisfied. Here $S^p=\emptyset$.
\end{examplenumber}

If $S\cap\Sigma\neq\emptyset$, more information is needed. The set $S\cap\Sigma$ corresponds to the set of vertices which have circles above and below. For each $\alpha\in S\cap\Sigma$, we identify these two circles with the elements of $\mathbf A(\alpha)$, the circle above to $\delta_\alpha^+$, where $\delta_\alpha^+$ is chosen such that $c(\delta_\alpha^+,\sigma)\in\{1,0-1\}$, for every spherical root $\sigma$. Then we join circles in different $\mathbf A(\alpha)$'s by lines, if they correspond to the same element in $\mathbf A$. Finally, for every spherical root $\sigma$ not orthogonal to $\alpha$ such that $c(\delta_\alpha^+, \sigma) = -1$, we add an arrow of the form $<$ or $>$, starting from the circle corresponding to $\delta_\alpha^+$, and pointing toward $\sigma$. This can always be done, and the set $\mathbf A$ and the Cartan pairing $c\colon \mathbf A \times\Sigma\to \mathbb Z$ is then determined by Axiom~A2.

\bigskip
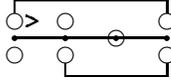
\begin{examplenumber}
\[\begin{picture}(6000,2700)\multiput(300,1350)(1800,0){3}{\usebox{\segm}}\multiput(0,450)(1800,0){2}{\usebox{\aone}}\put(5400,450){\usebox{\aone}}\put(3900,1350){\circle{600}}\multiput(300,2700)(5400,0){2}{\line(0,-1){450}}\put(300,2700){\line(1,0){5400}}\multiput(2100,0)(3600,0){2}{\line(0,1){450}}\put(2100,0){\line(1,0){3600}}\put(700,1750){\usebox{\toe}}\end{picture}\]
Here $G$ is of type $\mathsf A_4$. There are 3 spherical roots, all of them are simple roots: $\Sigma=\{\alpha_1,\alpha_2,\alpha_4\}$; $S^p=\emptyset$.
The set $\mathbf A$ has 4 elements: $\delta_{\alpha_1}^+=\delta_{\alpha_4}^+\in\mathbf A(\alpha_1)\cap\mathbf A(\alpha_4)$, $\delta_{\alpha_1}^-\in\mathbf A(\alpha_1)$, $\delta_{\alpha_2}^+\in\mathbf A(\alpha_2)$, $\delta_{\alpha_2}^-=\delta_{\alpha_4}^-\in\mathbf A(\alpha_2)\cap\mathbf A(\alpha_4)$. 
Since $\delta_{\alpha_1}^+$ belongs to $\mathbf A(\alpha_1)\cap\mathbf A(\alpha_4)$, $c(\delta_{\alpha_1}^+,\alpha_1)=c(\delta_{\alpha_1}^+,\alpha_4)=1$; since there is an arrow from $\delta_{\alpha_1}^+$ to $\alpha_2$, $c(\delta_{\alpha_1}^+,\alpha_2)=-1$. Then $c(\delta_{\alpha_1}^-)$ is determined by Axiom~A2, since $\{\delta_{\alpha_1}^+,\delta_{\alpha_1}^-\}=\mathbf A(\alpha_1)$ and $c(\delta_{\alpha_1}^+,-)+c(\delta_{\alpha_1}^-,-)=\langle\alpha^\vee,-\rangle$. Analogously, $c(\delta_{\alpha_4}^-)$ is determined since $\delta_{\alpha_1}^+=\delta_{\alpha_4}^+$. Finally, once $c(\delta_{\alpha_4}^-)$ is determined, so is $c(\delta_{\alpha_2}^+)$. The (restricted) Cartan pairing is then as follows:
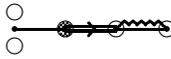
\[\begin{array}{r|rrr}c(-,-)&\alpha_1&\alpha_2&\alpha_4\\\hline \delta_{\alpha_1}^+&1&-1&1\\\delta_{\alpha_1}^-&1&0&-1\\\delta_{\alpha_2}^+&0&1&-1\\\delta_{\alpha_2}^-&-1&1&1\\
\end{array}\]
\end{examplenumber}

\bigskip
\begin{examplenumber}
\[\begin{picture}(6000,1800)\put(300,900){\usebox{\dynkinf}}\put(3600,600){\usebox{\atwo}}\put(0,0){\usebox{\aone}}\put(1800,600){\usebox{\GreyCircle}}\end{picture}\]
Here $G$ is of type $\mathsf F_4$. There are 3 spherical roots, $\sigma_1=\alpha_1$, $\sigma_2=\alpha_2+\alpha_3$ and $\sigma_3=\alpha_3+\alpha_4$, $S^p=\emptyset$, $\mathbf A=\{\delta_{\alpha_1}^+,\delta_{\alpha_1}^-\}$. Since there is no arrow in the diagram, $c(\delta_{\alpha_1}^+,\sigma_i)=0$, for $i=2,3$. 
\end{examplenumber}

At this stage, we should remark that diagrams not only allow to visualize spherical systems, but directly the set of all colors $\Delta$ and the (full) Cartan
pairing \mbox{$c\colon \Delta\times\Sigma\to\mathbb Z$}. Indeed, the way we have defined them, there is a natural bijection between the set $\Delta$ and the set of equivalence classes of circles of a diagram (two circles are equivalent if they are joined by a line); moreover, the Cartan pairing can be read off easily from the diagram.

Let us look again at the examples above.

In Example 1, there are 3 colors, 2 of type $2a$, $\delta_{2\alpha_1}\in\Delta(\alpha_1)$, $\delta_{2\alpha_2}\in\Delta(\alpha_2)$, and 1 of type $b$, $\delta_{\alpha_3}\in\Delta(\alpha_3)$. By \ref{ssscartan}, one has $c(\delta_{2\alpha_i},-)={1\over 2}\langle\alpha_i^\vee,-\rangle$, $i=1,2$, and $c(\delta_{\alpha_3},-)=\langle\alpha_3^\vee,-\rangle$. So the Cartan pairing is as follows:
\[\begin{array}{r|rrr}c(-,-)&\sigma_1&\sigma_2&\sigma_3\\\hline \delta_{2\alpha_1}&2&-1&0\\\delta_{2\alpha_2}&-1&2&-1\\\delta_{\alpha_3}&0&-2&2\end{array}\]

In Example 4, there are 5 colors, 2 of type $a$, $\delta_{\alpha_1}^+,\delta_{\alpha_1}^-\in\Delta(\alpha_1)$, and 3 of type $b$, $\delta_{\alpha_i}\in\Delta(\alpha_i)$, $i=2,3,4$. For $\delta_{\alpha_i}$ of type $b$, $c(\delta_{\alpha_i},-)=\langle\alpha^\vee,-\rangle$. The full Cartan pairing is then as follows:
\[\begin{array}{r|rrr}c(-,-)&\sigma_1&\sigma_2&\sigma_3\\\hline \delta_{\alpha_1}^+&1&0&0\\\delta_{\alpha_1}^-&1&-1&0\\\delta_{\alpha_2}&-1&1&-1\\\delta_{\alpha_3}&0&0&1\\\delta_{\alpha_4}&0&-1&1\end{array}\]

The reader is invited to determine colors and Cartan pairings for the two other diagrams above. In Chapter~3, we will often leave this to the reader.

As another exercise, the reader should also determine, for every spherical systems $\mathcal S$ given by one of the diagrams above, all the diagrams corresponding to the different localizations of $\mathcal S$.

\vspace{6ex}\section{Type $\mathsf F_4$}\label{sstype}

In this section, unless otherwise stated, $G$ will denote a simple group of
type $\mathsf F_4$. This group is unique (up to isomorphism), since adjoint and
simply connected groups of type $\mathsf F_4$ coincide, and has rank 4 and
dimension 52.

The Dynkin diagram of $G$ is
\[\begin{picture}(5400,600)\multiput(0,300)(3600,0){2}{\usebox{\segm}}\put(1800,300){\usebox{\rightbisegm}}\end{picture}\]

and the 24 positive roots of G are
\begin{itemize}
\item[] \mbox{$\alpha_1$},\quad \mbox{$\alpha_2$},\quad \mbox{$\alpha_3$},\quad \mbox{$\alpha_4$},\quad \mbox{$\alpha_1+\alpha_2$},\quad \mbox{$\alpha_2+\alpha_3$},\quad \mbox{$\alpha_3+\alpha_4$},\quad \mbox{$\alpha_1+\alpha_2+\alpha_3$},\quad\\ \mbox{$\alpha_2+\alpha_3+\alpha_4$},\quad \mbox{$\alpha_1+\alpha_2+\alpha_3+\alpha_4$},
\item[] \mbox{$\alpha_2+2\alpha_3$},\quad \mbox{$\alpha_1+\alpha_2+2\alpha_3$},\quad \mbox{$\alpha_1+2\alpha_2+2\alpha_3$},
\item[] \mbox{$\alpha_2+2\alpha_3+\alpha_4$},\quad \mbox{$\alpha_2+2\alpha_3+2\alpha_4$},
\item[] \mbox{$\alpha_1+\alpha_2+2\alpha_3+\alpha_4$},\quad \mbox{$\alpha_1+2\alpha_2+2\alpha_3+\alpha_4$},\quad \mbox{$\alpha_1+2\alpha_2+3\alpha_3+\alpha_4$},\quad\\ \mbox{$\alpha_1+\alpha_2+2\alpha_3+2\alpha_4$},\quad \mbox{$\alpha_1+2\alpha_2+2\alpha_3+2\alpha_4$},\quad \mbox{$\alpha_1+2\alpha_2+3\alpha_3+2\alpha_4$},\quad\\ \mbox{$\alpha_1+2\alpha_2+4\alpha_3+2\alpha_4$},\quad \mbox{$\alpha_1+3\alpha_2+4\alpha_3+2\alpha_4$},\quad \mbox{$2\alpha_1+3\alpha_2+4\alpha_3+2\alpha_4$}.
\end{itemize}

Since $G$ is adjoint and simply connected, the fundamental weights belong
to the root lattice. They are
\begin{itemf}
\item[] $\omega_1=2\alpha_1+3\alpha_2+4\alpha_3+2\alpha_4$, 
\item[] $\omega_2=3\alpha_1+6\alpha_2+8\alpha_3+4\alpha_4$, 
\item[] $\omega_3=2\alpha_1+4\alpha_2+6\alpha_3+3\alpha_4$, 
\item[] $\omega_4=\alpha_1+2\alpha_2+3\alpha_3+2\alpha_4$.
\end{itemf}

The corresponding fundamental representations have dimensions
respectively 52, 1274, 273 and 26.

\subsection{Structure of the maximal parabolic subgroups}\label{ssstructure}

Let for a moment $G$ be again an arbitrary semisimple group.

Let $Q$ be a maximal parabolic subgroup of $G$ containing $B$. Let us write (as before) $\m Q$ for the opposite parabolic subgroup with respect to $T$, and $M=Q\cap\m Q$ for the common Levi subgroup of $Q$ and $\m Q$.

Then $\mathrm{Lie}(Q^u)$ is graded 
\[\mathrm{Lie}(Q^u)=\bigoplus_{i=1}^s\mathfrak n_i\]
by the action of $M^r$, which is $\cong\mathbb C^\times$. Moreover, the $\mathfrak n_i$'s are simple $M$-modules and one has $[\mathfrak n_i,\mathfrak n_j]\subset\mathfrak n_{i+j}$, so in particular 
\[\mathrm{Lie}((Q^u,Q^u))=\bigoplus_{i=2}^s\mathfrak n_i.\]
Furthermore,
\[\mathrm{Lie}(\m Q^u)=\bigoplus_{i=1}^s\mathfrak n_{-i},\]
where the $\mathfrak n_{-i}$ are isomorphic to the duals $\mathfrak n_i^\ast$ as $M$-modules.

The maximal parabolic subgroups of $G$ can be indexed by simple roots in the following way. To every $\alpha\in S$, one associates $Q_\alpha$, the unique maximal parabolic subgroup of $G$ containing $B$ such that $P_\alpha\not\subset Q_\alpha$. Then the unique color of $G/\m Q_\alpha$ is $D_\alpha$ (a color of type $b$).

Let us return now to the case when $G$ is of type $\mathsf F_4$. The general picture above particularizes then to the data given in the following table.

\begin{center}
\begin{tabular}{l|cccccc}
&semisimple type of $M$&$s$&$\dim \mathfrak n_1$&$\dim \mathfrak n_2$&$\dim \mathfrak n_3$&$\dim \mathfrak n_4$\\
\hline
$Q_{\alpha_1}$&$\mathsf C_3$&2&14&1\\
$Q_{\alpha_2}$&$\mathsf A_1\times\mathsf A_2$&3&12&6&2\\
$Q_{\alpha_3}$&$\mathsf A_2\times\mathsf A_1$&4&6&9&2&3\\
$Q_{\alpha_4}$&$\mathsf B_3$&2&8&7
\end{tabular}
\end{center}

\begin{remark}
When $G$ is of type $\mathsf F_4$, every parabolic subgroup $Q$ is conjugated to $\m Q$. But since we want to present the case of type $\mathsf F_4$ as example of the general case, we will not use this to simplify our notations.
\end{remark}

\subsection{Spherical roots of type $\mathsf F_4$}\label{ssssphericalroots}

Using the table in \ref{ssswonderfulvarietiesofrank1}, one obtains the list of all spherical roots of type $\mathsf F_4$:

\begin{itemize}
\item[] $\alpha_1$,\quad $\alpha_2$,\quad $\alpha_3$,\quad $\alpha_4$,\quad $2\alpha_1$,\quad $2\alpha_2$,\quad $2\alpha_3$,\quad $2\alpha_4$,
\item[] $\alpha_1+\alpha_3$,\quad $\alpha_1+\alpha_4$,\quad $\alpha_2+\alpha_4$,
\item[] $\alpha_1+\alpha_2$,\quad $\alpha_3+\alpha_4$,
\item[] $\alpha_2+\alpha_3$,\quad $2\alpha_2+2\alpha_3$,
\item[] $\alpha_1+\alpha_2+\alpha_3$,\quad $2\alpha_1+2\alpha_2+2\alpha_3$,\quad $\alpha_1+2\alpha_2+3\alpha_3$,
\item[] $\alpha_2+2\alpha_3+\alpha_4$,
\item[] $\alpha_1+2\alpha_2+3\alpha_3+2\alpha_4$.
\end{itemize}

Since compatible couples $(\sigma,S^p)$ are in bijective correspondence with spherical systems of rank 1, we refer to Tables~\ref{tr1a1}--\ref{tr1f4} below (where we list all these systems of type $\mathsf F_4$) for an explicit description of these couples.

\subsection{Spherical systems of type $\mathsf F_4$}\label{ssssphericalsystems}

In the tables below, we will give the complete list of the diagrams of all spherical systems $(\Sigma,S^p,\mathbf A)$ of type $\mathsf F_4$, ordered by their rank and by the support of $\Sigma$. Recall that the spherical systems of rank 0 correspond to the generalized flag varieties, and those of rank 1 to compatible couples.

The spherical systems of rank 4 are subdivided into two tables: Table~\ref{tr4ss} for spherical systems admitting a morphism to the spherical system of the full flag variety ($\emptyset$, $\emptyset$, $\emptyset$), called \textit{strongly solvable} (see \ref{sssquotient2}); Table~\ref{tr4o} for the remaining cases.

\input{tables.1}

%\clearpage
\chapter{Further notions concerning wonderful varieties}

In this chapter, we will go beyond the basic notions of Chapter~1. We will introduce the notion of wonderful subgroup, which allows a Lie theory point of view on wonderful varieties. Then, after mentioning briefly facts on equivariant automorphisms, we will introduce and study a natural notion of morphism between wonderful varieties. This will be our main tool for analyzing the examples of type $\mathsf F_4$ in Chapter~3. Finally, we will examine the relations which spherical orbits in simple projective spaces have with wonderful varieties.

In what follows, $X$ will always denote a wonderful $G$-variety, $\mathcal S=(\Sigma, S^p, \mathbf A)$ its spherical system, $\Delta$ its set of colors and $H$ will be an isotropy group of $G$ at a point in the open orbit of $X$. The Cartan pairing $c\colon \Delta\times\Sigma\to\mathbb Z$ will also be considered as a $\mathbb Z$-bilinear pairing $c\colon\mathbb Z\Delta\times\mathbb Z\Sigma\to\mathbb Z$.

We put $d(\mathcal S)=\mathrm{card}(\Delta)-\mathrm{card}(\Sigma)$, integer we will call the \textit{defect} of $\mathcal S$; this integer is also equal to the rank of $\Xi(H)$.

\vspace{6ex}\section{Wonderful subgroups}

We will call \textit{wonderful subgroups} of $G$ those subgroups which are isotropy groups at points in wonderful $G$-varieties.

If $H$ is a wonderful subgroup of $G$, a wonderful $G$-variety whose open orbit is isomorphic to $G/H$ is called a \textit{wonderful completion} of $G/H$. This wonderful completion exists (by definition) and is unique up to $G$-isomorphism.

Every wonderful subgroup $H$ of $G$ is a spherical subgroup of $G$ (i.e.\ $B$ has an open orbit in $G/H$) and $\mathrm{N}_G(H)/H$ is finite. The converse is not true in general; but every spherical subgroup equal to its normalizer is wonderful (this result is due to F.~Knop \cite{Kn96}).

\smallskip
Let us place here some general considerations. We have now three (equivalent) levels in our subject of study:

\begin{itemf}
\item[-] wonderful varieties, which is the level of highest geometric content; but these varieties, with the exception of some of low rank, can rarely be described and studied explicitly, since they appear in a natural way only as subvarieties of very high dimensional projective spaces;
\item[-] wonderful subgroups, which is the level of Lie theory; these subgroups, more precisely their conjugacy classes, are more accessible as we will see in the examples of Chapter~3;
\item[-] spherical systems, which is the combinatorial level; this is the most accessible level, since these invariants can be described most easily, as we have seen already in Chapter~1.
\end{itemf}

The difficulty is to go from one level to another. For instance, from a point $x$ on the open orbit of $G$ in a wonderful variety $X$ (where one is near $H=G_x$, a wonderful subgroup), the closed orbit of $G$ in $X$ looks to be very far, ``at infinity''. But the spherical system of $X$ can be read off most easily at the points $z$ and $z_\alpha^\pm$ ($\alpha\in\Sigma\cap S$), which are on or near the closed orbit (see Chapter~1). How to go (in general) from the spherical system to the wonderful subgroup (and back), is for the moment not completely understood.

When $H$ is a wonderful subgroup, by \textit{its} spherical system we will of course mean the spherical system of the wonderful completion of $G/H$. Conversely, when we talk about the wonderful subgroup of (or associated to) a spherical system $\mathcal S$, we will mean of course the generic isotropy group of the wonderful variety having $\mathcal S$ as spherical system. Similar remarks apply to wonderful subgroups and diagrams.

Here are some (first) relations between properties of $H$ and properties of $\mathcal S$:
\begin{itemf}
\item[-] $H$ is reductive if and only if there exists $\sigma\in\mathbb N\Sigma$ such that $c(\delta,\sigma)>0$ for all $\delta\in\Delta$. 
\item[-] $H$ is very reductive in $G$ if and only if $d(\mathcal S)=0$
\end{itemf}
($H$ \textit{very reductive} in $G$ means that $H$ is not contained in any proper parabolic subgroup of $G$; very reductive implies reductive, and even semisimple if $H$ is connected).

\bigskip
\begin{example}\label{exa2}
Consider the following diagrams:
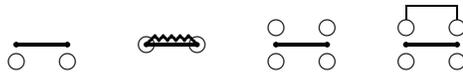
\[\begin{picture}(2400,1800)\put(300,900){\usebox{\segm}}\multiput(0,0)(1800,0){2}{\usebox{\aprime}}\end{picture}\quad\quad
\begin{picture}(2400,1800)\put(0,600){\usebox{\atwo}}\end{picture}\quad\quad
\begin{picture}(2400,1800)\put(300,900){\usebox{\segm}}\multiput(0,0)(1800,0){2}{\usebox{\aone}}\end{picture}\quad\quad
\begin{picture}(2400,2250)\put(300,900){\usebox{\segm}}\multiput(0,0)(1800,0){2}{\usebox{\aone}}\multiput(300,1800)(1800,0){2}{\line(0,1){450}}\put(300,2250){\line(1,0){1800}}\end{picture}\]
Here $G=\mathrm{SL}(3)$, and corresponding wonderful subgroups are (from left to right):
\begin{itemf}
\item[-] $H=\mathrm{SO}(3)\cdot C_G$ (where $C_G$ denotes the center of $G$); this $H$ is very reductive in $G$;
\item[-] $H=\mathrm{GL(2)}$, which is reductive but not very reductive in $G$;
\item[-] $H=T\,U_{\alpha}$, where $\alpha$ is any root and $U_{\alpha}$ denotes the corresponding unipotent subgroup of dimension 1; these $H's$ are conjugated in $G$; 
\item[-] $H=L\,H^u$, where $L$ is a subtorus of dimension 1 of $T$ and $H^u$ is a unipotent subgroup of dimension 2 of $G$ such that $L=\mathrm N_T(H^u)$; these $H$'s are conjugated in $G$.
\end{itemf} 
\end{example}

\vspace{6ex}\section{Equivariant automorphisms}\label{sssequivariantautomorphisms}

Remember that $X$ denotes a wonderful $G$-variety, $\mathcal S = (\Sigma, S^p , \mathbf A)$ its spherical system and $\Delta$ its set of colors. Let us choose a point $x$ in the open orbit of $G$ in $X$, and put $H = G_x$.

We will say that a spherical root $\sigma\in\Sigma$ is \textit{loose} in $\mathcal S$, if
\begin{itemf}
\item[-] either $\sigma=\alpha\in S\cap\Sigma$ and $c(\delta^+_\alpha,\sigma')=c(\delta^-_\alpha,\sigma')$, for all $\sigma'\in\Sigma$;
\item[-] or $\sigma\in\Sigma\setminus S$, $2\sigma\in\Sigma(G)$ and the couple $(2\sigma, S^p)$ is compatible (in the sense of spherical systems, see \ref{ssswonderfulvarietiesofrank1}).
\end{itemf}

We denote by $\Sigma_\ell(\mathcal S)$ the set of spherical roots that are loose in $\mathcal S$.

Let us denote by $\Gamma=\Gamma_X=\mathrm{Aut}_G(X)$ the group of $G$-automorphisms of $X$. 
Restriction to $G/H=G.x\subset X$ induces an isomorphism between $\Gamma$ and $\mathrm{Aut}_G(G/H)=(\mathrm{N}_G(H)/H)^{\mathrm{opp}}$.

For every $\sigma\in\Sigma_\ell(\mathcal S)$, there exists a unique element $\gamma(\sigma)$ of order 2 in $\Gamma$ that fixes the points of the divisor $D^\sigma$ (this follows from the conjecture of \ref{question}; a direct proof has been given by I.~Losev in \cite{Lo07}). 
Moreover, these $\gamma(\sigma)$, $\sigma\in\Sigma_\ell(\mathcal S)$, commute and generate $\Gamma$. For every $\sigma\in\Sigma_\ell(\mathcal S)$, the variety $X/\sigma$ is again wonderful, and its spherical roots are $\{2\sigma\}\cup(\Sigma\setminus\{\sigma\})$.

As we will see in the following section, the presence of nontrivial equivariant automorphisms introduces some complications. But it is clear by the analysis above, that for a majority of wonderful varieties $X$, $\Gamma_X$ is reduced to the identity (i.e.\ $\Sigma_\ell(\mathcal S)=\emptyset$).

Here are some simple examples of type $\mathsf B_3$ with $\Sigma_\ell(\mathcal S)\neq\emptyset$ (we will give $\mathcal S$ by its diagram):

\[\begin{picture}(4200,1800)\put(300,900){\usebox{\segm}}\put(2100,900){\usebox{\rightbisegm}}\put(300,900){\circle{600}}\put(1800,0){\usebox{\aone}}\put(2500,1300){\usebox{\toe}}\put(3600,0){\usebox{\aprime}}\end{picture}\qquad\begin{picture}(7200,1800)\put(0,600){$\Sigma_\ell=\{\alpha_2\}$}\end{picture}\]

\[\begin{picture}(4200,1800)\put(300,900){\usebox{\segm}}\put(2100,900){\usebox{\rightbisegm}}\put(0,0){\usebox{\aone}}\put(1800,600){\usebox{\GreyCircle}}\end{picture}\qquad\begin{picture}(7200,1800)\put(0,600){$\Sigma_\ell=\{\alpha_2+\alpha_3\}$}\end{picture}\]

\[\begin{picture}(4200,1800)\put(300,900){\usebox{\segm}}\put(2100,900){\usebox{\rightbisegm}}\put(0,600){\usebox{\atwo}}\put(1800,600){\usebox{\GreyCircle}}\put(3600,0){\usebox{\aone}}\put(3100,1300){\usebox{\tow}}\end{picture}\qquad\begin{picture}(7200,1800)\put(0,600){$\Sigma_\ell=\{\alpha_3\}$}\end{picture}\]

\vspace{6ex}\section{Wonderful morphisms}

1) A $G$-morphism $\phi\colon X\to \p X$ will be called \textit{wonderful} if $X$ and $\p X$ are wonderful $G$-varieties, and if $\phi$ is dominant (i.e.\ surjective) and has connected fibers.

2) An inclusion $H\subset \p H$ of wonderful subgroups of $G$ is called \textit{co-connected} if $\p H/H$ is connected.

3) Let $\mathcal S=(\Sigma, S^p, \mathbf A)$ be a spherical system and $\Delta$ its set of colors. A subset $\Delta^\ast$ of $\Delta$ is called \textit{distinguished} in $\Delta$, if there exists $\delta\in\mathbb N_{>0}\Delta^\ast$ such that $c(\delta,\sigma)\geq0$, for all $\sigma\in\Sigma$. 

For every distinguished subset $\Delta^\ast$ of $\Delta$, the \textit{quotient system}\\ $\mathcal S/\Delta^\ast=(\Sigma/\Delta^\ast, S^p/\Delta^\ast, \mathbf A/\Delta^\ast)$ is defined as follows:
\begin{itemf}
\item[-] $\Sigma/\Delta^\ast$ is the set of minimal generators of the (free) semigroup\\ $\{\sigma\in\mathbb N\Sigma,\, c(\delta,\sigma)=0\mathrm{\ for\ all\ }\delta\in\Delta^\ast\}$;
\item[-] $S^p/\Delta^\ast=S^p\cup\{\alpha\in S,\, \Delta(\alpha)\subset\Delta^\ast\}$;
\item[-] $\mathbf A/\Delta^\ast$ is the union of the $\mathbf A(\alpha)$'s, $\alpha\in S\cap(\Sigma/\Delta^\ast)$, and the Cartan pairing $\mathbf A/\Delta^\ast\times\Sigma/\Delta^\ast\to\mathbb Z$ is obtained from $\mathbf A\times\Sigma\to\mathbb Z$ in the obvious way.
\end{itemf}

In this section we will study these three notions and their interrelations. 

\subsection{Quotient systems and wonderful subgroups}\label{sssquotient1}

Let $H\subset \p H$ be two wonderful subgroups of $G$. Denote by $\phi\colon G/H\to G/\p H$ the natural map, and $\Delta_\phi$ the set of $D\in\Delta=\Delta_{G/H}$ such that $\phi(D)$ is dense in $G/\p H$. Then
\begin{itemf}
\item[-] $\Delta_\phi$ is distinguished in $\Delta$; 
\item[-] conversely, for every distinguished subset $\Delta^\ast$ of $\Delta$, there exists a unique wonderful subgroup $\p H$ having the following properties: $H\subset \p H$, $\p H/H$ is connected, and $\Delta_\phi=\Delta^\ast$;
\item[-] moreover, if $\p H/H$ is connected, the spherical system of $\p H$ is given by the quotient system $\mathcal S/\Delta_\phi$ (where $\mathcal S$ is the spherical system of $H$)
\end{itemf}
(these results are close to statements of F.~Knop in \cite{Kn91}).

\subsection{Quotient systems and wonderful morphisms}\label{sssquotient2}

Remember that $X$ denotes a wonderful $G$-variety, $\mathcal S=(\Sigma,S^p,\mathbf A)$ its spherical system and $\Delta$ its set of colors. The group $\Gamma_X$ acts naturally in the set of distinguished subsets of $\Delta$. Let $\p X$ be another wonderful variety. The group $\Gamma_X\times\Gamma_{\p X}$ acts naturally in the set of wonderful morphisms $\phi\colon X\to \p X$. If $\phi$ is such a morphism, denote by $\Delta_\phi$ the set of $D\in\Delta$ such that $\phi(D)=\p X$. Then
\begin{itemf}
\item[-] $\Delta_\phi$ is distinguished in $\Delta$; 
\item[-] conversely, for every orbit $\Gamma_X.\Delta^\ast$ of distinguished subsets in $\Delta$, there exist wonderful morphisms $\phi\colon X\to \p X$, defined (and unique) up to isomorphism, such that $\Delta_\phi\in\Gamma_X.\Delta^\ast$;
\item[-] moreover, the spherical system of $\p X$ is given by the quotient system $\mathcal S/\Delta_\phi$.
\end{itemf}

This statement follows from \ref{sssquotient1} and from the fact that if $x\in X$ and $\p x\in \p X$ are such that $H=G_x\subset \p H=G_{\p x}$, then $\phi\colon G/H=G.x\to G/\p H=G.\p x$ extends always to a $G$-morphism $\phi\colon X\to \p X$ such that $\phi(x)=\p x$.

Notice that there may exist morphisms $\psi\colon X\to \p X$ which are not obtained by $\psi=\p \gamma\circ\phi\circ\gamma$, $\gamma\in\Gamma_X$ and $\p \gamma\in\Gamma_{\p X}$ (they correspond to the existence of subgroups $K$ of $G$ containing $H$ and conjugated to $\p H$, but such that the two inclusions $H\subset \p H$ and $H\subset K$ cannot be conjugated simultaneously).

\begin{remark} 
Most of the wonderful varieties we will see in the examples in Chapter~3 have trivial automorphism groups; in that case the statements above have simpler form.
\end{remark}

Here are some relations between properties on different levels:
\begin{itemf}
\item[-] 
the variety $\p X$ is homogeneous (in other words $\p H$ is a parabolic subgroup of $G$), if and only if $\Sigma/\Delta^\ast=\emptyset$; $\p H$ is then conjugated in $G$ to $\m Q$, where $Q$ is the parabolic subgroup of $G$ containing $B$ and corresponding to $S^p/\Delta^\ast$; and there exists a $\m Q$-variety $Y$ such that $X\cong G\ast_{\m Q}Y$;
\item[-]
moreover, $\m Q^r$ acts trivially on $Y$ (so that $X$ is obtained by parabolic induction from the wonderful $\m Q/\m Q^r$-variety $Y$), if and only if $\mathrm{supp}(\Sigma)\subset S^p/\Delta^\ast$.
\item[-]
in particular, $X$ cannot be obtained by nontrivial parabolic induction if and only if $\mathrm{supp}(\Sigma)\cup S^p=S$.
\item[-]
the group $H$ is \textit{strongly solvable} in $G$ (i.e.\ $H$ is contained in a Borel subgroup of $G$) if and only if there exists a distinguished subset $\Delta^\ast$ such that $\mathcal S/\Delta^\ast=(\emptyset,\emptyset,\emptyset)$ (this system being the spherical system of $G/B$).
\end{itemf}

\subsection{Quotient systems and diagrams}\label{qsdiagrams}

Let us compute explicitly some quotient systems. Consider the spherical system $\mathcal S=(\Sigma,S^p,\mathbf A)$ of $G=\mathrm{SL}(4)$ having diagram:
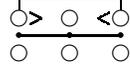
\[\begin{picture}(4200,2250)\multiput(300,900)(1800,0){2}{\usebox{\segm}}\multiput(0,0)(1800,0){3}{\usebox{\aone}}\multiput(300,1800)(3600,0){2}{\line(0,1){450}}\put(300,2250){\line(1,0){3600}}\put(700,1300){\usebox{\toe}}\put(3100,1300){\usebox{\tow}}\end{picture}\]
Here $\Sigma=S$, $S^p=\emptyset$, $\mathrm{card}(\mathbf A)=5$ and the Cartan pairing is:
\[\begin{array}{r|rrr}c(-,-)&\alpha_1&\alpha_2&\alpha_3\\\hline\delta_{\alpha_1}^+&1&-1&1\\\delta_{\alpha_1}^-&1&0&-1\\\delta_{\alpha_2}^+&0&1&0\\\delta_{\alpha_2}^-&-1&1&-1\\\delta_{\alpha_3}^-&-1&0&1\end{array}\]
The following subsets of colors are distinguished in $\Delta$:
\[\Delta^{1}=\{\delta_{\alpha_1}^+,\delta_{\alpha_2}^-\},\quad\Delta^{2}=\{\delta_{\alpha_1}^-,\delta_{\alpha_3}^-\},\quad\Delta^{3}=\{\delta_{\alpha_2}^+\},\]\[\Delta^{1,2}=\Delta^{1}\cup\Delta^{2},\quad\Delta^{2,3}=\Delta^{2}\cup\Delta^{3}.\]
Indeed, $c(\delta_{\alpha_1}^++\delta_{\alpha_2}^-,\sigma)=0$, for all $\sigma\in\Sigma$, so $\Delta^{1}$ is distinguished. Similarly, $c(\delta_{\alpha_1}^-+\delta_{\alpha_3}^-,\sigma)=0$, $\sigma\in\Sigma$. Since $\Delta^{1}$ and $\Delta^{2}$ are distinguished, so is their union $\Delta^{1,2}$. Again, $c(\delta_{\alpha_2}^+,\sigma)\geq0$, $\sigma\in\Sigma$, so $\Delta^{3}$ and $\Delta^{2,3}$ are distinguished. 

Let us now compute the quotients. Set $\sigma=m_1\alpha_1+m_2\alpha_2+m_3\alpha_3\in\mathbb N\Sigma$. Let us start with $\Delta^{1}$: $c(\delta_{\alpha_1}^+,\sigma)=0$ is equivalent to $m_2=m_1+m_3$ and $c(\delta_{\alpha_2}^-,\sigma)=0$ gives a proportional equation, so $\sigma=m_1(\alpha_1+\alpha_2)+m_3(\alpha_2+\alpha_3)$ and $\Sigma/\Delta^{1}=\{\alpha_1+\alpha_2,\alpha_2+\alpha_3\}$. For $\Delta^{2}$ one has $m_1=m_3$, so $\Sigma/\Delta^{2}=\{\alpha_1+\alpha_3,\alpha_2\}$. For $\Delta^{3}$ one has $m_2=0$, so $\Sigma/\Delta^{3}=\{\alpha_1,\alpha_3\}$. Similarly for $\Delta^{1,2}$ and $\Delta^{2,3}$. One gets the quotient spherical systems corresponding to the following diagrams: 
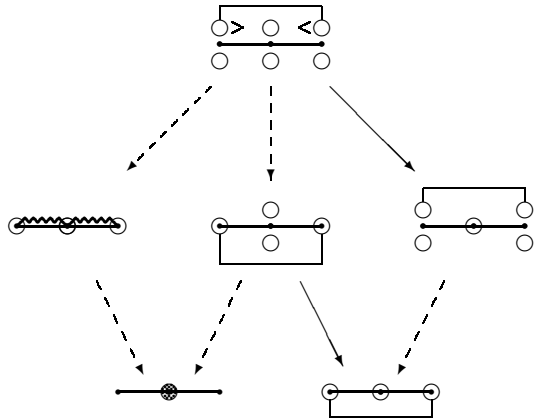
\[\begin{picture}(18600,14550)
\put(7200,12300){
\multiput(300,900)(1800,0){2}{\usebox{\segm}}
\multiput(0,0)(1800,0){3}{\usebox{\aone}}
\multiput(300,1800)(3600,0){2}{\line(0,1){450}}
\put(300,2250){\line(1,0){3600}}
\put(700,1300){\usebox{\toe}}
\put(3100,1300){\usebox{\tow}}
}
\put(14400,5850){
\multiput(300,900)(1800,0){2}{\usebox{\segm}}
\multiput(0,0)(3600,0){2}{\usebox{\aone}}
\multiput(300,1800)(3600,0){2}{\line(0,1){450}}
\put(300,2250){\line(1,0){3600}}
\put(2100,900){\circle{600}}
}
\put(7200,5850){
\multiput(300,900)(1800,0){2}{\usebox{\segm}}
\multiput(300,900)(3600,0){2}{\circle{600}}
\multiput(300,-450)(3600,0){2}{\line(0,1){1050}}
\put(300,-450){\line(1,0){3600}}
\put(1800,0){\usebox{\aone}}
}
\put(0,5850){
\multiput(0,600)(1800,0){2}{\usebox{\atwo}}
}
\put(11100,0){
\multiput(300,900)(1800,0){2}{\usebox{\segm}}
\multiput(300,900)(1800,0){3}{\circle{600}}
\multiput(300,0)(3600,0){2}{\line(0,1){600}}
\put(300,0){\line(1,0){3600}}
}
\put(3600,0){
\multiput(300,900)(1800,0){2}{\usebox{\segm}}
\put(1800,600){\usebox{\GreyCircle}}
}
\put(11400,11700){\vector(1,-1){3000}}
\multiput(9300,11700)(0,-600){5}{\line(0,-1){300}}
\put(9300,8700){\vector(0,-1){300}}
\put(7200,11700){
\multiput(0,0)(-400,-400){7}{\multiput(0,0)(-20,-20){10}{\line(-1,0){30}}}
\put(-2800,-2800){\vector(-1,-1){200}}
}
\put(3150,4800){
\multiput(0,0)(300,-600){5}{\multiput(0,0)(15,-30){10}{\line(-1,0){30}}}
\put(1500,-3000){\vector(1,-2){150}}
}
\multiput(8250,4800)(7200,0){2}{
\multiput(0,0)(-300,-600){5}{\multiput(0,0)(-15,-30){10}{\line(-1,0){30}}}
\put(-1500,-3000){\vector(-1,-2){150}}
}
\put(10350,4800){\vector(1,-2){1500}}
\end{picture}\]

An arrow between two diagrams denotes that the target arises as a quotient of the source. For ``minimal'' quotients we sometimes use a dashed arrow to give more information about the ``type'' of the quotient, as it will be explained in \ref{minimal}.

In Chapter~3, quotients of spherical systems by distinguished subsets of colors will be often represented only by arrows between diagrams. Indeed, the distinguished subset of colors can usually be recovered easily, given only the two spherical systems.

To explain this, let $\mathcal S$ be a spherical system with set of colors $\Delta$, and $\Delta^\ast$ be a distinguished subset of $\Delta$. Denote by $\p \mathcal S=\mathcal S/\Delta^\ast$ the quotient system, and by $\p \Delta$ its set of colors. Then $\p \Delta$ can be identified with $\Delta\setminus\Delta^\ast$, and for every $\alpha\in S$, one has $\p \Delta(\alpha)\subset\Delta(\alpha)$. So if one knows $\Delta$ and $\p \Delta$, one can usually recover $\Delta^\ast$ (the only problem is to understand what happens for $\alpha\in S\cap\Sigma$).

Consider for instance the diagrams
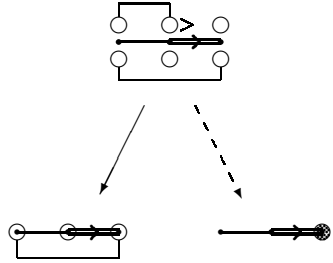
\[\begin{picture}(11400,9000)
%\multiput(0,0)(11400,0){2}{\line(0,1){8400}}\multiput(0,0)(0,8400){2}{\line(1,0){11400}}
\put(3600,6750){
\put(300,900){\usebox{\segm}}
\put(2100,900){\usebox{\rightbisegm}}
\multiput(0,0)(1800,0){3}{\usebox{\aone}}
\multiput(300,1800)(1800,0){2}{\line(0,1){450}}
\put(300,2250){\line(1,0){1800}}
\multiput(300,0)(3600,0){2}{\line(0,-1){450}}
\put(300,-450){\line(1,0){3600}}
\put(2500,1300){\usebox{\toe}}
}
\put(7200,0){
\put(300,900){\usebox{\segm}}
\put(2100,900){\usebox{\rightbisegm}}
\put(3600,600){\usebox{\GreyCircle}}
}
\put(0,0){
\put(300,900){\usebox{\segm}}
\put(2100,900){\usebox{\rightbisegm}}
\multiput(300,900)(1800,0){3}{\circle{600}}
\multiput(300,0)(3600,0){2}{\line(0,1){600}}
\put(300,0){\line(1,0){3600}}
}
\put(4800,5400){\vector(-1,-2){1575}}
\put(3600,300){\multiput(3000,5100)(300,-600){5}{\multiput(0,0)(15,-30){10}{\line(0,-1){30}}}
\put(4500,2100){\vector(1,-2){150}}}
\end{picture}\]

These are diagrams of type $\mathsf B_3$. Denote by $\mathcal S$ the spherical system of the diagram on top, and by ${}^{(1)}\!\mathcal S$ and ${}^{(2)}\!\mathcal S$ the spherical system of those on the bottom line. The Cartan pairing of $\mathcal S$ is 
\[\begin{array}{r|rrr}c(-,-)&\alpha_1&\alpha_2&\alpha_3\\\hline\delta_{\alpha_1}^+&1&1&-1\\\delta_{\alpha_1}^-&1&-2&1\\\delta_{\alpha_2}^-&-2&1&0\\\delta_{\alpha_3}^+&-1&0&1\end{array}\]
Then the two distinguished subsets of colors $\Delta^{(1)}$ and $\Delta^{(2)}$ of $\Delta$ such that $\mathcal S/\Delta^{(i)}={}^{(i)}\!\mathcal S$, $i=1,2$, are clearly given by $\Delta^{(1)}=\{\delta_{\alpha_1}^+,\delta_{\alpha_3}^+\}$ and $\Delta^{(2)}=\Delta\setminus\{\delta_{\alpha_3}^+\}$ (the reader is invited to check this). %These two subsets are actually the only two ``minimal'' distinguished subsets of $\Delta$.

As another exercise, the reader is invited to explicit in the same way the quotients given by the following figure:
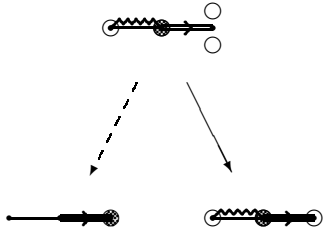
\[\begin{picture}(11400,8400)
%\multiput(0,0)(11400,0){2}{\line(0,1){8400}}\multiput(0,0)(0,8400){2}{\line(1,0){11400}}
\put(7200,-600){
\put(300,900){\usebox{\segm}}
\put(2100,900){\usebox{\rightbisegm}}
\put(0,600){\usebox{\atwo}}
\put(1800,600){\usebox{\GreyCircle}}
\put(3900,900){\circle{600}}
}
\put(0,-600){
\put(300,900){\usebox{\segm}}
\put(2100,900){\usebox{\rightbisegm}}
\put(3600,600){\usebox{\GreyCircle}}
}
\put(3600,6150){
\put(300,900){\usebox{\segm}}
\put(2100,900){\usebox{\rightbisegm}}
\put(0,600){\usebox{\atwo}}
\put(1800,600){\usebox{\GreyCircle}}
\put(3600,0){\usebox{\aone}}
}
\put(-3600,0){\multiput(8400,5100)(-300,-600){5}{\multiput(0,0)(-15,-30){10}{\line(0,-1){30}}}
\put(6900,2100){\vector(-1,-2){150}}}
\put(6600,5100){\vector(1,-2){1575}}
\end{picture}\]

This way of describing quotients of spherical systems is not without ambiguity, since it may happen that there are several distinguished subsets of colors giving the same quotient system. For instance the arrow
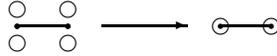
\[\begin{picture}(9000,1800)\multiput(300,900)(7200,0){2}{\usebox{\segm}}\multiput(0,0)(1800,0){2}{\usebox{\aone}}\multiput(7500,900)(1800,0){2}{\circle{600}}\put(3300,900){\vector(1,0){3000}}\end{picture}\]
may come from $\{\delta_{\alpha_1}^+,\delta_{\alpha_2}^+\}$, $\{\delta_{\alpha_1}^+,\delta_{\alpha_2}^-\}$ or $\{\delta_{\alpha_1}^-,\delta_{\alpha_2}^+\}$.

\subsection{Generic fibers of wonderful morphisms}

Let $\phi\colon X\to \p X$ be a wonderful morphism, and let $\p x$ be a point in the open orbit of $G$ in $\p X$. The fiber $Y=\phi^{-1}(\p x)$ is called a \textit{generic fiber} of $\phi$. It follows from the Theorem of Sard that $Y$ is a complete and smooth $\p H$-variety, where we have put $\p H=G_{\p x}$ (this group is not necessarily reductive nor connected). Let us denote by $\Sigma_\phi$ the set of spherical roots of $X$ which contribute to spherical roots of $\p X$ (recall that the second are sums of the first). For every $\Sigma'\subset \Sigma_X$, remember that $X_{\Sigma'}$ denotes the wonderful subvariety of $X$ having $\Sigma'$ as set of simple roots.

\begin{proposition}
\begin{itemf}
\item[]
\item[1)] One has $\phi(X_{\Sigma'})=\p X$ if and only if $\Sigma'\supset\Sigma_\phi$.
\item[2)] The $\p H$-variety $Y$ is wonderful, its rank is $\mathrm{card}(\Sigma\setminus\Sigma_\phi)$, and $Y$ is spherical under the action of any Levi subgroup of $\p H$. 
\end{itemf}
\end{proposition}

\begin{sketche}\makebox{}\\
1) A $G$-morphism $\psi\colon Z\to \p Z$ between spherical $G$-varieties is dominant if and only if $\Xi(\psi)\colon\Xi(\p Z)\to\Xi(Z)$ is injective. From this, (1) follows, since $\Sigma_{\p X}$ is a basis of $\Xi(\p X)$.

2) The intersection of $Y$ with each orbit of $G$ in $X$ is either empty, or is an orbit of $\p H$ in $Y$. This implies that, for every $\sigma\in\Sigma\setminus\Sigma_\phi$, $Y\cap D^\sigma$ is a divisor stable by $\p H$. One checks easily that these divisors have all the properties of the definition of wonderful varieties. The last assertion (although less easy) is left to the reader.
\end{sketche}

\subsection{Notions of minimality}\label{minimal}

Let $H\subset \p H$ be a co-connected inclusion of wonderful subgroups of $G$. We will say that this inclusion is \textit{minimal} if there exists no proper intermediate wonderful subgroup $K$ with $K/H$ connected. 

\begin{proposition}
If $H\subset \p H$ is such a minimal inclusion, then three possibilities can occur:
\begin{itemf}
\item[($\mathcal P$)] either $H^u$ contains $\p H^u$ strictly; then $H$ is a maximal parabolic subgroup of $\p H$;
\item[($\mathcal R$)] or $H^u=\p H^u$; then $H/H^r$ is very reductive and maximal in $\p H/\p H^r$;
\item[($\mathcal L$)] or $H^u$ is strictly contained in $\p H^u$; then $\mathrm{Lie}(\p H^u)/\mathrm{Lie}(H^u)$ is a simple $H$-module, $H=\mathrm{N}_{\p H}(H^u)$, and the Levi factors of $H$ and $\p H$ differ only by their connected centers.
\end{itemf}
\end{proposition}
\vspace{-4pt}
\begin{sketche}
Let $H=L\,H^u$ and $\p H=\p L\,\p H^u$ be Levi decompositions such that $L\subset \p L$. Denote by $q\colon H\to \p H/\p H^r$ the natural map. If $q$ is not surjective, it follows by minimality that $\p H^u$ is contained in $H^u$. Then either $q(H)$ is contained in a parabolic subgroup of $\p H/\p H^r$ and we have $\mathcal P$; or $q(H)$ is very reductive in $\p H/\p H^r$, which implies $H^u=\p H^u$, and we have $\mathcal R$. If $q$ is surjective, then $L$ and $\p L$ differ only by their connected centers, and we have $\mathcal L$.
\end{sketche}

Let $H$ be a wonderful subgroup of $G$. Assume that $Q$ is a parabolic subgroup of $G$, containing $B$, such that $\m Q$ is minimal among the parabolic subgroups of $G$ containing $H$. Then there exists a sequence of co-connected inclusions of wonderful subgroups
\vspace{-4pt}
\[H=H_0\subset H_1\subset\ldots\subset H_{m-1}\subset H_m\subset \m Q\]
\vspace{-4pt}
having the following properties:
\begin{itemf}
\item[-] $H_{i-1}\subset H_i$ is minimal of type $\mathcal L$ ($i=1,\ldots,m$);
\item[-] $H_m^r=\m Q^r$ and $H_m/H_m^r$ is very reductive in $\m Q/\m Q^r$
\end{itemf}
\vspace{-4pt}
(this is a reformulation of statements in \cite{Lo07}).

\bigskip
A wonderful morphism $\phi\colon X\to \p X$ will be called \textit{minimal} if it cannot be written as composition of two wonderful morphisms, both of which are not isomorphisms. 

Let $\p x$ be a point in the open orbit of $G$ in $\p X$, and consider the generic fiber $Y=\phi^{-1}(\p x)$ introduced in the preceding section. According to the analysis above, minimal wonderful morphisms are divided into the following three types: 
\begin{itemf}
\item[($\mathcal P$)] $G_{\p x}$ acts transitively on $Y$;
\item[($\mathcal R$)] $G_{\p x}$ does not act transitively on $Y$, but $G_{\p x}^r$ acts trivially on $Y$;
\item[($\mathcal L$)] $G_{\p x}^r$ has an open dense orbit in $Y$.
\end{itemf}

Let $\mathcal S=(\Sigma,S^p,\mathbf A)$ be the spherical system of $X$, and $\Delta$ its set of colors. It is clear that $\phi$ is minimal if and only if $\Delta_\phi$ is minimal among the distinguished subsets of $\Delta$.

Let us call a color \textit{negative} if its values on $\Sigma$ (under the Cartan pairing) are $\leq0$. All negative colors are of the form $\delta_\alpha$, for $\alpha\in S^b$ (which is uniquely determined by the color). A negative color $\delta_{\alpha}$ will be called \textit{interior} if $\alpha\in\mathrm{supp}(\Sigma)$, and \textit{exterior} if $\alpha\not\in\mathrm{supp}(\Sigma)$.

For every $\alpha\in S$, let us denote by $Q_\alpha$ a maximal parabolic associated to $\alpha$ (containing B). The following characterization of negative exterior colors follows from the properties of parabolic inductions: $D_{\alpha}$ is a negative exterior color for a wonderful variety $X$ if and only if there exists a generic isotropy group $H$ of $X$ verifying: $\m Q_\alpha^u\subset H\subset \m Q_\alpha$.

Is the following statement true in general? Let $D_\alpha$ be a color of type $b$ for a wonderful variety $X$. Then $D_\alpha$ is negative interior if and only if the derived subgroup $(\m Q_\alpha^u,\m Q_\alpha^u)$ is not trivial, and there exists a generic isotropy group $H$ of $X$, verifying: $(\m Q_\alpha^u,\m Q_\alpha^u)\subset H\subset \m Q_\alpha$ but $\m Q_\alpha^u\not\subset H$.

Let $X$ be a wonderful variety with spherical system $\mathcal S=(\Sigma,S^p,\mathbf A)$. Let $\phi\colon X\to \p X$ be a minimal wonderful morphism, and put $\Delta^\ast=\Delta_\phi$. Remember that the spherical system of $\p X$ is given by the quotient system $\mathcal S/\Delta^\ast=(\Sigma/\Delta^\ast,S^p/\Delta^\ast,\mathbf A/\Delta^\ast)$. How can one see the different types of $\phi$ on the combinatorial data?

From the preceding section follows that $\phi$ is of type $\mathcal P$ if and only if $\Sigma_\phi=\Sigma$.

The following facts are easy to check:
\begin{itemf}
\item[-] $d(\mathcal S)>d(\mathcal S/\Delta^\ast)$ if and only if we are in type $\mathcal P$; 
\item[-] if $d(\mathcal S)<d(\mathcal S/\Delta^\ast)$, then we are in type $\mathcal L$;
\item[-] if $d(\mathcal S)=d(\mathcal S/\Delta^\ast)$ and a new exterior negative color appears in $\mathcal S/\Delta^\ast$, then we are in type $\mathcal L$.
\end{itemf}

The following statement is true in all examples we know: if $d(\mathcal S)=d(\mathcal S/\Delta^\ast)$ and no new negative color appears in $\mathcal S/\Delta^\ast$, then we are in type $\mathcal R$.

If $d(\mathcal S)=d(\mathcal S/\Delta^\ast)$ and the only new negative color which appears in $\mathcal S/\Delta^\ast$ is interior, we can either be in type $\mathcal L$ or in type $\mathcal R$ (the first happens very often, see for instance Example~2 below, or Section~\ref{ss1521}; for the second see Section~\ref{notf4}, Example~7).

\bigskip
\textit{Examples}

1) Let $\mathcal S=(\Sigma,S^p,\mathbf A)$ be any spherical system and $\Delta$ its set of colors. A color $\delta\in\mathbf A$ is called \textit{projective}, if $c(\delta,\sigma)\geq0$, for all $\sigma\in\Sigma$. If $\delta$ is any projective color, the set $\{\delta\}$ is clearly distinguished in $\Delta$, and the corresponding quotient is of type $\mathcal L$.

2) 
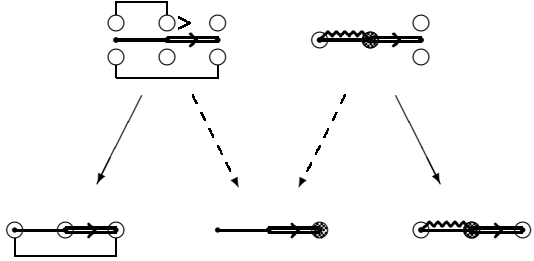
\[\begin{picture}(18600,9000)
\put(3600,6750){
\put(300,900){\usebox{\segm}}
\put(2100,900){\usebox{\rightbisegm}}
\multiput(0,0)(1800,0){3}{\usebox{\aone}}
\multiput(300,1800)(1800,0){2}{\line(0,1){450}}
\put(300,2250){\line(1,0){1800}}
\multiput(300,0)(3600,0){2}{\line(0,-1){450}}
\put(300,-450){\line(1,0){3600}}
\put(2500,1300){\usebox{\toe}}
}
\put(7200,0){
\put(300,900){\usebox{\segm}}
\put(2100,900){\usebox{\rightbisegm}}
\put(3600,600){\usebox{\GreyCircle}}
}
\put(0,0){
\put(300,900){\usebox{\segm}}
\put(2100,900){\usebox{\rightbisegm}}
\multiput(300,900)(1800,0){3}{\circle{600}}
\multiput(300,0)(3600,0){2}{\line(0,1){600}}
\put(300,0){\line(1,0){3600}}
}
\put(4800,5700){\vector(-1,-2){1575}}
\put(3600,600){\multiput(3000,5100)(300,-600){5}{\multiput(0,0)(15,-30){10}{\line(0,-1){30}}}
\put(4500,2100){\vector(1,-2){150}}}
\put(7200,600){
\put(7200,-600){
\put(300,900){\usebox{\segm}}
\put(2100,900){\usebox{\rightbisegm}}
\put(0,600){\usebox{\atwo}}
\put(1800,600){\usebox{\GreyCircle}}
\put(3900,900){\circle{600}}
}
\put(3600,6150){
\put(300,900){\usebox{\segm}}
\put(2100,900){\usebox{\rightbisegm}}
\put(0,600){\usebox{\atwo}}
\put(1800,600){\usebox{\GreyCircle}}
\put(3600,0){\usebox{\aone}}
}
\put(-3600,0){\multiput(8400,5100)(-300,-600){5}{\multiput(0,0)(-15,-30){10}{\line(0,-1){30}}}
\put(6900,2100){\vector(-1,-2){150}}}
\put(6600,5100){\vector(1,-2){1575}}
}
\end{picture}\] 

These diagrams have already been analyzed in \ref{qsdiagrams}. Each arrow represents a minimal quotient. As a general rule we put dashed arrows to represent minimal quotients with strictly decreasing defect, and continuous arrows otherwise. Indeed, the two spherical systems with diagram on the top line have defect 1, both having $\mathrm{card}(\Delta)=4$ and $\mathrm{card}(\Sigma)=3$. For those on the second line, the defect is respectively $2-1=1$, $1-1=0$, $3-2=1$.

3) Here is an example of type $\mathcal R$:
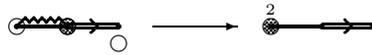
\[\begin{picture}(13200,1800)\multiput(300,900)(9000,0){2}{\usebox{\segm}}\multiput(2100,900)(9000,0){2}{\usebox{\rightbisegm}}\put(0,600){\usebox{\atwo}}\put(1800,600){\usebox{\GreyCircle}}\put(3600,0){\usebox{\aprime}}\put(9000,600){\usebox{\GreyCircleTwo}}\put(5100,900){\vector(1,0){3000}}\end{picture}\]

\vspace{6ex}\section{Spherical closure}\label{sssphericalclosure}

We will call \textit{simple projective $G$-spaces} the $G$-varieties of the form $\mathbb P(V)$, where $V$ is any simple (rational) $G$-module. In this section we will see that spherical orbits in simple projective spaces play an important role in the theory of wonderful varieties.

\subsection{Definition and first properties}\label{sphericalclosure1}

Let $H$ be any spherical subgroup of $G$. The group $\mathrm{N}_G(H)$ acts naturally on $\Delta_{G/H}$, the set of colors of $G/H$. We call \textit{spherical closure} of $H$ in $G$ (denoted by $\overline H^{sph}$) the kernel of the action of $\mathrm{N}_G(H)$ on $\Delta_{G/H}$; if $H=\overline H^{sph}$, we say that $H$ is \textit{spherically closed} in $G$.

Spherical subgroups equal to their normalizer are of course spherically closed. Spherically closed subgroups are wonderful (result due to F.~Knop \cite{Kn96}).

Let $X$ be a wonderful $G$-variety, $\mathcal S=(\Sigma,S^p,\mathbf A)$ its spherical system, $\Delta$ its set of colors, and let $H$ be a generic isotropy group of $G$ in $X$. Here are some properties relating these notions:
\begin{itemf}
\item[-] $H$ is spherically closed if and only if $\Sigma_\ell(\mathcal S)\subset S$; 
\item[-] in particular, $H=\mathrm{N}_G(H)$ if and only if $\Sigma_\ell(\mathcal S)=\emptyset$;
\item[-] the variety $X$ is strict if and only if $S\cap\Sigma=\emptyset=\Sigma_\ell(\mathcal S)$ (\textit{strict} means that $\mathrm{N}_G(G_x)=G_x$, for all $x\in X$);
\item[-] from the combinatorial characterization above follows in particular that, if a generic isotropy group of X is spherically closed, then so are all the
other isotropy groups;
\item[-] if $H$ is spherically closed, then $\mathrm{N}_G(H)/H=\Gamma_X$ can be identified to the group $\Gamma=\Gamma_{\mathcal S}$ of permutations of $\Delta$ stabilizing each $\Delta(\alpha)$, $\alpha\in S$, and leaving invariant the Cartan pairing (i.e.\ by those permutations which may exchange the 2 colors of $\Delta(\alpha)$, when $\alpha\in\Sigma_\ell(\mathcal S)$, but fix all other colors).
\end{itemf}

The wonderful subgroups not spherically closed are somewhat exceptional. 
In type $\mathsf F_4$, only three spherical systems have wonderful subgroups which are not spherically closed:
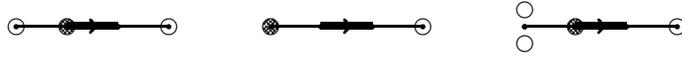
\[\begin{picture}(24000,1800)
\put(0,0){
\put(300,900){\usebox{\dynkinf}}
\multiput(300,900)(5400,0){2}{\circle{600}}
\put(1800,600){\usebox{\GreyCircle}}
}
\put(9000,0){
\put(300,900){\usebox{\dynkinf}}
\put(5700,900){\circle{600}}
\put(0,600){\usebox{\GreyCircle}}
}
\put(18000,0){
\put(300,900){\usebox{\dynkinf}}
\put(5700,900){\circle{600}}
\put(1800,600){\usebox{\GreyCircle}}
\put(0,0){\usebox{\aone}}
}
\end{picture}\]

\subsection{Further properties of the spherical closure}\label{sphericalclosure2}

Let $H$ be a spherical subgroup of $G$. It is well known that $H$ fixes at most a finite number of points in every simple projective $G$-space $\mathbb P(V)$ (otherwise one would be able to construct nonconstant $B$-invariant rational functions on $G/H$, which is impossible, since $B$ has an open orbit in $G/H$).

Let us denote by $q_H \colon G \to G/H$ the canonical map. The following lemma gives alternative characterizations of the spherical closure.

\begin{lemma} 
Let $g\in G$. The following three conditions are equivalent:
\begin{itemf}
\item[(1)] $g\in\overline H^{sph}$;
\item[(2)] $q_H^{-1}(D).g=q_H^{-1}(D)$, for every $D\in\Delta_{G/H}$;
\item[(3)] $g$ fixes all points in simple projective $G$-spaces that are fixed by $H$.
\end{itemf}
\end{lemma}

\begin{proof}
$(1)\Rightarrow (2)$ is clear.

$(2)\Rightarrow (3)$ Let $V$ be a simple $G$-module, and let $[x]\in\mathbb P(V)$ be a point fixed by $H$ (where $x$ denotes a point in $V\setminus\{0\}$ ``over'' $[x]$); such an $x$ is an eigenvector of $H$. Choose an eigenvector $v$ of $B$ in $V^\ast$. Then $v\otimes x$ gives, by means of the map $V^\ast\otimes V\to \mathbb C[G]$ corresponding to the natural map $G\to \mathrm{End}(V)=V\otimes V^\ast$, a function $f\in\mathbb C[G]$ which is an eigenvector of $B$ (acting by left translations), and an eigenvector of $H$ (acting by right translations). Since we assume $G$ simply connected, $\mathbb C[G]$ is factorial, so $f$ can be written in a unique way as $\prod_{D\in\Delta}f_D^{n(f,D)}$, where the $f_D$'s are equations of the $q_H^{-1}(D)$'s, $D\in\Delta_{G/H}$. The assumption (2) implies then that the $f_D$'s (and so also $f$) are again eigenvectors of $g$, which gives $g.[x]=[x]$.

$(3)\Rightarrow (1)$ One has to show that $g\in\mathrm{N}_G(H)$. Since $H$ is an algebraic group, there exist $G$-modules $U$ (which are in general not simple) and $H$-eigenvectors $u\in U$, such that $H=G_{[u]}$. We can assume that $U$ has a direct sum decomposition into simple $G$-modules $U=V_1\oplus\ldots\oplus V_n$ such that $u$ can be written $u=v_1+\ldots+v_n$, where the $v_i$ are $H$-eigenvectors in $V_i$ ($i=1,\ldots,n$). Since by (3) these $v_i$ are also eigenvectors of $g$, one can find $a\in\mathrm{Aut}_G(U)$ such that $g.[u]=a.[u]$. From this follows that $G_{g.[u]}=G_{[u]}=H$, which implies $g\in\mathrm{N}_G(H)$.
\end{proof}

\begin{corollarynumber}
If $H\subset \p H$ are spherical subgroups of $G$, then $\overline H^{sph}\subset \overline {\p H}^{sph}$.
\end{corollarynumber}

\begin{proof} 
This follows for instance immediately from the property (3) of the preceding lemma.
\end{proof}

\begin{corollarynumber}\label{corollarysc}
A spherical subgroup $H$ of $G$ is spherically closed if and only if it occurs as isotropy group in some simple projective $G$-space.
\end{corollarynumber}

\begin{proof}
Let $V$ be a simple $G$-module and let $x\in V\setminus\{0\}$ be such that $H=G_{[x]}$ is spherical in $G$. Then the property (3) of the preceding lemma implies $\overline H^{sph}\subset G_{[x]}=H$, so $H$ is spherically closed.

Conversely, let $H$ be any spherically closed subgroup of $G$. Choose pairwise different integers $n(D)>0$, $D\in\Delta_{G/H}$ and put $f=\prod_{D\in\Delta}f_D^{n(D)}$. Then there exists a simple $G$-module $V$ and $[x]\in\mathbb P(V)^H$, such that $f$ is obtained from $x$ as in the proof $(2)\Rightarrow (3)$ above. We know already that $\p H=G_{[x]}$ is spherically closed in $G$. By definition, $H\subset \p H$. Because of the choice of $f$, the natural map $G/H\to G/\p H$ induces a bijection $\Delta_{G/H}\to\Delta_{G/\p H}$. So the sets $\{q_H^{-1}(D), D\in\Delta_{G/H}\}$ and $\{q_{\p H}^{-1}(D), D\in\Delta_{G/\p H}\}$ are the same, which implies $H=\p H$.
\end{proof}

The following result will be used several times in Chapter~3.

\begin{corollarynumber}\label{argument} 
Let $K$ be a spherically closed subgroup of $G$, and let $H$ be a subgroup of $K$ such that $\mathrm{N}_K(H) = H$. Then $H$ is wonderful in $G$ if (and only if) it is spherical in $G$.
\end{corollarynumber}

\begin{proof}
If $H$ is spherical in $G$, Corollary~\ref{corollarysc} implies $\overline H^{sph}\subset\overline K^{sph} = K$. So we have $\overline H^{sph}\subset\mathrm{N}_G(H)\cap K=\mathrm{N}_K(H)$ which is equal to $H$ by assumption. This implies $H$ spherically closed in $G$, so $H$ is wonderful in $G$.
\end{proof}

\subsection{Spherical orbits in simple projective spaces}\label{ssorbits}

Let $V$ be any simple (rational) $G$-module. It is known that only finitely many orbits of $G$ in $\mathbb P(V)$ are spherical. In what follows, we will explain how one can classify these with the help of spherical systems.

For any $\delta\in\mathbb N\Delta$, we write $n(\delta,D)$ for the coefficient of $\delta$ at $D\in\Delta$, and define $\mathrm{supp}_\Delta(\delta)$ as the set of colors $D$ such that $n(\delta,D)>0$.

We will say that a couple ($\mathcal S$,$\delta$) -- where $\mathcal S$ is a spherical system and $\delta\in\mathbb N\Delta$ -- is \textit{faithful}, if the following conditions are fulfilled:
\begin{itemf}
\item[] (1)\hspace{6pt}$\mathcal S$ is spherically closed;
\item[] (2)\hspace{6pt}any (nonempty) distinguished subset of $\Delta$ meets $\mathrm{supp}_\Delta(\delta)$;
\item[] (3)\hspace{6pt}$n(\delta,D^+_\alpha)\neq n(\delta,D^-_\alpha)$, for all $\alpha\in\Sigma_\ell(\mathcal S)$.
\end{itemf}

Remember the natural map $\omega\colon\mathbb N\Delta\to\mathbb N\Omega$ (where $\mathbb N\Omega$ is the set of dominant weights) introduced in \ref{ssscartan}. On the combinatorial level, $\omega$ is given as follows:
\begin{itemf}
\item[] -\hspace{6pt}if $\delta\in\Delta^a\cup\Delta^b$, then $\omega(\delta)$ is the sum of those fundamental weights $\omega_{\alpha}$\\  
\rule{10pt}{0pt}such that $\delta\in\Delta(\alpha)$;
\item[] -\hspace{6pt}if $\delta\in\Delta^{2a}$, then $\omega(\delta) = 2\omega_\alpha$ where $\delta\in\Delta(\alpha)$.
\end{itemf}

Remember also that we have introduced above in \ref{sphericalclosure1} a permutation group $\Gamma=\Gamma_{\mathcal S}$ of $\Delta$, which is acting also in $\mathbb N\Delta$. It is clear that for every faithful couple ($\mathcal S$,$\delta$) and every $\gamma\in\Gamma$, the couple ($\mathcal S$,$\gamma.\delta$) is again faithful, and that $\omega(\gamma.\delta)=\omega(\delta)$.

If a simple $G$-module $V$ has a dominant weight $\pi$, we denote by $\pi^\ast$ the dominant weight of its dual module $V^\ast$.

\begin{proposition}
Assume the conjecture of \ref{question} true for $G$. Let $V_\pi$ be a simple $G$-module of highest weight $\pi$. Then there exists a natural bijection between the set of spherical orbits of $G$ in $\mathbb P(V_\pi)$ and the set of $\Gamma$-orbits of faithful couples \mbox{($\mathcal S$, $\Gamma.\delta$)} such that $\omega(\delta)=\pi^\ast$.
\end{proposition}

\begin{proof}
Let ($\mathcal S$, $\delta$) be a faithful couple such that $\omega(\delta)=\pi^\ast$. Let $H$ be a spherically closed subgroup of $G$ having $\mathcal S$ as spherical system. Denote by $f_D\in\mathbb C[G]$ equations of the $q^{-1}(D)$'s, $D\in\Delta_{G/H}$, and put $f=\prod_{D\in\Delta} f_D^{n(\delta,D)}$, where $\delta=\sum_{D\in\Delta} n(\delta,D)D$. Let $V$ be a simple $G$-module and $x$ be an eigenvector of $H$ in $V$, such that $f$ is obtained from $x$ as in the previous subsection. The character $\omega(\delta)$ is also the $B$-character of $f$, which is the dominant weight of $V^\ast$. So $V$ has dominant wight $\pi$ and we have associated to the couple ($\mathcal S$, $\delta$) a spherical orbit $G.[x]$ in $\mathbb P(V_\pi)$.

Conversely, let $G.[x]$ be a spherical orbit in $\mathbb P(V_\pi)$. By Corollary~\ref{corollarysc} of \ref{sphericalclosure2}, $H=G_{[x]}$ is a spherically closed subgroup of $G$. Denote by $\mathcal S$ its spherical system, and by $f=\prod_{D\in\Delta}f_D^{n(f,D)}$ the function on $G$ given by $x$ as in the previous section. Put $\delta=\sum_{D\in\Delta}n(f,D)D$. By construction, the couple ($\mathcal S$,$\delta$) is faithful (otherwise the natural map $G/H\to G.[x]$ would not be injective), and one has $\omega(\delta)=\pi^\ast$.

If ($\mathcal S',\delta'$) is another couple going to the same (spherical) orbit in $\mathbb P(V_\pi)$, then $\mathcal S'=\mathcal S$ (since this system does depend only on the orbit). But $\delta$ does depend on the choice of the point on the orbit fixed by $H$, points which are exchanged by $\mathrm{N}_G(H)$, so $\delta'\in\Gamma.\delta$.
\end{proof}

\begin{example}
Let $G = \mathrm{SL}(4)$ and let $\pi = \omega_1+\omega_3$ be the dominant weight of the adjoint representation $\mathfrak{sl}(4)$. The following table contains all faithful couples \mbox{($\mathcal S$, $\delta$)} such that $\omega(\delta) = \omega_1+\omega_3$:

\[\begin{array}{c|ccccc}
\begin{picture}(1800,2700)\put(300,1050){$\mathcal S$}\end{picture} &
\begin{picture}(6000,2700)\put(900,0){\multiput(300,1350)(1800,0){2}{\usebox{\segm}}\multiput(300,1350)(3600,0){2}{\circle{600}}}\end{picture} &
\begin{picture}(6000,2700)\put(900,0){\multiput(300,1350)(1800,0){2}{\usebox{\segm}}\multiput(300,1350)(1800,0){3}{\circle{600}}\multiput(300,1050)(3600,0){2}{\line(0,-1){600}}\put(300,450){\line(1,0){3600}}}\end{picture} &
\begin{picture}(6000,2700)\put(900,0){\multiput(300,1350)(1800,0){2}{\usebox{\segm}}\multiput(0,450)(1800,0){3}{\usebox{\aone}}\multiput(300,2250)(3600,0){2}{\line(0,1){450}}\put(300,2700){\line(1,0){3600}}\put(2500,1750){\usebox{\toe}}}\end{picture} &
\begin{picture}(6000,2700)\put(900,0){\multiput(300,1350)(1800,0){2}{\usebox{\segm}}\multiput(300,1350)(3600,0){2}{\circle{600}}\put(0,1050){\multiput(300,300)(25,25){13}{\circle*{70}}\put(600,600){\multiput(0,0)(300,0){10}{\multiput(0,0)(25,-25){7}{\circle*{70}}}\multiput(150,-150)(300,0){10}{\multiput(0,0)(25,25){7}{\circle*{70}}}}\multiput(3900,300)(-25,25){13}{\circle*{70}}}}\end{picture} &
\begin{picture}(6000,2700)\put(900,0){\multiput(300,1350)(1800,0){2}{\usebox{\segm}}\multiput(300,1350)(3600,0){2}{\circle{600}}\multiput(300,1050)(3600,0){2}{\line(0,-1){1050}}\put(300,0){\line(1,0){3600}}\put(1800,450){\usebox{\aprime}}}\end{picture} \\
%\hline
\rule{0pt}{3ex}\delta & \delta_{\alpha_1}+\delta_{\alpha_3} & \delta_{\alpha_1} & \delta_{\alpha_1}^+ & \delta_{\alpha_1}+\delta_{\alpha_3} & \delta_{\alpha_1} 
\end{array}\]
Notice that the group $\Gamma$ is trivial for all systems in this table. The first three couples correspond to nilpotent orbits, the last two couples to semisimple orbits. Notice that the third couple corresponds to the nilpotent matrices having Jordan form 
\[\left(\begin{array}{rrrr}0&1&0&0\\0&0&1&0\\0&0&0&0\\0&0&0&0\end{array}\right)\]
orbit which is spherical in $\mathbb P(\mathfrak{sl}(4))$, but not in $\mathfrak{sl}(4)$.
\end{example}

\bigskip
\textit{Remarks}

1) Spherical orbits in general adjoint representations are well known (see \cite{Pa94,Pa03,CCC05,Co08}).

2) Let us use again the notation of the above proposition. G.~Pezzini has shown in \cite{Pe07} that the closure of $G.[x]$ in $\mathbb P(V_\pi)$ is wonderful if and only if $\mathcal S$ is strict and $\mathrm{supp}_\Delta(\delta)=\Delta$.

\bigskip
As a nontrivial exercise, the reader is invited to characterize the inclusion relations of spherical orbit closures in simple projective spaces, in terms of faithful couples.

\clearpage
\chapter{Examples of type $\mathsf F_4$}

We already have listed in Chapter~1 the 266 spherical systems of type $\mathsf F_4$. We will now give more information on the corresponding wonderful varieties. The non-trivial ones have dimension between 15 (the smallest generalized flag variety) and 28 (the dimension of the Borel subgroup), and are not very accessible to study. One knows that they can be theoretically realized as closures of orbits in some high dimensional projective spaces, but even the set-theoretical description of these closures seems in general to be out of reach. So we have to be content mainly with describing their isotropy groups.

In some of the following sections, we start with a spherical system $\mathcal S$ and obtain the structure of the generic isotropy group $H$ of the corresponding wonderful variety $X$. Then we list all spherical systems $\mathcal S_i$ that have $\mathcal S$ as minimal quotient, and determine the corresponding generic isotropy groups $H_i$ as subgroups of $H$. Even if in this way we do not get a direct picture of the wonderful varieties themselves, we obtain some understanding of the wonderful subgroups of type $\mathsf F_4$, and of their co-connected minimal inclusions. These can be seen as explicit examples of minimal wonderful morphisms.

We also give an explicit (but abstract) classification of the spherical orbits in the projective fundamental representations of type $\mathsf F_4$. 

Finally, in the last section, we discuss some examples not of type $\mathsf F_4$, for the additional information they contain.

%\clearpage
\subsection*{Contents of Chapter 3}
\begin{itemf}
\item[\ref{ss1521}] A remarkable example of rank 2
\item[\ref{ss811}] A remarkable example of rank 1
\item[\ref{projectivecsss}] Projective colors and strongly solvable systems
\item[\ref{sb3}] Examples coming from type $\mathsf B_3$
\item[\ref{ss711}] Examples coming from type $\mathsf C_3$
\item[\ref{ss1514}] Another remarkable example of rank 2
\item[\ref{constant}] $\mathcal L$-type minimal morphisms with constant defect
\item[\ref{increasing}] $\mathcal L$-type minimal morphisms with strictly increasing defect
\item[\ref{fibered}] An example of a fiber product
\item[\ref{ssfundamental}] Spherical orbits in fundamental representations
\item[\ref{notf4}] Examples not of type $\mathsf F_4$
\end{itemf}

%\clearpage
\vspace{6ex}\section{A remarkable example of rank 2}\label{ss1521}

\[\begin{picture}(6000,1200)
\put(300,900){\usebox{\dynkinf}}
\multiput(300,900)(5400,0){2}{\circle{600}}
\multiput(300,0)(5400,0){2}{\line(0,1){600}}
\put(300,0){\line(1,0){5400}}
\put(3900,900){\circle{600}}
\put(1800,600){\usebox{\GreyCircle}}
\end{picture}\]
The corresponding spherical system $\mathcal S=(\Sigma,S^p,\mathbf A)$ has $\Sigma=\{\alpha_1+\alpha_4,\alpha_2+\alpha_3\}$, $S^p=\emptyset$ and the following Cartan pairing:
\[\begin{array}{r|rr}c(-,-)&\sigma_1&\sigma_2\\\hline\delta_{\alpha_1}&2&-1\\\delta_{\alpha_2}&-1&1\\\delta_{\alpha_3}&-1&0\end{array}\]

Let us briefly explain how one can prove that this $\mathcal S$ comes from a
wonderful variety $X$, unique up to isomorphism.

If $X$ exists, since $S^p = \emptyset$, it should have dimension $\mathrm{card}(\Sigma) + \dim(B^u) = 2 + 24 = 26$; in particular, a generic stabilizer $H$ should have again dimension $26 = 52 - 26 = \dim(G) - \dim(X)$. The spherical system $\mathcal S$ has a (unique) distinguished subset $\{\delta_{\alpha_1},\delta_{\alpha_2}\}$. Indeed, $c(\delta_{\alpha_1}+\delta_{\alpha_2},\sigma)\geq0$, for both $\sigma\in\Sigma$. Then $\Sigma/\{\delta_{\alpha_1},\delta_{\alpha_2}\}=\emptyset$ and the quotient $\mathcal S /\{\delta_{\alpha_1},\delta_{\alpha_2}\}$ has the following diagram:
\[\begin{picture}(5400,600)
\put(0,300){\usebox{\dynkinf}}
\put(3600,300){\circle{600}}
\end{picture}\]
So this quotient is homogeneous, and the corresponding wonderful variety is isomorphic to $G/^-Q$, where $Q=Q_{\alpha_3}$ (with notation as in \ref{ssstructure}); this $Q$ is of semisimple type $\mathsf A_2\times \mathsf A_1$. Moreover, the corresponding minimal wonderful morphism $X\to G/^-Q$ would be of type $\mathcal L$. Indeed, the defect remains constant and the negative color $\delta_{\alpha_3}$, which is interior in $\mathcal S$, becomes exterior in the quotient. Since the dimension of $Q$ is 32, it follows that if we choose $H$ included in $^-Q$, then we can suppose that it has the same Levi part as $^-Q$, say $M$, and $H^u\subset {}^-Q^u$ has codimension 6. Since $\mathrm{Lie}(^-Q^u)/\mathrm{Lie}((^-Q^u, {}^-Q^u))$ is a simple $M$-module of dimension 6 (see Section~\ref{ssstructure}), it follows that $H$ is necessarily of the form $H = (^-Q^u, {}^-Q^u)\,M$. So $H$ is determined up to conjugation in $G$, which gives the uniqueness of $X$ (up to $G$-isomorphism).

Conversely, let $H$ be a subgroup of $G$ given by $H = (^-Q^u,{}^-Q^u)\,M$ as above. One can check that $H$ is spherical (by using for instance Corollary~1.4 in \cite{Pa94}). Since $\mathrm{N}_{\m Q}(H)=H$, $H$ is wonderful in $G$ by Corollary~\ref{argument} of \ref{sphericalclosure2}. The wonderful completion $X'$ of $G/H$ has dimension 26 and comes with a (minimal) wonderful morphism $X'\to G/^-Q$. This implies that the spherical system of $X'$ is equal to $\mathcal S$ (there is simply no other system having the corresponding properties on the combinatorial level), which proves the existence of $X$.\hspace{\stretch{1}}$\square$

\bigskip
Henceforth, in Chapter 3, we will admit that all spherical systems considered are geometrically realizable.

With a little combinatorial effort, one obtains all spherical systems having $\mathcal S$ as quotient:

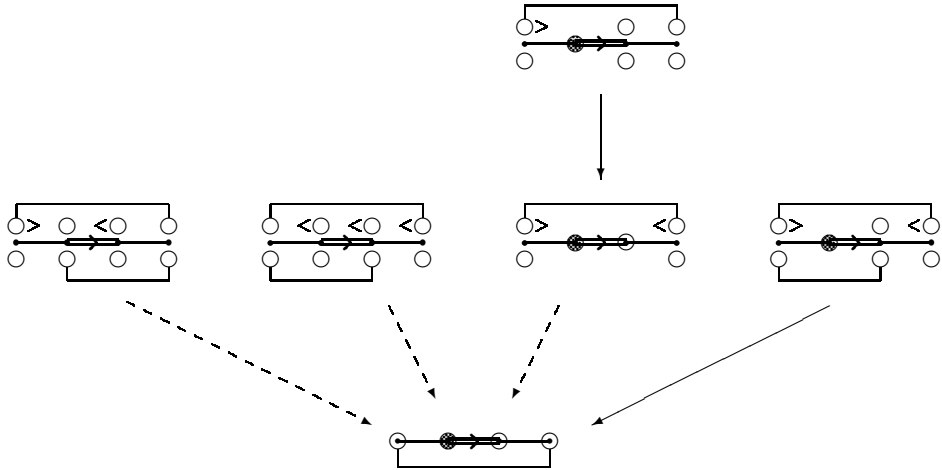
\[\strtwoffourtwoone\]

The reader is expected to explicit on his own the corresponding
spherical roots, Cartan pairings, distinguished subset of
colors\ldots\ for all the diagrams in the figure above, and to check that these (minimal) quotients are well defined. In what follows we will explain only part of this information.

Let $\mathcal S_i$, $i=1,2,3,4$, be the spherical systems corresponding to the diagrams on the second line of the figure (numbered from left to right), and by $H_i$ the generic isotropy groups of the corresponding wonderful varieties (we admit their existence). One can choose the $H_i$'s contained in $H$ (the generic isotropy group of the wonderful variety $X$ introduced above).

The group $H$ is connected of semisimple type $\mathsf A_2\times\mathsf A_1$. One checks easily that $H$ has (up to conjugation) exactly three parabolic subgroups of dimension $\geq24=\dim(B^u)$ (two having semisimple type $\mathsf A_1\times\mathsf A_1$, and the last having semisimple type $\mathsf A_2$). Since $d(\mathcal S_i)=2$ ($i=1,2,3$) and $d(\mathcal S)=1$, the quotients $\mathcal S_i\to\mathcal S$ ($i=1,2,3$) are of type $\mathcal P$. It follows that the $H_i$ ($i=1,2,3$) are the parabolic subgroups of $H$ mentioned above. We will come back to these systems $\mathcal S_i$ ($i=1,2,3$) in Section~\ref{constant}.

Let us now consider the spherical system $\mathcal S_4$. Here $d(\mathcal S_4) = d(\mathcal S) = 1$, and the only new negative color $\delta_{\alpha_3}$ which appears in $\mathcal S$ is interior. So a priori we don't know if the quotient  $\mathcal S_4\to\mathcal S$ is of type $\mathcal L$ or of type $\mathcal R$. But we have seen above that $H = (\m Q^u,\m Q^u)M$, where $Q = Q_{\alpha_3} = M\,Q^u$, and one can check that $H_4$ has codimension 2 in $H$. Since $(M,M)$ does not contain any semisimple subgroup of codimension 2, the quotient must be of type $\mathcal L$. We know that $\mathrm{Lie}(\m Q^u)$ decomposes into simple $M$-modules
\[\mathrm{Lie}(\m Q^u) = \mathfrak n_{-1} +\mathfrak n_{-2} +\mathfrak n_{-3} +\mathfrak n_{-4}\] 
having dimensions 6, 9, 2, 3 respectively (cf.\ \ref{ssstructure}). Since $H_4$ has codimension 2 in $H$, it follows that $H_4 = H_4^uM$, where $\mathrm{Lie}(H_4^u) = \mathfrak n_{-2} + \mathfrak n_{-4}$.

\bigskip
There is another minimal distinguished subset in $\Delta_4$ (the set of colors of $\mathcal S_4$), namely $\Delta_4\setminus\{\delta_{\alpha_4}^-\}$. This subset defines a minimal quotient (of type $\mathcal P$)
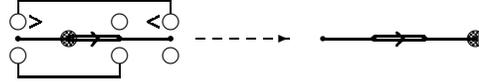
\[\begin{picture}(16800,2700)
\put(0,0){
\put(300,1350){\usebox{\dynkinf}}
\multiput(0,450)(3600,0){2}{\usebox{\aone}}
\put(5400,450){\usebox{\aone}}
\multiput(300,2250)(5400,0){2}{\line(0,1){450}}
\put(300,2700){\line(1,0){5400}}
\multiput(300,0)(3600,0){2}{\line(0,1){450}}
\put(300,0){\line(1,0){3600}}
\put(700,1750){\usebox{\toe}}
\put(4900,1750){\usebox{\tow}}
\put(1800,1050){\usebox{\GreyCircle}}
}
\multiput(6600,1350)(600,0){5}{\line(1,0){300}}
\put(9600,1350){\vector(1,0){300}}
\put(11100,1350){\usebox{\dynkinf}}
\put(16200,1050){\usebox{\GreyCircle}}
\end{picture}\]
which allows us a good transition to the next section.

%\clearpage
\vspace{6ex}\section{A remarkable example of rank 1}\label{ss811}
\[\begin{picture}(5700,600)
\put(0,300){\usebox{\dynkinf}}
\put(5100,0){\usebox{\GreyCircle}}
\end{picture}\]
The spherical system of the diagram above is $\mathcal S=(\Sigma,S^p,\mathbf A)$ where $\Sigma=\{\sigma\}$ with $\sigma=\alpha_1+2\alpha_2+3\alpha_3+2\alpha_4$, and $S^p=\mathbf A=\emptyset$.

This system has only one color (of type $b$), $\delta_{\alpha_4}$, and $c(\delta_{\alpha_4},\sigma)=1$. 

The corresponding wonderful variety is of rank 1 and has dimension
16. Since the defect of $\mathcal S$ is 0, the corresponding generic
isotropy group $H$ is very reductive, and of dimension 36. There is
only one such subgroup in $G$ (up to conjugation), which is simple,
simply connected and of type $\mathsf B_4$ (so $H$ is isomorphic to $\mathrm{Spin}(9)$), a well known symmetric subgroup of $G$. Since $\mathrm{N}_G(H)=H$, $H$ is spherically closed in $G$.

Among the 16 parabolic subgroups of $H$, only 5 have dimension $\geq
24 = \dim(B^u)$, the 4 maximal parabolic subgroups of $H$ (which have
semisimple types respectively $\mathsf B_3$, $\mathsf A_1\times
\mathsf B_2$, $\mathsf A_2\times\mathsf A_1$, $\mathsf A_3$) and
another one (which has semisimple type $\mathsf B_2$).

With a little (combinatorial) effort, one obtains all spherical systems having as quotient the system $\mathcal S$ (arrows correspond to minimal quotients):
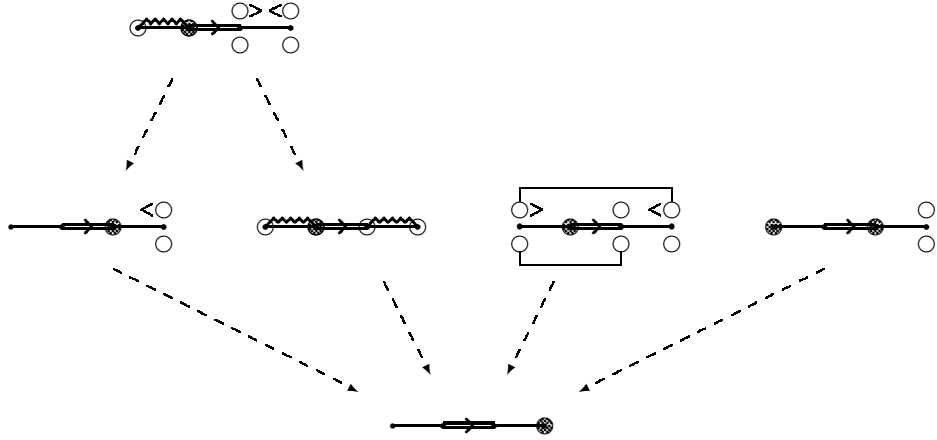
\[\stroneffouroneone\]
The reader is invited to explicit the corresponding spherical roots, Cartan pairings, distinguished subset of colors\ldots\ of all the diagrams above, and to check that all these quotients are well defined and minimal.

Let us denote by $\mathcal S_{1,2}$ the system corresponding to the
diagram on the first line of the figure, and by $\mathcal S_i$
($i=1,2,3,4$) the systems of those on the second line. Let us denote
by $H_1$, $H_2$, $H_3$, $H_4$, $H_{1,2}$ generic isotropy groups of
the corresponding wonderful varieties (which we assume to exist). They
have dimensions respectively 29, 25, 24, 26, 24. Since $d(\mathcal
S_{1,2}) = 2$, $d(\mathcal S_i) =1$ ($i=1,2,3,4$) and $d(\mathcal S )=
0$, all arrows above are of type $\mathcal P$. Since there are exactly
5 systems having $\mathcal S$ as quotient of type $\mathcal P$, we obtain that $H_1$, $H_2$, $H_3$, $H_4$, $H_{1,2}$ are exactly the 5 parabolic subgroups of $H$ we
have considered above (having semisimple types respectively $\mathsf B_3$, $\mathsf A_1\times \mathsf B_2$, $\mathsf A_2\times\mathsf A_1$, $\mathsf A_3$ and $\mathsf B_2$).

We will come back to these 5 systems in Section~\ref{constant}.

For the convenience of the reader, let us give some details in the case of the system $\mathcal S_{1,2}=(\Sigma_{1,2},S^p_{1,2},\mathbf A_{1,2})$. In this case $\Sigma_{1,2}=\{\alpha_1+\alpha_2,\alpha_2+\alpha_3,\alpha_3,\alpha_4\}$, $S^p_{1,2}=\emptyset$ and the Cartan pairing is as follows:
\[\begin{array}{r|rrrr}c(-,-)&\sigma_1&\sigma_2&\sigma_3&\sigma_4\\\hline\delta_{\alpha_1}&1&-1&0&0\\\delta_{\alpha_2}&1&1&-1&0\\\delta_{\alpha_3}^+&0&0&1&-1\\\delta_{\alpha_3}^-&-2&0&1&0\\\delta_{\alpha_4}^+&0&0&-1&1\\\delta_{\alpha_4}^-&0&-1&0&1\end{array}\]

Since the diagram of $\mathcal S_{1,2}$ might be ambiguous, let us remark that $c(\delta_{\alpha_4}^+,\alpha_3)=-1$ but $c(\delta_{\alpha_4}^+,\alpha_2+\alpha_3)=0$.

The subset of colors $\Delta_{1,2}'=\{\delta_{\alpha_1},\delta_{\alpha_2},\delta_{\alpha_3}^-\}$ is minimal distinguished. 
%Indeed, $c(\delta_{\alpha_1}+\delta_{\alpha_2}+\delta_{\alpha_3}^-,\sigma)=0$, for any $\sigma\in\Sigma_{1,2}$. 
Moreover, for any $\sigma\in\mathbb N\Sigma_{1,2}$, say $\sigma=\sum_{i=1}^4m_i\sigma_i$, $c(\delta,\sigma)=0$ for all $\delta\in\Delta_{1,2}'$ implies $2m_1=2m_2=m_3$. Therefore, one gets $\sigma=m_3(\alpha_1+2\alpha_2+3\alpha_3)+m_4(\alpha_4)$ and obtains as quotient the spherical system $\mathcal S_1$, with defect 1. 

Another minimal distinguished subset of colors is $\Delta_{1,2}''=\{\delta_{\alpha_3}^+,\delta_{\alpha_4}^+\}$. 
%Indeed, $c(\delta_{\alpha_3}^++\delta_{\alpha_4}^+,\sigma)=0$, for any $\sigma\in\Sigma_{1,2}$. 
Moreover, for any $\sigma=\sum_{i=1}^4m_i\sigma_i$, $c(\delta_{\alpha_3}^+,\sigma)=0$ implies $m_3=m_4$. Therefore, one gets $\sigma=m_1(\alpha_1+\alpha_2)+m_2(\alpha_2+\alpha_3)+m_3(\alpha_3+\alpha_4)$ and obtains as quotient the spherical system $\mathcal S_2$, with defect 1. 

\bigskip
There is another well known algebraic subgroup of $G$ whose dimension
is $\geq24$, the symmetric subgroup $K$ of $H$ which is simple of type
$\mathsf D_4$. This $K$ has dimension 28, but is \textit{not}
spherical in $G$. Indeed, if $K$ were spherical in $G$, it would be
wonderful in $G$ (by Corollary~3 of \ref{sphericalclosure2}). But
there is no spherical system of type $\mathsf F_4$ having as quotient
the system $\mathcal S$ and corresponding to a wonderful subgroup
which is reductive. The fact that $K$ is not spherical in $G$ has been known for a long time, in the context of multiplicity-free homogeneous spaces of compact Lie groups (see for instance \cite{Kr79}). This also follows from the fact that $N_G(K)/K$, being isomorphic to the exterior automorphism group of $K$, is not commutative.

%\clearpage
\vspace{6ex}\section{Projective colors and strongly solvable\\ systems}\label{projectivecsss}

In this section we will explain some facts about strongly solvable spherical systems, and illustrate them with examples of type $\mathsf F_4$. 

Let us begin with some details on projective colors.

\subsection{Projective colors}\label{projective}

Let $X$ be a wonderful variety, $\mathcal S=(\Sigma,S^p,\mathbf A)$ its spherical system and $\Delta$ its set of colors. Remember also that a color $\delta\in\mathbf A$ is called \textit{projective} if $c(\delta,\sigma)\geq0$, for all $\sigma\in\Sigma$. If $\delta$ is a projective color, the set $\{\delta\}$ is clearly distinguished in $\Delta$. Let us denote by $\phi\colon X\to \p X$ a wonderful morphism associated to $\{\delta\}$. The term projective color comes from the fact that $\phi$ is a \textit{projective fibration} (i.e.\ is smooth and its fibers are isomorphic to $\mathbb P^n$'s).

Denote by $S_\delta$ the set of $\alpha\in S$ such that $\delta\in\Delta(\alpha)$. Then $\delta$ is also called an \textit{$n$-comb}, where $n=\mathrm{card}(S_\delta)$. Here is an example with 2-comb $\delta_{\alpha_1}^+=\delta_{\alpha_3}^+$:
\[\begin{picture}(6000,2250)\put(300,900){\usebox{\dynkinf}}\multiput(0,0)(3600,0){2}{\usebox{\aone}}\multiput(300,1800)(3600,0){2}{\line(0,1){450}}\put(300,2250){\line(1,0){3600}}\put(1800,600){\usebox{\GreyCircle}}\put(5400,0){\usebox{\aone}}\put(4900,1300){\usebox{\tow}}\end{picture}\]

Wonderful morphisms given by a projective color are the simplest examples of $\mathcal L$-type minimal morphisms. All colors $\delta_\alpha$, $\alpha\in S_\delta$, appear as negative colors in $\p X$, and the difference $d(\p X)-d(X)$ is equal to $n-1$.

We will now explain how one can reduce, in some sense, $n$-combs to 1-combs. This will prepare a similar reduction of $\mathcal L$-type minimal morphisms of strictly increasing defect, to those of constant defect (see Section~\ref{increasing}).

Let $\delta$ be an $n$-comb. Remember that $S_\delta\subset\Sigma$. For each $\alpha\in S_\delta$, denote by $\mathcal S_{\alpha}$ the spherical system obtained from $\mathcal S$ by localization with respect to $(\Sigma\setminus S_\delta)\cup\{\alpha\}$. The corresponding variety is a wonderful $G$-subvariety of $X$ of codimension $n-1$, which we will denote by $X_{\alpha}$. The restriction of $\phi$ to $X_{\alpha}$ gives wonderful morphisms $\phi\colon X_{\alpha}\to \p X$ which come from a 1-comb on $X_{\alpha}$, for all $\alpha\in S_\delta$. In our example above, $S_\delta=\{\alpha_1,\alpha_3\}$, and on the digram level, the (two) morphisms $\phi\colon X_\alpha\to\p X$ are as follows:

\[\begin{picture}(16800,7800)
\put(0,6000){
\put(300,900){\usebox{\dynkinf}}\multiput(0,0)(5400,0){2}{\usebox{\aone}}\put(1800,600){\usebox{\GreyCircle}}\put(3900,900){\circle{600}}
}
\put(5400,0){
\put(300,900){\usebox{\dynkinf}}\multiput(300,900)(3600,0){2}{\circle{600}}\put(1800,600){\usebox{\GreyCircle}}\put(5400,0){\usebox{\aone}}
}
\put(10800,6000){
\put(300,900){\usebox{\dynkinf}}\multiput(3600,0)(1800,0){2}{\usebox{\aone}}\put(1800,600){\usebox{\GreyCircle}}\put(300,900){\circle{600}}\put(4900,1300){\usebox{\tow}}
}
\put(3600,5400){\vector(1,-1){3000}}
\put(13200,5400){\vector(-1,-1){3000}}
\end{picture}\]

We will now explain how one can understand $\phi\colon X\to \p X$ if one knows the $\phi\colon X_{\alpha}\to \p X$, $\alpha\in S_\delta$. Let us choose generic stabilizers $H_{\alpha}$ and $\p H$ in $X_{\alpha}$ and $\p X$. One can assume these $H_{\alpha}$ contained in $\p H$, and that all these groups have a same Levi subgroup $\p L$; moreover, for each $\alpha\in S_\delta$, $\mathrm{Lie}(H_{\alpha}^u)$ is a subspace of codimension 1 in $\mathrm{Lie}(\p H^u)$, containing $\mathrm{Lie}(\p H^u,\p H^u)$, and stable by $\p L$. This gives in particular characters $\tilde\alpha$ of $\p L$ (since $\p L$ acts on the 1-dimensional space $\mathrm{Lie}(\p H^u)/\mathrm{Lie}(H_{\alpha}^u)$). 

Then one can choose $H$ in the following form: $H=L\,H^u$, where $L$ is the subgroup of $\p L$ where all the characters $\tilde\alpha$, $\alpha\in S_\delta$, coincide; and $\mathrm{Lie}(H^u)$ is a codimension 1 subspace of $\mathrm{Lie}(\p H^u)$, which contains the intersection of all $\mathrm{Lie}(H_{\alpha}^u)$, but otherwise is in general position (all these $\mathrm{Lie}(H^u)$ are conjugated under $\p L$). 

On the other hand, if one introduces $H$ by means of this definition, one can check that $H$ is wonderful, and that the wonderful completion of $G/H$ has spherical system $\mathcal S$. 

In other words, in order to prove the uniqueness and existence for geometric realizations of $\mathcal S$, it is necessary and sufficient to do it for the systems $\mathcal S_{\alpha}$, $\alpha\in S_\delta$.

\subsection{Strongly solvable systems}

Remember that $\mathcal S$ is called strongly solvable, if there exists a distinguished subset of colors $\p \Delta$ such that $\mathcal S/\p \Delta=(\emptyset,\emptyset,\emptyset)$ (which is the spherical system of the flag variety $G/B$). By definition, $\mathcal S$ is strongly solvable if and only if a generic stabilizer $H$ of $G$ in $X$ is contained in $B$. In Table~\ref{tr4o}, one can find all 38 strongly solvable rank 4 spherical systems of type $\mathsf F_4$.

If $\mathcal S$ is strongly solvable, then $\Sigma\subset S$ and $S^p=\emptyset$, but the converse is not true in general. The following conditions are equivalent:
\begin{itemf}
\item[(1)] $\mathcal S$ is strongly solvable;
\item[(2)] there exists a sequence of successive quotients by projective colors starting at $\mathcal S$ and ending at $(\emptyset,\emptyset,\emptyset)$;
\item[(3)] any sequence of successive quotients by projective colors starting at $\mathcal S$ can be prolonged until $(\emptyset,\emptyset,\emptyset)$.
\end{itemf}

It has been known for some time that every strongly solvable spherical system (for any semisimple group) is the spherical system of a wonderful variety, uniquely determined up to isomorphism, see \cite{Lu93}.

Here are two examples of strongly solvable systems of type $\mathsf F_4$ with sequences of successive quotients by projective colors:

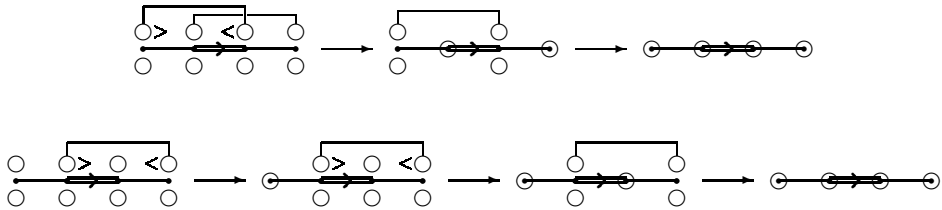
\[\stronglysolvabletwo\]

\[\stronglysolvablethree\]

%\clearpage
\vspace{6ex}\section{Examples coming from type $\mathsf B_3$}\label{sb3}

Remember that $G$ denotes in this chapter a simple group of type
$\mathsf F_4$. In this section, let $Q = Q_{\alpha_4}$ be the parabolic subgroup of $G$ containing $B$ and having semisimple type $\mathsf B_3$. Then $\m Q/\m Q^r = \underline G$ is a simple group of type $\mathsf B_3$. For every wonderful $\underline G$-variety $\underline X$, remember that $X = G \ast_{\m Q} \underline X$ is the wonderful $G$-variety obtained from $\underline X$ by parabolic induction by means of $Q$ (see \ref{ssparabolic}).

On the combinatorial level, this means that the spherical system $\mathcal S=(\Sigma,S^p,\mathbf A)$ of $X$ has an exterior negative color at $\alpha_4$. In the tables of Chapter~1, one can easily find all spherical systems having this property; for instance, Table~\ref{tr3b3} contains all such spherical systems of rank 3.

In this section, we will give some examples of wonderful subgroups
$\underline H$ of $\underline G$ (essentially two reductive spherical
subgroups of $\underline G$, together with those of their parabolic subgroups which are wonderful in $\underline G$), and explicit the corresponding spherical systems of type $\mathsf B_3$.

\subsection{}\label{ss611}

Consider the following figure of diagrams of type $\mathsf B_3$
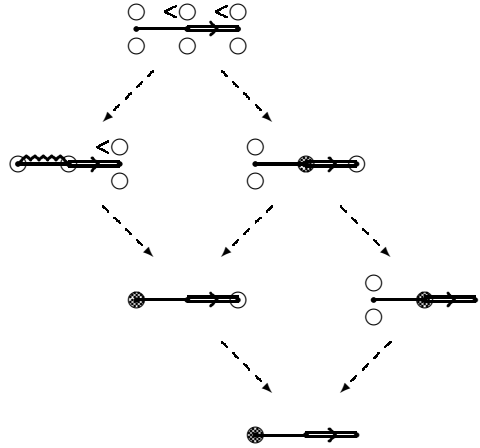
\[\bthreeone\]

The reader is invited to explicit spherical roots, Cartan diagrams, \ldots\ of all the spherical systems corresponding to the diagrams in this figure, and to check that all arrows are associated to well defined minimal quotients of type $\mathcal P$.

The spherical system of the diagram on the bottom line of the figure, together with its double
\[\begin{picture}(3900,1200)
\put(0,-600){
\put(300,900){\usebox{\segm}}
\put(2100,900){\usebox{\rightbisegm}}
\put(0,600){\usebox{\GreyCircle}}
}
\end{picture}
\hspace{2cm}
\begin{picture}(3900,1200)
\put(0,-600){
\put(300,900){\usebox{\segm}}
\put(2100,900){\usebox{\rightbisegm}}
\put(0,600){\usebox{\GreyCircleTwo}}
}
\end{picture}\]
are well known spherical systems of rank 1, with spherical roots $\sigma= \alpha_1+\alpha_2+\alpha_3$ and $2\sigma$. Their corresponding wonderful subgroups $\underline K$ and $\mathrm{N}_{\underline G}(\underline K)$ are symmetric subgroups of $\underline G$, of type $\mathsf D_3 = \mathsf A_3$.

The spherical system of the diagram on the top line of the figure
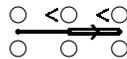
\[\begin{picture}(3900,1800)
\put(300,900){\usebox{\segm}}
\put(2100,900){\usebox{\rightbisegm}}
\multiput(0,0)(1800,0){3}{\usebox{\aone}}
\multiput(1300,1300)(1800,0){2}{\usebox{\tow}}
\end{picture}\]
is strongly solvable. It follows that Borel subgroups $B(\underline K)$ of $\underline K$ are wonderful subgroups in $\underline G$, and have this system as spherical system.

This implies that all parabolic subgroups of $\underline K$ are wonderful in $\underline G$, and that their spherical systems are those of the diagrams in the figure above (there are only 6 systems and 8 parabolic subgroups of $\underline K$ containing $B(\underline K)$, because two couples of the latter set are conjugated in $\underline G$).

The normalizer of a parabolic subgroup of $\underline K$ is not always co-connected in $\mathrm{N}_{\underline G}(\underline K)$. The following figure gives the diagrams of those which are:
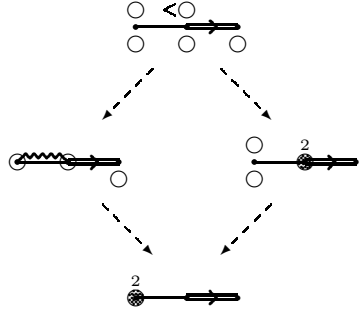
\[\stronebthreeonethree\]

\subsection{}\label{ss614}

\[\begin{picture}(3600,600)
\put(-300,-600){
\put(300,900){\usebox{\segm}}
\put(2100,900){\usebox{\rightbisegm}}
\put(3600,600){\usebox{\GreyCircle}}
}
\end{picture}\]
is the diagram of another well known spherical system of rank 1 of type $\mathsf B_3$, with spherical root $\sigma = \alpha_1+ 2\alpha_2+ 3\alpha_3$. The corresponding wonderful subgroup $\underline H$ of $\underline G$ is simple of type $\mathsf G_2$ (and not symmetric).

The two maximal parabolic subgroups of $\underline H$ have dimension $9 = \dim(B^u(\underline G))$. On the other hand, the following figure contains all spherical systems having as quotient (of type $\mathcal P$) the spherical system of $H$:
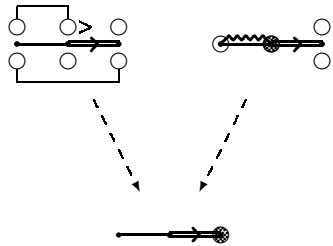
\[\stronebthreeonefour\]

Since we assume the spherical systems above geometrically realizable, it follows that these two parabolic subgroups of $\underline H$ are wonderful in $\underline G$.

\textit{Exercise.} These two parabolic subgroups have same dimension and same semisimple type, but are not isomorphic as groups. Which subgroup corresponds to which diagram?

%\clearpage
\vspace{6ex}\section{Examples coming from type $\mathsf C_3$}\label{ss711}

This section is analogous to Section~\ref{sb3}: the role of $\mathsf B_3$ is now played by $\mathsf C_3$. So now $Q=Q_{\alpha_1}$ will be the parabolic subgroup of $G$ having semisimple type $\mathsf C_3$, and we will set $\m Q/\m Q^r=\underline G$, which is now a simple group of type $\mathsf C_3$. The contrast between sections \ref{sb3} and \ref{ss711} is striking: although simple groups of type $\mathsf B_3$ and type $\mathsf C_3$ are very similar (they are dual groups in the sense of Langlands), the structure of their subgroups is quite different.

We will draw the Dynkin diagram of type $\mathsf C_3$ as follows
\[\begin{picture}(3600,600)
\put(0,300){\usebox{\rightbisegm}}
\put(1800,300){\usebox{\segm}}
\end{picture}\]
keeping the numbering of simple roots $\alpha_2,\alpha_3,\alpha_4$ induced by $\mathsf F_4$.

Consider the following figure of diagrams of type $\mathsf C_3$
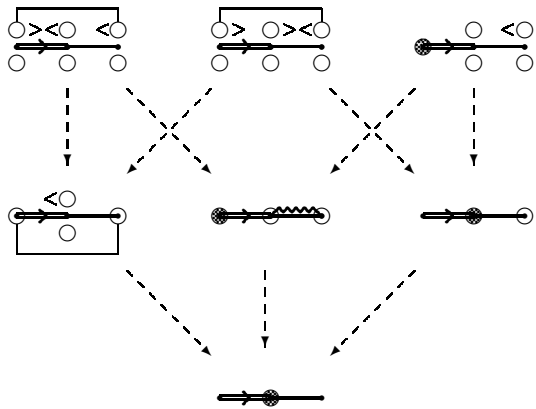
\[\sevenoneone\]

The reader is invited to explicit spherical roots, Cartan diagrams, \ldots\ of all the spherical systems corresponding to the diagrams in this figure, and to check that all arrows come from well defined minimal quotients of type $\mathcal P$.

The spherical system of the diagram on the bottom line has a (unique) spherical root $\sigma = \alpha_2+ 2\alpha_3+ \alpha_4$. It corresponds to a well known symmetric subgroup $\underline K$ of $\underline G$, which is connected and of semisimple type $\mathsf B_2\times\mathsf A_1$.

All parabolic subgroups of $\underline K$, with the exception of the Borel subgroups, have dimension $\geq 9 = \dim(B^u(\underline K))$. Since the figure above contains 7 spherical systems having as quotient of type $\mathcal P$ the system of $\underline K$, and since we assume all these spherical systems geometrically realizable, all parabolic subgroups of $\underline K$, except the Borel subgroups, are wonderful in $\underline G$, and their diagrams are those of the figure above.

%\clearpage
\vspace{6ex}\section{Another remarkable example of rank 2}\label{ss1514}

Consider the figure
\[\begin{picture}(16800,600)\put(300,300){\usebox{\dynkinf}}\multiput(0,0)(3600,0){2}{\usebox{\GreyCircle}}\put(5700,300){\circle{600}}\put(6900,300){\vector(1,0){3000}}\put(10800,0){\put(300,300){\usebox{\dynkinf}}\put(0,0){\usebox{\GreyCircle}}\put(5700,300){\circle{600}}}\end{picture}\]

The first diagram comes from a spherical system $\mathcal S$ given by $\Sigma=\{\alpha_1+\alpha_2+\alpha_3,\alpha_2+2\alpha_3+\alpha_4\}$, $S^p=\{\alpha_2\}$ and $\mathbf A=\emptyset$.

The distinguished set of colors corresponding to the arrow is
$\{\delta_{\alpha_3}\}$, and this quotient is of type $\mathcal
L$. The second diagram has already been analyzed in \ref{ss611}. Let $Q = Q_{\alpha_4}$ be the parabolic subgroup of $G$, and $\underline K$ the subgroup of $\underline G = \m Q/\m Q^u$, as in Section~\ref{sb3}. Put $L = M\cap q^{-1}(\underline K)$, where $M$ is the Levi factor of $Q$ containing $T$, and where $q\colon \m Q\to \underline G$ denotes the canonical map. We know that
\[\mathrm{Lie}(\m Q^u) = \mathfrak n_{-1}+ \mathfrak n_{-2}\]
where $\mathfrak n_{-1}$ and $\mathfrak n_{-2}$ are simple $M$-modules of dimension 8 and 7 (see \ref{ssstructure}). Under the action of $L$, $\mathfrak n_{-1} = \mathfrak n_{-1}'\oplus \mathfrak n_{-1}''$ splits into two $L$-submodules of dimension 4. 

Let us define $H$ by $H=H^uL$, where $\mathrm{Lie}(H^u)$ is either $\mathfrak n_{-1}'+ \mathfrak n_{-2}$ or $\mathfrak n_{-1}'' +\mathfrak n_{-2}$ (these two choices are conjugated in $G$). It is not difficult to check that $H$ is a wonderful subgroup of G having $\mathcal S$ as spherical system.

There are only three spherical systems admitting $\mathcal S$ as quotient, with diagrams as follows:
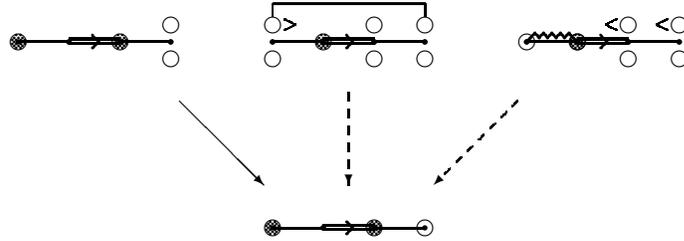
\[\strtwoffouronefour\]
Let us denote them $\mathcal S_i$, $i=1,2,3$ from left to right.

The first quotient $\mathcal S_1\to\mathcal S$ is a quotient by a projective color. The other two quotients $\mathcal S_i\to\mathcal S$ ($i=2,3$) are of type $\mathcal P$. It follows that the two maximal parabolic subgroups of $H$ having semisimple type $\mathsf A_2$ and dimension $24 = \dim(B^u)$, are wonderful in $G$ and have as spherical systems the $\mathcal S_i$ ($i=2,3$) (in particular they are not conjugated in $G$).

%\clearpage
\vspace{6ex}\section{$\mathcal L$-type minimal morphisms with constant defect}\label{constant}

In this section we gather examples of minimal wonderful morphisms $\phi\colon X\to \p X$ of type $\mathcal L$, such that $d(\p X) = d(X)$. 

Remember that this is equivalent to the existence of generic isotropy groups $\p H$ of $\p X$ and $H$ of $X$, having a common Levi subgroup $L$, such that $H^u \subset \p H^u$ and $\mathrm{Lie}(\p H^u)/\mathrm{Lie}(H^u)$ is a simple $L$-module. 

We will introduce these wonderful varieties and morphisms only by their diagrams, without going into details. The reader is of course invited to explicit the corresponding spherical roots, Cartan pairings, distinguished subset of
colors, and to verify that these morphisms are indeed minimal, of type $\mathcal L$ and have constant defect.

We have already seen several such examples, in Section~\ref{ss1521}:

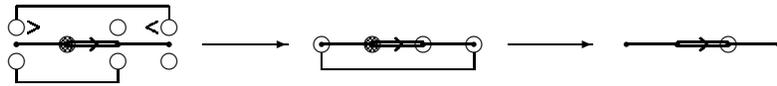
\[\constantone\]

and in Section~\ref{ss1514}:

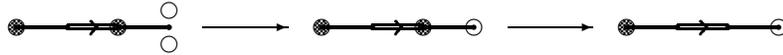
\[\constanttwo\]

In what follows, we will give similar diagrams in particular for all the spherical systems we have come across in Sections \ref{ss1521}, \ref{ss811} and \ref{ss1514} (arrows going from left to right will always correspond to minimal wonderful morphisms of type $\mathcal L$ with constant defect).

\subsection{}

Here is an example close to the spherical system of Section~\ref{ss1521} (having an $\mathcal L$-type minimal quotient which is homogeneous):

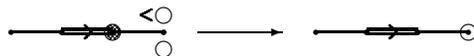
\[\constantthree\]

Notice that the wonderful variety corresponding to the diagram of rank 2 above has isotropy groups having semisimple type $\mathsf B_3$, $\mathsf G_2$ and $\mathsf B_2$.

\subsection{}

Here are examples admitting $\mathcal L$-type minimal wonderful morphisms onto wonderful varieties coming from type $\mathsf B_3$:

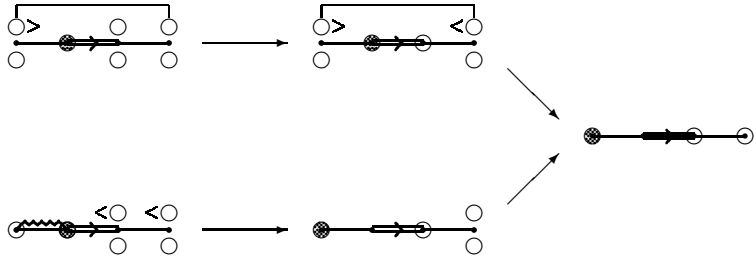
\[\constantfour\]

We have met already the diagrams of rank 4 above in Section~\ref{ss1514}, and that of rank 3 on the top line in Section~\ref{ss1521}.

\subsection{}

Here are examples admitting $\mathcal L$-type minimal wonderful morphisms onto wonderful varieties coming from type $\mathsf C_3$:

1) 
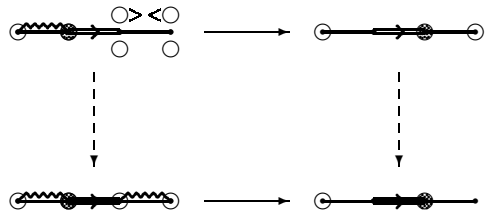
\[\constantfive\]

The vertical arrows (corresponding to morphisms of type $\mathcal P$) have appeared already in sections \ref{ss811} and \ref{ss711}.

2) 
\[\constantsix\]

Of the two diagrams of rank 4 above, one appeared already in Section~\ref{ss1521}, while the other is new in Chapter~3.

\newpage
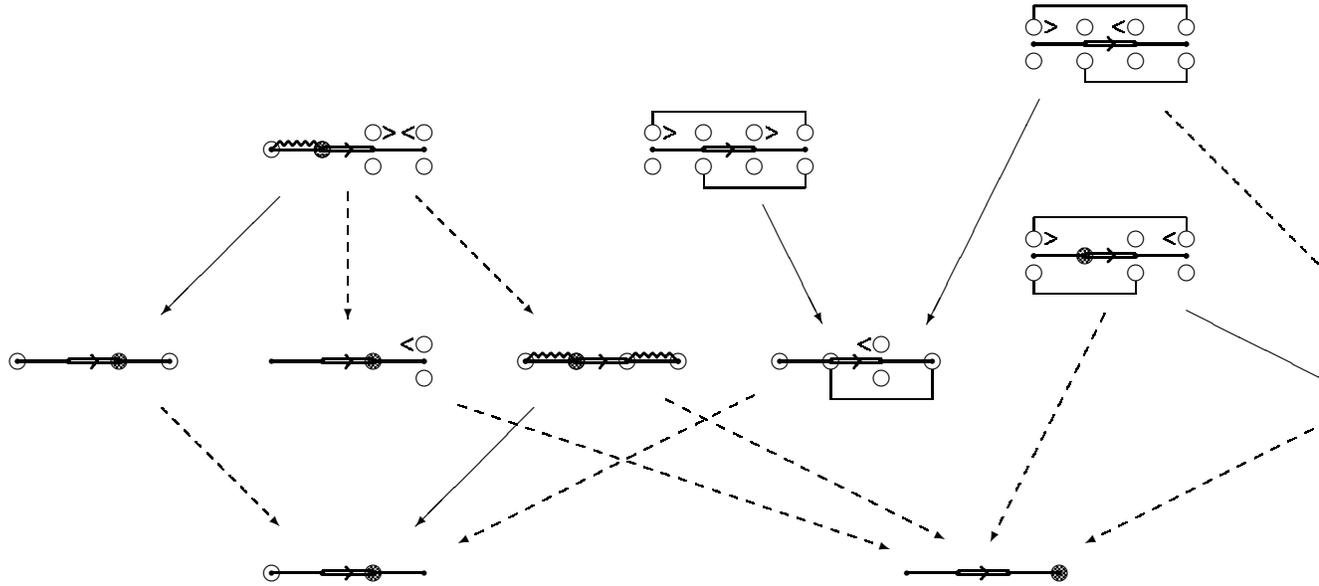
\begin{picture}(90000,25500)\put(3000,0){\thebigone}\end{picture}

\bigskip
Here is some information on the quotients of the last figure on the previous page (notice that this figure appears also as part of the ``big'' figure above), two quotients which seem to be closely similar. Let us choose wonderful subgroups $H_i$ ($i=1,2$) associated to the two diagrams of rank 4, and $H$ a wonderful subgroup associated to their common quotient of rank 2, all three having a common Levi factor $L$. Then $\mathrm{Lie}(H^u)/\mathrm{Lie}(H_i^u)$ ($i=1,2$) are simple $L$-modules of dimension 2. What happens here is that $L$, which is of semisimple type $\mathsf A_1\times \mathsf A_1$, acts on the two $\mathrm{Lie}(H^u)/\mathrm{Lie}(H_i^u)$ ($i=1,2$) via homomorphisms $L\to \mathrm{GL}(2)$ which are trivial on different semisimple factors of $L$.

\bigskip
Let us mention also two other minimal $\mathcal L$-type quotients:
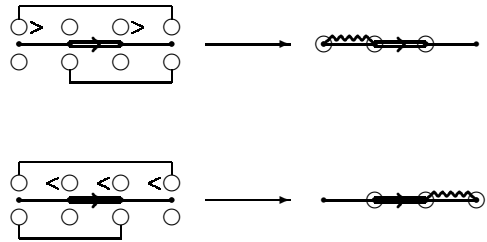
\[\constantseven\]

\[\constanteight\]

We have met the second diagram of rank 4 already in Section~\ref{ss1521}.

\newpage
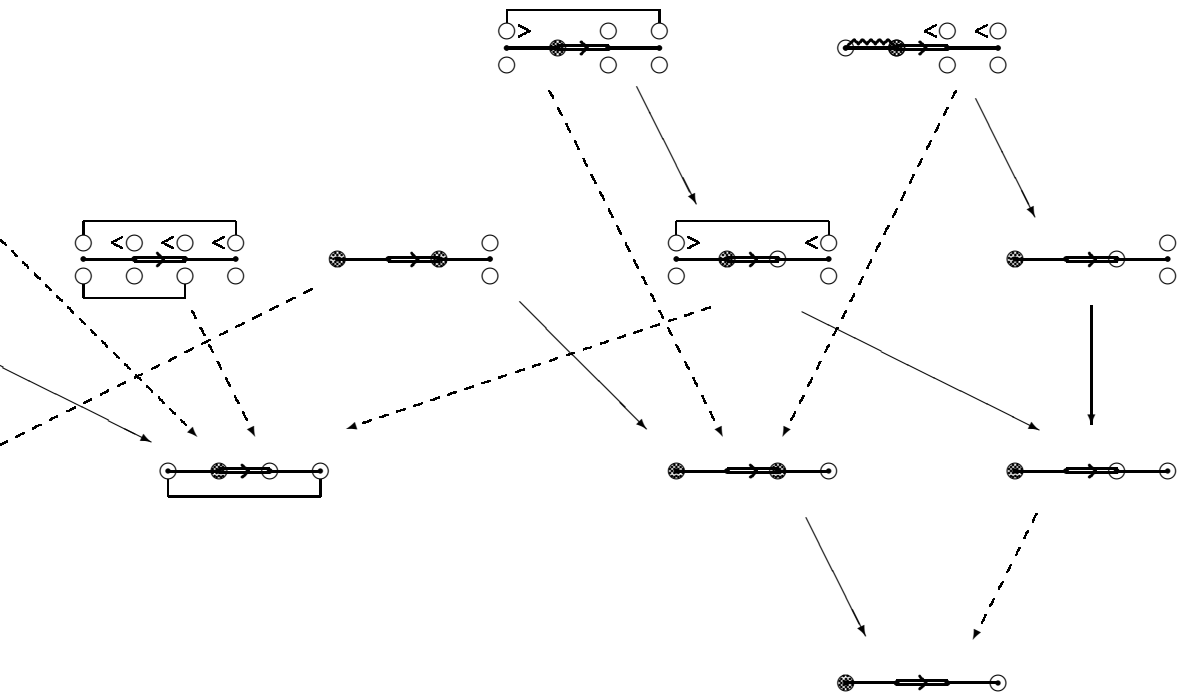
\begin{picture}(0,25500)\put(-53466,0){\thebigone}\end{picture}

\subsection{}

In the ``big'' two-pages figure above we have assembled $\mathcal P$-type minimal wonderful morphisms introduced in Sections \ref{ss1521}, \ref{ss811} and \ref{ss1514}, together with some of the $\mathcal L$-type minimal wonderful morphisms seen in this section.

%\clearpage
\vspace{6ex}\section{$\mathcal L$-type minimal morphisms with strictly\\ increasing defect}\label{increasing}

In the preceding section we have given examples of $\mathcal L$-type minimal (wonderful) morphisms, with same defect on source and target. In this section, we give examples of $\mathcal L$-type minimal morphisms where the defect increases by 1. To each of these examples, we will attach other $\mathcal L$-type minimal morphisms with constant defects, in a way somewhat similar to the one we have used in Section~\ref{projective} to reduce $n$-combs to 1-combs.

Remember that $X$ denotes a wonderful variety, that $\mathcal S=(\Sigma,S^p,\mathsf A)$ is its spherical system and $\Delta$ its set of colors.

\subsection{}

Let us start with the example of the spherical system $\mathcal S$ given by the following diagram:

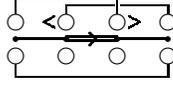
\[\begin{picture}(6000,2850)
\put(300,1350){\usebox{\dynkinf}}
\multiput(0,450)(1800,0){4}{\usebox{\aone}}
\multiput(300,2250)(3600,0){2}{\line(0,1){600}}
\multiput(2100,2250)(3600,0){2}{\line(0,1){300}}
\put(300,2850){\line(1,0){3600}}
\put(2100,2550){\line(1,0){1700}}
\put(4000,2550){\line(1,0){1700}}
\multiput(300,0)(5400,0){2}{\line(0,1){450}}
\put(300,0){\line(1,0){5400}}
\put(1300,1750){\usebox{\tow}}
\put(4300,1750){\usebox{\toe}}
\end{picture}\]

and with the morphisms $\phi:X\to \p X$ given by the distinguished subset of colors $\Delta^\ast=\{\delta_{\alpha_1}^+,\delta_{\alpha_2}^+\}$. Here $d(X)=1$ and $d(\p X)=2$. Let us denote by $S^\ast$ the set of $\alpha\in S$ such that $\delta_{\alpha}$ is a negative exterior color of $\p X$ (but not of $X$). In this case, we have $S^\ast=\{\alpha_2,\alpha_3\}\subset\Sigma$. For each $\alpha\in S^\ast$, denote  by $\mathcal S_{\alpha}$ the spherical system obtained from $\mathcal S$ by localization with respect to $(\Sigma\setminus S^\ast)\cup\{\alpha\}$. The corresponding variety is a wonderful $G$-subvariety of $X$ of codimension $1$, which we will denote by $X_{\alpha}$. For every $\alpha\in S^\ast$, the restriction of $\phi$ to $X_{\alpha}$ gives a wonderful morphism $\phi\colon X_{\alpha}\to \p X$ which is minimal of type $\mathcal L$ with constant defect.

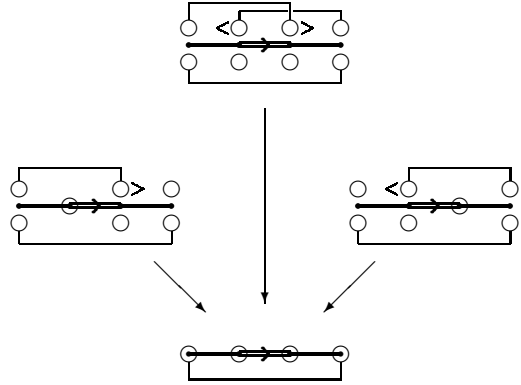
\[\increasingone\]

As in Section~\ref{projective}, let us choose generic stabilizers $H_{\alpha}$ and $\p H$ in $X_{\alpha}$ and $\p X$. One can assume these $H_{\alpha}$ contained in $\p H$, and that all these groups have a same Levi subgroup $\p L$; moreover, for each $\alpha\in S^\ast$, $\mathrm{Lie}(H_{\alpha}^u)$ is here a subspace of codimension 2 in $\mathrm{Lie}(\p H^u)$, containing $\mathrm{Lie}(\p H^u,\p H^u)$ and stable by $\p L$, and the quotient $\mathrm{Lie}(\p H^u)/\mathrm{Lie}(H_{\alpha}^u)$ is a simple $\p L$-module. This last fact gives in particular a character $\tilde\alpha$ of $\ \p L^r$ (which is the connected center of $\ \p L$). What is also important here is that the two $(\p L,\p L)$-modules $\mathrm{Lie}(\p H^u)/\mathrm{Lie}(H^u_\alpha)$, $\alpha\in S^\ast$, are isomorphic. 

Then one can choose $H$ in the following form: $H=\p T(\p L,\p L)\,H^u$, where $\p T$ is the subgroup of $\p L^r$ where the characters $\tilde\alpha$, $\alpha\in S^\ast$, coincide; and $\mathrm{Lie}(H^u)$ is a codimension 2 subspace of $\mathrm{Lie}(\p H^u)$, stable by $(\p L,\p L)$, which contains the intersection of the two $\mathrm{Lie}(H_{\alpha}^u)$, but otherwise is in general position (all these $\mathrm{Lie}(H^u)$ are conjugated by $\p L$). 

On the other hand, if one introduces $H$ by means of this definition, one can check that $H$ is wonderful, and that the wonderful completion of $G/H$ has spherical system $\mathcal S$.

\subsection{}

Here are two other similar examples of type $\mathsf F_4$:

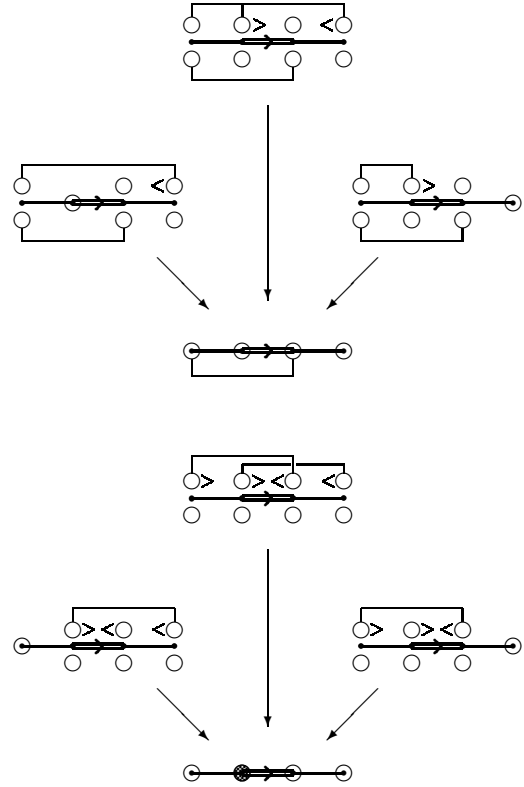
\[\increasingtwo\]

\[\increasingthree\]

\bigskip
$\mathcal L$-type minimal morphisms with strictly increasing defect also exist when $\Sigma\cap S=\emptyset$. To show this let us leave type $\mathsf F_4$. Consider for instance the following example of type $\mathsf B_5$:

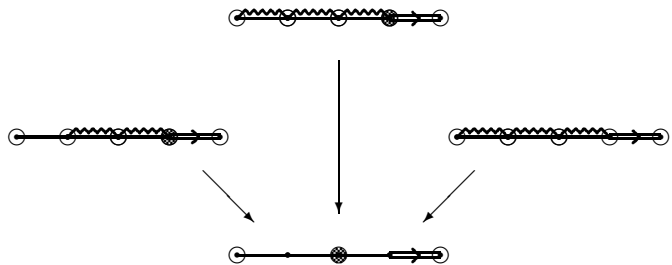
\[\increasingfour\]

%\clearpage
\vspace{6ex}\section{An example of a fiber product}\label{fibered}

Let $\mathcal S$ be a spherical system, $\Delta$ its set of colors, and $X$ a corresponding wonderful variety.

Let $\Delta_1$ and $\Delta_2$ be two distinguished subsets of $\Delta$. We set $\Delta_{1\,2} = \Delta_1\cup\Delta_2$, subset which is again distinguished in $\Delta$. These three distinguished subsets of $\Delta$ give a (commutative) diagram of wonderful morphisms:

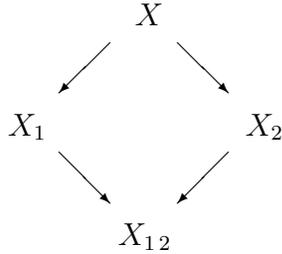
\[\begin{picture}(9600,8700)
\put(4500,7800){$X$}
\put(0,3900){$X_1$}
\put(8400,3900){$X_2$}
\put(3900,0){$X_{1\,2}$}
\multiput(1800,3300)(4200,3900){2}{\vector(1,-1){1800}}
\multiput(7800,3300)(-4200,3900){2}{\vector(-1,-1){1800}}
\end{picture}\]

diagram which induces a $G$-morphism $\psi\colon X\to
X_1\times_{X_{1\,2}} X_2$. In general, this fiber product has no
reason to be a wonderful variety, and $\psi$ has no reason to be an
isomorphism. But sometimes this happens, and then we say that $\Delta_1$ and $\Delta_2$ decompose $\mathcal S$. We will not discuss here the combinatorial conditions on $\mathcal S$, $\Delta_1$, $\Delta_2$ corresponding to this notion of decomposition. We will give only a simple example of type $\mathsf F_4$. 

Consider the spherical system  $\mathcal S = (\Sigma,S^p,\emptyset)$, where $\Sigma = \{\alpha_1+\alpha_2,\alpha_3+\alpha_4\}$.
Here $\Delta = \{\delta_{\alpha_1},\delta_{\alpha_2},\delta_{\alpha_3},\delta_{\alpha_4}\}$. If we choose $\Delta_1=\{\delta_{\alpha_1}\}$ and $\Delta_2 = \{\delta_{\alpha_4}\}$, on the combinatorial level we get the following figure: 

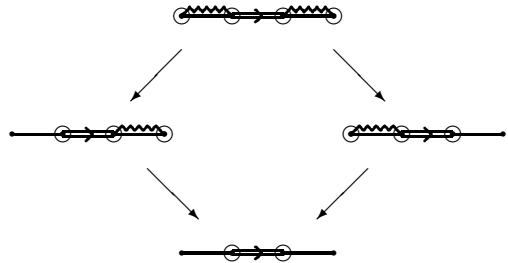
\[\begin{picture}(18000,9600)
\put(6000,300){
\put(300,300){\usebox{\dynkinf}}
\multiput(2100,300)(1800,0){2}{\circle{600}}
}
\put(0,4500){
\put(300,300){\usebox{\dynkinf}}
\put(3600,0){\usebox{\atwo}}
\put(2100,300){\circle{600}}
}
\put(12000,4500){
\put(300,300){\usebox{\dynkinf}}
\put(0,0){\usebox{\atwo}}
\put(3900,300){\circle{600}}
}
\put(6000,8700){
\put(300,300){\usebox{\dynkinf}}
\multiput(0,0)(3600,0){2}{\usebox{\atwo}}
}
\multiput(5100,3600)(6600,4200){2}{\vector(1,-1){1800}}
\multiput(12900,3600)(-6600,4200){2}{\vector(-1,-1){1800}}
\end{picture}\]

In this case the variety $X_{1\,2}$ is homogeneous, isomorphic to $G/\m Q$, where $Q$ is a parabolic subgroup of $G$ generated by $P_{\alpha_1}$ and $P_{\alpha_4}$, and one checks easily that
\[X\to X_1\times_{G/\m Q}X_2\]
is an isomorphism.

%\clearpage
\vspace{6ex}\section{Spherical orbits in fundamental\\ representations}\label{ssfundamental}

In this section we will apply results of Chapter~2 to describe all spherical orbits in the projective fundamental representations of type $\mathsf F_4$.

We will use the terminology and the notation of Section~\ref{ssorbits}. Remember that $\omega_i$ (i=1,2,3,4) are the fundamental weights (of type $\mathsf F_4$).

\begin{proposition}
The faithful couples $(\mathcal S , \delta )$ such that $\omega(\delta) = \omega_i$ are:

if $i=1$
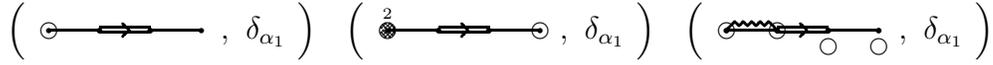
\[\begin{picture}(34800,1800)
\put(0,600){$\Big($}
\put(1200,0){
\put(300,900){\usebox{\dynkinf}}
\put(300,900){\circle{600}}
}
\put(7200,600){$\ ,\ \delta_{\alpha_1}\ \Big)$}
\put(12000,0){
\put(0,600){$\Big($}
\put(1200,0){
\put(300,900){\usebox{\dynkinf}}
\put(0,600){\usebox{\GreyCircleTwo}}
\put(5700,900){\circle{600}}
}
\put(7200,600){$\ ,\ \delta_{\alpha_1}\ \Big)$}
}
\put(24000,0){
\put(0,600){$\Big($}
\put(1200,0){
\put(300,900){\usebox{\dynkinf}}
\put(0,600){\usebox{\atwo}}
\multiput(3600,0)(1800,0){2}{\usebox{\aprime}}
}
\put(7200,600){$\ ,\ \delta_{\alpha_1}\ \Big)$}
}
\end{picture}\]

if $i=2$
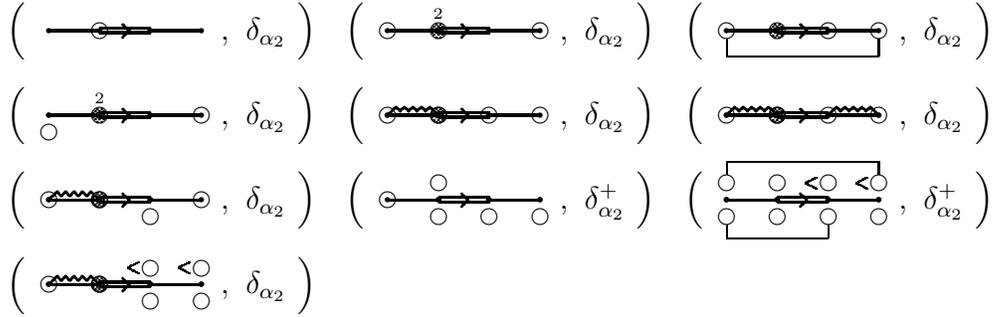
\[\begin{picture}(34800,10800)
\put(0,9000){
\put(0,600){$\Big($}
\put(1200,0){
\put(300,900){\usebox{\dynkinf}}
\put(2100,900){\circle{600}}
}
\put(7200,600){$\ ,\ \delta_{\alpha_2}\ \Big)$}
}
\put(12000,9000){
\put(0,600){$\Big($}
\put(1200,0){
\put(300,900){\usebox{\dynkinf}}
\multiput(300,900)(5400,0){2}{\circle{600}}
\put(1800,600){\usebox{\GreyCircleTwo}}
}
\put(7200,600){$\ ,\ \delta_{\alpha_2}\ \Big)$}
}
\put(24000,9000){
\put(0,600){$\Big($}
\put(1200,0){
\put(300,900){\usebox{\dynkinf}}
\multiput(300,900)(5400,0){2}{\circle{600}}
\multiput(300,0)(5400,0){2}{\line(0,1){600}}
\put(300,0){\line(1,0){5400}}
\put(1800,600){\usebox{\GreyCircle}}
\put(3900,900){\circle{600}}
}
\put(7200,600){$\ ,\ \delta_{\alpha_2}\ \Big)$}
}
\put(0,6000){
\put(0,600){$\Big($}
\put(1200,0){
\put(300,900){\usebox{\dynkinf}}
\put(0,0){\usebox{\aprime}}
\put(1800,600){\usebox{\GreyCircleTwo}}
\put(5700,900){\circle{600}}
}
\put(7200,600){$\ ,\ \delta_{\alpha_2}\ \Big)$}
}
\put(12000,6000){
\put(0,600){$\Big($}
\put(1200,0){
\put(300,900){\usebox{\dynkinf}}
\put(0,600){\usebox{\atwo}}
\put(1800,600){\usebox{\GreyCircle}}
\multiput(3900,900)(1800,0){2}{\circle{600}}
}
\put(7200,600){$\ ,\ \delta_{\alpha_2}\ \Big)$}
}
\put(24000,6000){
\put(0,600){$\Big($}
\put(1200,0){
\put(300,900){\usebox{\dynkinf}}
\multiput(0,600)(3600,0){2}{\usebox{\atwo}}
\put(1800,600){\usebox{\GreyCircle}}
}
\put(7200,600){$\ ,\ \delta_{\alpha_2}\ \Big)$}
}
\put(0,3000){
\put(0,600){$\Big($}
\put(1200,0){
\put(300,900){\usebox{\dynkinf}}
\put(0,600){\usebox{\atwo}}
\put(1800,600){\usebox{\GreyCircle}}
\put(3600,0){\usebox{\aprime}}
\put(5700,900){\circle{600}}
}
\put(7200,600){$\ ,\ \delta_{\alpha_2}\ \Big)$}
}
\put(12000,3000){
\put(0,600){$\Big($}
\put(1200,0){
\put(300,900){\usebox{\dynkinf}}
\put(300,900){\circle{600}}
\put(1800,0){\usebox{\aone}}
\multiput(3600,0)(1800,0){2}{\usebox{\aprime}}
}
\put(7200,600){$\ ,\ \delta_{\alpha_2}^+\ \Big)$}
}
\put(24000,3000){
\put(0,600){$\Big($}
\put(1200,0){
\put(300,900){\usebox{\dynkinf}}
\multiput(0,0)(1800,0){4}{\usebox{\aone}}
\multiput(300,1800)(5400,0){2}{\line(0,1){450}}
\put(300,2250){\line(1,0){5400}}
\multiput(300,0)(3600,0){2}{\line(0,-1){450}}
\put(300,-450){\line(1,0){3600}}
\multiput(3100,1300)(1800,0){2}{\usebox{\tow}}
}
\put(7200,600){$\ ,\ \delta_{\alpha_2}^+\ \Big)$}
}
\put(0,0){
\put(0,600){$\Big($}
\put(1200,0){
\put(300,900){\usebox{\dynkinf}}
\put(0,600){\usebox{\atwo}}
\put(1800,600){\usebox{\GreyCircle}}
\multiput(3600,0)(1800,0){2}{\usebox{\aone}}
\multiput(3100,1300)(1800,0){2}{\usebox{\tow}}
}
\put(7200,600){$\ ,\ \delta_{\alpha_2}\ \Big)$}
}
\end{picture}\]

if $i=3$
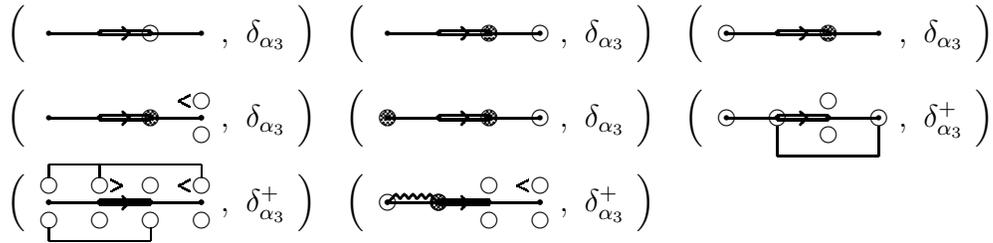
\[\begin{picture}(34800,7800)
\put(0,6000){
\put(0,600){$\Big($}
\put(1200,0){
\put(300,900){\usebox{\dynkinf}}
\put(3900,900){\circle{600}}
}
\put(7200,600){$\ ,\ \delta_{\alpha_3}\ \Big)$}
}
\put(12000,6000){
\put(0,600){$\Big($}
\put(1200,0){
\put(300,900){\usebox{\dynkinf}}
\put(5700,900){\circle{600}}
\put(3600,600){\usebox{\GreyCircle}}
}
\put(7200,600){$\ ,\ \delta_{\alpha_3}\ \Big)$}
}
\put(24000,6000){
\put(0,600){$\Big($}
\put(1200,0){
\put(300,900){\usebox{\dynkinf}}
\put(300,900){\circle{600}}
\put(3600,600){\usebox{\GreyCircle}}
}
\put(7200,600){$\ ,\ \delta_{\alpha_3}\ \Big)$}
}
\put(0,3000){
\put(0,600){$\Big($}
\put(1200,0){
\put(300,900){\usebox{\dynkinf}}
\put(5400,0){\usebox{\aone}}
\put(4900,1300){\usebox{\tow}}
\put(3600,600){\usebox{\GreyCircle}}
}
\put(7200,600){$\ ,\ \delta_{\alpha_3}\ \Big)$}
}
\put(12000,3000){
\put(0,600){$\Big($}
\put(1200,0){
\put(300,900){\usebox{\dynkinf}}
\put(5700,900){\circle{600}}
\multiput(0,600)(3600,0){2}{\usebox{\GreyCircle}}
}
\put(7200,600){$\ ,\ \delta_{\alpha_3}\ \Big)$}
}
\put(24000,3000){
\put(0,600){$\Big($}
\put(1200,0){
\put(300,900){\usebox{\dynkinf}}
\put(300,900){\circle{600}}
\multiput(2100,900)(3600,0){2}{\circle{600}}
\multiput(2100,600)(3600,0){2}{\line(0,-1){1050}}
\put(2100,-450){\line(1,0){3600}}
\put(3600,0){\usebox{\aone}}
}
\put(7200,600){$\ ,\ \delta_{\alpha_3}^+\ \Big)$}
}
\put(0,0){
\put(0,600){$\Big($}
\put(1200,0){
\put(300,900){\usebox{\dynkinf}}
\multiput(0,0)(1800,0){4}{\usebox{\aone}}
\multiput(300,1800)(5400,0){2}{\line(0,1){450}}
\put(2100,1800){\line(0,1){450}}
\put(300,2250){\line(1,0){5400}}
\multiput(300,0)(3600,0){2}{\line(0,-1){450}}
\put(300,-450){\line(1,0){3600}}
\put(2500,1300){\usebox{\toe}}
\put(4900,1300){\usebox{\tow}}
}
\put(7200,600){$\ ,\ \delta_{\alpha_3}^+\ \Big)$}
}
\put(12000,0){
\put(0,600){$\Big($}
\put(1200,0){
\put(300,900){\usebox{\dynkinf}}
\put(0,600){\usebox{\atwo}}
\put(1800,600){\usebox{\GreyCircle}}
\multiput(3600,0)(1800,0){2}{\usebox{\aone}}
\put(4900,1300){\usebox{\tow}}
}
\put(7200,600){$\ ,\ \delta_{\alpha_3}^+\ \Big)$}
}
\end{picture}\]

if $i=4$
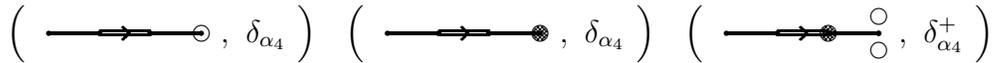
\[\begin{picture}(34800,1800)
\put(0,0){
\put(0,600){$\Big($}
\put(1200,0){
\put(300,900){\usebox{\dynkinf}}
\put(5700,900){\circle{600}}
}
\put(7200,600){$\ ,\ \delta_{\alpha_4}\ \Big)$}
}
\put(12000,0){
\put(0,600){$\Big($}
\put(1200,0){
\put(300,900){\usebox{\dynkinf}}
\put(5400,600){\usebox{\GreyCircle}}
}
\put(7200,600){$\ ,\ \delta_{\alpha_4}\ \Big)$}
}
\put(24000,0){
\put(0,600){$\Big($}
\put(1200,0){
\put(300,900){\usebox{\dynkinf}}
\put(3600,600){\usebox{\GreyCircle}}
\put(5400,0){\usebox{\aone}}
}
\put(7200,600){$\ ,\ \delta_{\alpha_4}^+\ \Big)$}
}
\end{picture}\]

\end{proposition}

All spherical systems mentioned in this proposition have trivial $\Gamma$.

\begin{sketche}
From the combinatorial definition of the map $\omega\colon\mathbb N\Delta\to\mathbb N\Omega$ given in \ref{ssorbits} follows that $\omega(\delta) = \omega_i$ implies either $\delta = \delta_{\alpha_i}$ or $\delta = \delta_{\alpha_i}^+$. A careful combinatorial analysis then gives the lists above.
\end{sketche}

\bigskip
\textit{Remarks}

1) The representation $V(\omega_1)$ is the adjoint representation. It is well known that there are exactly 3 spherical adjoint orbits, all of them nilpotent (see \cite{Pa03} and \cite{CCC05}).

2) Let us give explicitly the spherical orbit in $\mathbb P(V(\omega_2))$ associated to the faithful couple 
\[\begin{picture}(10800,1800)
\put(0,600){$\Big(\ $}
\put(1200,-450){
\put(300,1350){\usebox{\dynkinf}}
\multiput(300,1350)(5400,0){2}{\circle{600}}
\put(3900,1350){\circle{600}}
\multiput(300,450)(5400,0){2}{\line(0,1){600}}
\put(300,450){\line(1,0){5400}}
\put(1800,1050){\usebox{\GreyCircle}}
}
\put(7200,600){$\ ,\ \delta_{\alpha_2}\ \Big)$}
\end{picture}\]
Consider the vector
\[v = X_{-(\alpha_1+2\alpha_2+3\alpha_3+\alpha_4)} \wedge X_{-(\alpha_1+2\alpha_2+3\alpha_3+2\alpha_4)}\in\bigwedge^2 V(\omega_1)\]
(where for every root $\beta$, $X_\beta$ denotes an associated root vector). If $Q$ is a maximal parabolic subgroup of $G$ associated to $\alpha_3$ and containing $B$, it is not difficult to check that $H = (^-Q^u, ^-Q^u)\,M$ stabilizes the line $\mathbb Cv$. We know that this $H$ is a generic isotropy group of the wonderful variety having the above diagram (see Section~\ref{ss1521}). Since \[\bigwedge^2 V(\omega_1) \cong V(\omega_2) \oplus V(\omega_1),\]
and since $(\alpha_1+2\alpha_2+3\alpha_3+\alpha_4)+(\alpha_1+2\alpha_2+3\alpha_3+2\alpha_4)$ is not a root of $G$, necessarily $v\in V(\omega_2)$. It follows that the spherical orbit in $\mathbb P(V(\omega_2))$ given by $G.[v]$ is the one associated to the above faithful couple.

%\clearpage
\vspace{6ex}\section{Examples not of type $\mathsf F_4$}\label{notf4}

In this section, we give several examples of general type to illustrate phenomena which do not occur or cannot be well illustrated in type $\mathsf F_4$.

\subsection{}

There are not so many examples of minimal wonderful $\mathcal R$-type morphisms
when $G$ is of type $\mathsf F_4$. In what follows, we give more examples for other types:
\begin{itemf}

\item[1)] $G$ of type $\mathsf B_{2m}$
\[\begin{picture}(27600,1200)\put(0,0){\diagrambcprimen}
\put(12300,900){\vector(1,0){3000}}
\put(16200,0){
\multiput(300,900)(1800,0){2}{\usebox{\segm}}
\put(3900,900){\usebox{\susp}}
\put(7500,900){\usebox{\segm}}
\put(9300,900){\usebox{\rightbisegm}}
\put(0,600){\usebox{\GreyCircle}}
}
\end{picture}
\]

\item[2)] $G$ of type $\mathsf D_4$
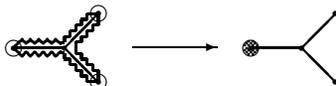
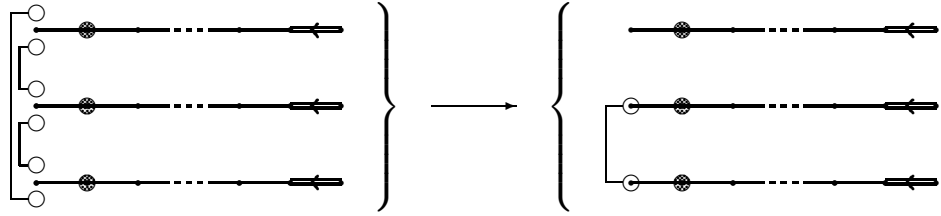
\[\begin{picture}(12000,3000)
\put(0,0){\diagramdsastfour}
\put(4500,1500){\vector(1,0){3000}}
\put(8700,1500){\usebox{\segm}}
\put(10500,300){\usebox{\bifurc}}
\put(8400,1200){\usebox{\GreyCircle}}
\end{picture}\]

\item[3)] $G$ of type $\mathsf C_l\times\mathsf C_m\times\mathsf C_n$%, with $S=\{\alpha_1,\ldots,\alpha_l;\alpha'_1,\ldots,\alpha'_m;\alpha''_1,\ldots,\alpha''_n\}$ ($\alpha_l,\alpha'_m, \alpha''_n$ long roots). $\mathcal S=(\Sigma,S^p,\mathbf A)$ with $\Sigma=\{\sigma_1,\ldots,\sigma_6\}$, where\\
%$\sigma_1=\alpha_1$, $\sigma_2=\alpha'_1$, $\sigma_3=\alpha''_1$;\\
%$\sigma_4=\alpha_1+2\alpha_2+\ldots+2\alpha_{l-1}+\alpha_l$,\\
%$\sigma_5=\alpha'_1+2\alpha'_2+\ldots+2\alpha'_{m-1}+\alpha'_m$,\\
%$\sigma_6=\alpha''_1+2\alpha''_2+\ldots+2\alpha''_{n-1}+\alpha''_n$,\\ 
%$S^p=S\setminus\{\alpha_1,\alpha_2,\alpha'_1,\alpha'_2,\alpha''_1,\alpha''_2\}$ and (restricted) Cartan pairing:
%\[\begin{array}{r|rrr|rrr}c(-,-)&\sigma_1&\sigma_2&\sigma_3&\sigma_4&\sigma_5&\sigma_6\\\hline\delta_{\alpha_1}^+&1&1&-1&0&0&0\\\delta_{\alpha'_1}^+&-1&1&1&0&0&0\\\delta_{\alpha''_1}^+&1&-1&1&0&0&0\end{array}\]
%Take for instance the quotient given by the subset $\{\delta_{\alpha_1}^+,\delta_{\alpha'_1}^+\}$.
\[\left.
\begin{picture}(11700,3900)(0,3300)
\multiput(600,0)(0,2700){3}{
\multiput(300,900)(1800,0){2}{\usebox{\segm}}
\put(3900,900){\usebox{\susp}}
\put(7500,900){\usebox{\segm}}
\put(9300,900){\usebox{\leftbisegm}}
\put(1800,600){\usebox{\GreyCircle}}
\put(0,0){\usebox{\aone}}
}
\multiput(0,300)(0,6600){2}{\line(1,0){600}}
\put(0,300){\line(0,1){6600}}
\multiput(0,0)(0,2700){2}{
\multiput(300,1500)(0,1500){2}{\line(1,0){300}}
\put(300,1500){\line(0,1){1500}}
}
\end{picture}
\quad\right\}
\begin{picture}(4800,3900)(0,3300)
\put(900,3600){\vector(1,0){3000}}
\end{picture}
\left\{\quad
\begin{picture}(11700,3900)(0,3300)
\multiput(600,0)(0,2700){3}{
\multiput(300,900)(1800,0){2}{\usebox{\segm}}
\put(3900,900){\usebox{\susp}}
\put(7500,900){\usebox{\segm}}
\put(9300,900){\usebox{\leftbisegm}}
\put(1800,600){\usebox{\GreyCircle}}
}
\multiput(900,900)(0,2700){2}{\circle{600}}
\multiput(0,900)(0,2700){2}{\line(1,0){600}}
\put(0,900){\line(0,1){2700}}
\end{picture}
\right.\]
\end{itemf}

All the examples above have defect 0, so correspond to (very) reductive wonderful subgroups, which are well known (see \cite{B87}). The following examples are related to wonderful model varieties (see \cite{Lu07}), they have defect 1.

\begin{itemf}
\item[4)] $G$ of type $\mathsf A_{2m}$
\[\begin{picture}(27600,600)\put(0,0){\diagramacastn}
\put(12300,300){\vector(1,0){3000}}
\put(16200,0){
\multiput(0,0)(7200,0){2}{\multiput(300,300)(1800,0){2}{\usebox{\segm}}}
\put(3900,300){\usebox{\susp}}
\multiput(300,300)(10800,0){2}{\circle{600}}
\multiput(300,300)(25,25){13}{\circle*{70}}
\multiput(600,600)(5700,0){2}{\multiput(0,0)(300,0){15}{\multiput(0,0)(25,-25){7}{\circle*{70}}}\multiput(150,-150)(300,0){15}{\multiput(0,0)(25,25){7}{\circle*{70}}}}
\multiput(5400,600)(300,0){3}{\circle*{70}}
\multiput(11100,300)(-25,25){13}{\circle*{70}}
}
\end{picture}
\]

\item[5)] $G$ of type $\mathsf B_{2m}$
\[\begin{picture}(27600,600)\put(0,0){\diagrambcastn}
\put(12300,300){\vector(1,0){3000}}
\put(16200,0){
\multiput(300,300)(1800,0){2}{\usebox{\segm}}
\put(3900,300){\usebox{\susp}}
\put(7500,300){\usebox{\segm}}
\put(9300,300){\usebox{\rightbisegm}}
\put(0,0){\usebox{\GreyCircle}}
\put(11100,300){\circle{600}}
}
\end{picture}\]

\item[6)] $G$ of type $\mathsf D_{2m+1}$
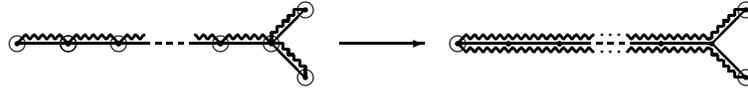
\[\begin{picture}(26400,3000)\put(0,0){\diagramdcastn}
\put(11700,1500){\vector(1,0){3000}}
\put(15600,0){
\multiput(300,1500)(1800,0){2}{\usebox{\segm}}
\put(3900,1500){\usebox{\susp}}
\put(7500,1500){\usebox{\segm}}
\put(9300,300){\usebox{\bifurc}}
\put(300,1500){\circle{600}}
\multiput(10500,300)(0,2400){2}{\circle{600}}
\multiput(300,1500)(25,25){13}{\circle*{70}}
\put(600,1800){\multiput(0,0)(300,0){15}{\multiput(0,0)(25,-25){7}{\circle*{70}}}\multiput(150,-150)(300,0){15}{\multiput(0,0)(25,25){7}{\circle*{70}}}}
\multiput(5400,1800)(300,0){3}{\circle*{70}}
\put(6300,1800){\multiput(0,0)(300,0){10}{\multiput(0,0)(25,-25){7}{\circle*{70}}}\multiput(150,-150)(300,0){9}{\multiput(0,0)(25,25){7}{\circle*{70}}}}
\multiput(300,1500)(25,-25){13}{\circle*{70}}
\put(600,1200){\multiput(0,0)(300,0){15}{\multiput(0,0)(25,25){7}{\circle*{70}}}\multiput(150,150)(300,0){15}{\multiput(0,0)(25,-25){7}{\circle*{70}}}}
\multiput(5400,1200)(300,0){3}{\circle*{70}}
\put(6300,1200){\multiput(0,0)(300,0){10}{\multiput(0,0)(25,25){7}{\circle*{70}}}\multiput(150,150)(300,0){9}{\multiput(0,0)(25,-25){7}{\circle*{70}}}}
\thicklines
\put(3600,0){\put(6900,2700){\line(-1,0){400}}\multiput(6500,2700)(-200,-200){5}{\line(0,-1){200}}\multiput(6500,2500)(-200,-200){4}{\line(-1,0){200}}\multiput(5700,1700)(-30,-10){5}{\line(-1,0){30}}\put(6900,300){\line(-1,0){400}}\multiput(6500,300)(-200,200){5}{\line(0,1){200}}\multiput(6500,500)(-200,200){4}{\line(-1,0){200}}\multiput(5700,1300)(-30,10){5}{\line(-1,0){30}}
}
}
\end{picture}\]
\end{itemf}

In the following example, the three quotients are of type $\mathcal L$, $\mathcal R$ and $\mathcal P$, respectively. Notice that a new interior negative color appears in the quotient of type $\mathcal R$.

\begin{itemf}
\item[7)] $G$ of type $\mathsf A_1\times\mathsf B_2$
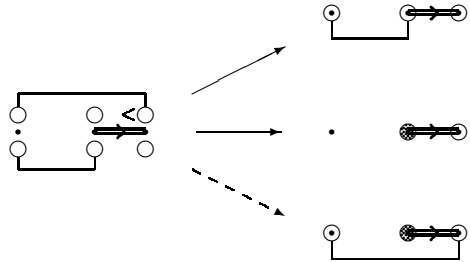
\[\begin{picture}(16200,9000)
\put(0,3150){
\put(0,450){\usebox{\aone}}
\put(3000,1350){\usebox{\rightbisegm}}
\multiput(2700,450)(1800,0){2}{\usebox{\aone}}
\multiput(300,2700)(4500,0){2}{\line(0,-1){450}}
\put(300,2700){\line(1,0){4500}}
\multiput(300,0)(2700,0){2}{\line(0,1){450}}
\put(300,0){\line(1,0){2700}}
\put(4000,1750){\usebox{\tow}}
}
\put(11100,7800){
\put(300,900){\circle*{150}}
\put(3000,900){\usebox{\rightbisegm}}
\multiput(300,900)(2700,0){2}{\circle{600}}
\put(4800,900){\circle{600}}
\multiput(300,0)(2700,0){2}{\line(0,1){600}}
\put(300,0){\line(1,0){2700}}
}
\put(11100,4200){
\put(300,300){\circle*{150}}
\put(3000,300){\usebox{\rightbisegm}}
\put(2700,0){\usebox{\GreyCircle}}
\put(4800,300){\circle{600}}
}
\put(11100,0){
\put(300,900){\circle*{150}}
\put(3000,900){\usebox{\rightbisegm}}
\multiput(300,900)(4500,0){2}{\circle{600}}
\put(2700,600){\usebox{\GreyCircle}}
\multiput(300,0)(4500,0){2}{\line(0,1){600}}
\put(300,0){\line(1,0){4500}}
}
\put(6450,5850){\vector(2,1){3300}}
\put(6600,4500){\vector(1,0){3000}}
\put(6450,3150){
\multiput(0,0)(600,-300){5}{\multiput(0,0)(30,-15){10}{\line(1,0){30}}}
\put(3000,-1500){\vector(2,-1){300}}
}
\end{picture}\]
\end{itemf}

\subsection{}

When $G$ is of type $\mathsf F_4$ (or $\mathsf B_3$, or $\mathsf C_3$), and $K$ is a (connected) spherically closed subgroup of $G$, we have seen (in Sections \ref{ss1521}, \ref{ss811}, \ref{sb3}, \ref{ss711} and \ref{ss1514}) many examples where all parabolic subgroups $H$, of dimension $\geq\dim(B^u)$, are spherical in $G$.

In what follows, we give examples to show that this is not true in general. The argument will always be the same: if $H$ were spherical in $G$, $H$ would be spherically closed and hence wonderful in $G$ (see Corollary~\ref{argument} in \ref{sphericalclosure2}); but in each case, one can check that on the combinatorial level there is simply no spherical system (of right properties to come from $H$) having as quotient the spherical system corresponding to $K$.

\begin{itemf}
\item[1)]
If $G$ is of type $\mathsf D_4$, and $K$ ($\cong \mathrm{GL}(4)$) is the spherically closed subgroup of $G$ having diagram

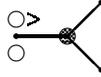
\[\begin{picture}(3300,2400)
\put(300,1200){\usebox{\segm}}
\put(2100,0){\usebox{\bifurc}}
\put(0,300){\usebox{\aone}}
\put(700,1600){\usebox{\toe}}
\put(1800,900){\usebox{\GreyCircle}}
\end{picture}\]

then the parabolic subgroup $H$ of $K$ having semisimple type $\mathsf A_1 \times \mathsf A_1$ and dimension $12 = \dim(B^u)$, is not spherical in $G$. On the other hand, the two parabolic subgroups of $K$ having semisimple type $\mathsf A_2$ (and dimension 13), which are conjugated in $G$, are spherical and wonderful in $G$, and the corresponding diagram is

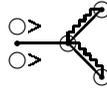
\[\begin{picture}(3600,3000)
\put(300,1500){\usebox{\segm}}
\put(2100,300){\usebox{\bifurc}}
\put(0,600){\usebox{\aone}}
\multiput(700,1900)(0,-1200){2}{\usebox{\toe}}
\put(2100,1500){\circle{600}}
\multiput(3300,300)(0,2400){2}{\circle{600}}
\thicklines\put(-7200,0){\put(9300,1500){\line(0,1){400}}\multiput(9300,1900)(200,200){4}{\line(1,0){200}}\multiput(9500,1900)(200,200){4}{\line(0,1){200}}\put(10500,2700){\line(-1,0){400}}\put(9300,1500){\line(1,0){400}}\multiput(9700,1500)(200,-200){4}{\line(0,-1){200}}\multiput(9700,1300)(200,-200){4}{\line(1,0){200}}\put(10500,300){\line(0,1){400}}}
\end{picture}\]

\item[2)]
If $G$ is of type $\mathsf A_5$ ($\cong\mathrm{SL}(6)$), and $K$ ($\cong (\mathrm{GL}(2)\!\times\!\mathrm{GL}(4))\cap\mathrm{SL}(6)$) is the spherically closed subgroups of $G$ having diagram

\[\begin{picture}(7800,1200)
\multiput(300,900)(1800,0){4}{\usebox{\segm}}
\multiput(300,900)(1800,0){2}{\multiput(0,0)(5400,0){2}{\circle{600}}}
\multiput(300,0)(7200,0){2}{\line(0,1){600}}
\put(300,0){\line(1,0){7200}}
\put(1800,600){\multiput(300,300)(25,25){13}{\circle*{70}}\put(600,600){\multiput(0,0)(300,0){10}{\multiput(0,0)(25,-25){7}{\circle*{70}}}\multiput(150,-150)(300,0){10}{\multiput(0,0)(25,25){7}{\circle*{70}}}}\multiput(3900,300)(-25,25){13}{\circle*{70}}}
\end{picture}\]

then the two parabolic subgroups of $K$ having semisimple type respectively $\mathsf A_2$ and $\mathsf A_1\times\mathsf A_1 \times \mathsf A_1$ and dimension $15 = \dim(B^u)$ are not spherical in $G$. On the other hand, the three (pairwise not $G$-conjugated) parabolic subgroups of $K$ having semisimple type $\mathsf A_3$ or $\mathsf A_1\times\mathsf A_2$ (and dimension greater than 15) are spherical and wonderful in $G$, and the corresponding diagrams are

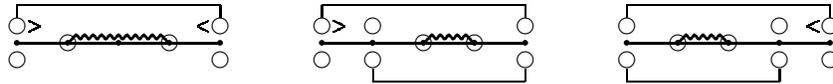
\[\begin{picture}(29400,2700)
\multiput(300,1350)(1800,0){4}{\usebox{\segm}}
\multiput(0,450)(7200,0){2}{\usebox{\aone}}
\multiput(2100,1350)(3600,0){2}{\circle{600}}
\multiput(300,2700)(7200,0){2}{\line(0,-1){450}}
\put(300,2700){\line(1,0){7200}}
\put(700,1750){\usebox{\toe}}
\put(6700,1750){\usebox{\tow}}
\put(1800,1050){\multiput(300,300)(25,25){13}{\circle*{70}}\put(600,600){\multiput(0,0)(300,0){10}{\multiput(0,0)(25,-25){7}{\circle*{70}}}\multiput(150,-150)(300,0){10}{\multiput(0,0)(25,25){7}{\circle*{70}}}}\multiput(3900,300)(-25,25){13}{\circle*{70}}}
\put(10800,0){
\multiput(300,1350)(1800,0){4}{\usebox{\segm}}
\multiput(0,450)(7200,0){2}{\usebox{\aone}}
\put(1800,450){\usebox{\aone}}
\put(3600,1050){\usebox{\atwo}}
\multiput(300,2700)(7200,0){2}{\line(0,-1){450}}
\put(300,2700){\line(1,0){7200}}
\multiput(2100,0)(5400,0){2}{\line(0,1){450}}
\put(2100,0){\line(1,0){5400}}
\put(700,1750){\usebox{\toe}}
}
\put(21600,0){
\multiput(300,1350)(1800,0){4}{\usebox{\segm}}
\multiput(0,450)(7200,0){2}{\usebox{\aone}}
\put(5400,450){\usebox{\aone}}
\put(1800,1050){\usebox{\atwo}}
\multiput(300,2700)(7200,0){2}{\line(0,-1){450}}
\put(300,2700){\line(1,0){7200}}
\multiput(300,0)(5400,0){2}{\line(0,1){450}}
\put(300,0){\line(1,0){5400}}
\put(6700,1750){\usebox{\tow}}
}
\end{picture}\]

\item[3)]
If $G$ is of type $\mathsf E_6$, and $K$ (of type $\mathsf F_4$) is the spherically closed subgroup of $G$ having diagram

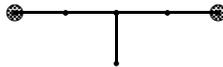
\[\diagramefsix\]

then the parabolic subgroup $H$ of $K$ having semisimple type $\mathsf B_3$ and dimension $37 = \dim(B^u)+1$, is not spherical in $G$. On the other hand, the parabolic subgroup of $K$ having semisimple type $\mathsf C_3$ (and \textit{same} dimension 37) is spherical and wonderful in $G$, and the corresponding diagram is

\[\begin{picture}(7800,2400)\multiput(300,2100)(1800,0){4}{\usebox{\segm}}\put(3900,300){\usebox{\vsegm}}\multiput(0,1800)(1800,0){4}{\multiput(300,300)(1800,0){2}{\circle{600}}\multiput(300,300)(25,25){13}{\circle*{70}}\multiput(600,600)(300,0){4}{\multiput(0,0)(25,-25){7}{\circle*{70}}}\multiput(750,450)(300,0){4}{\multiput(0,0)(25,25){7}{\circle*{70}}}\multiput(2100,300)(-25,25){13}{\circle*{70}}}\put(3900,300){\circle{600}}\put(3600,0){\multiput(300,300)(25,25){13}{\circle*{70}}\multiput(600,600)(0,300){4}{\multiput(0,0)(-25,25){7}{\circle*{70}}}\multiput(450,750)(0,300){4}{\multiput(0,0)(25,25){7}{\circle*{70}}}\multiput(300,2100)(25,-25){13}{\circle*{70}}}\end{picture}\]

\end{itemf}

%\clearpage

\vspace{3ex}
\textit{Acknowledgments.} The first named author would like to thank the Institut Fourier of Grenoble for its hospitality and the European Commission for financial support under the program ``Marie Curie Actions FP7-PEOPLE-2007-2-1-IEF''.

\vspace{3ex}
{\it
P.~Bravi and D.~Luna

Institut Fourier\\
Universit\'e Grenoble 1\\
100 rue des Maths\\
38402 St Martin d'H\`eres\\
France}

\verb|bravi@fourier.ujf-grenoble.fr|\\
\verb|dluna@fourier.ujf-grenoble.fr|
\end{document}

%% file: lunadiagrams.tex
%%%%%%%%%%% saveboxes

\newsavebox{\aone}
\newsavebox{\aprime}
\newsavebox{\GreyCircle}
\newsavebox{\segm}
\newsavebox{\susp}
\newsavebox{\shortsusp}
\newsavebox{\bifurc}
\newsavebox{\longbifurc}
\newsavebox{\dthree}
\newsavebox{\dm}
\newsavebox{\shortdm}
\newsavebox{\atwo}
\newsavebox{\tosw}
\newsavebox{\tose}
\newsavebox{\tonw}
\newsavebox{\tone}
\newsavebox{\toe}
\newsavebox{\tow}
\newsavebox{\longam}
\newsavebox{\mediumam}
\newsavebox{\shortam}
\newsavebox{\plusaaoneone}
\newsavebox{\plusdm}
\newsavebox{\vsegm}

\newsavebox{\rightbisegm}
\newsavebox{\bsecondtwo}
\newsavebox{\shortbprimem}
\newsavebox{\btwo}
\newsavebox{\shortbm}
\newsavebox{\leftbisegm}
\newsavebox{\shortcm}
\newsavebox{\cprimetwo}
\newsavebox{\shortbsecondm}
\newsavebox{\shortcsecondm}
\newsavebox{\bthirdthree}
\newsavebox{\ffour}
\newsavebox{\lefttrisegm}
\newsavebox{\gtwo}
\newsavebox{\gprimetwo}
\newsavebox{\gsecondtwo}

\newsavebox{\GreyCircleTwo}
\newsavebox{\plusbm}
\newsavebox{\plusbprimem}
\newsavebox{\pluscsecondm}

\newsavebox{\dynkinf}

\setlength{\unitlength}{700sp}

\savebox{\aone}{\begin{picture}(600,1800)\put(300,900){\circle*{150}}\put(300,1500){\circle{600}}\put(300,300){\circle{600}}\end{picture}}
\savebox{\aprime}{\begin{picture}(600,900)\put(300,900){\circle*{150}}\put(300,300){\circle{600}}\end{picture}}
\savebox{\GreyCircle}(600,600){\begin{picture}(600,600)(300,300)\put(600,600){\circle{600}}\put(30,25){\multiput(500,350)(150,0){2}{\circle*{70}}\multiput(425,425)(150,0){3}{\circle*{70}}\multiput(350,500)(150,0){4}{\circle*{70}}\multiput(425,575)(150,0){3}{\circle*{70}}\multiput(350,650)(150,0){4}{\circle*{70}}\multiput(425,725)(150,0){3}{\circle*{70}}\multiput(500,800)(150,0){2}{\circle*{70}}}\end{picture}}
\savebox{\segm}{\begin{picture}(1800,0)\multiput(0,0)(1800,0){2}{\circle*{150}}\thicklines\put(0,0){\line(1,0){1800}}\end{picture}}
\savebox{\susp}{\begin{picture}(3600,0)\multiput(0,0)(3600,0){2}{\circle*{150}}\thicklines\multiput(0,0)(2500,0){2}{\line(1,0){1100}}\multiput(1300,0)(400,0){3}{\line(1,0){200}}\end{picture}}
\savebox{\shortsusp}{\begin{picture}(1800,0)\multiput(0,0)(1800,0){2}{\circle*{150}}\thicklines\multiput(0,0)(400,0){5}{\line(1,0){200}}\end{picture}}
\savebox{\bifurc}{\begin{picture}(1200,2400)\multiput(1200,0)(0,2400){2}{\circle*{150}}\thicklines\put(0,1200){\line(1,1){1200}}\put(0,1200){\line(1,-1){1200}}\end{picture}}
\savebox{\longbifurc}{\begin{picture}(1800,3600)\multiput(1800,0)(0,3600){2}{\circle*{150}}\thicklines\put(0,1800){\line(1,1){1800}}\put(0,1800){\line(1,-1){1800}}\end{picture}}
\savebox{\dthree}{\begin{picture}(3600,600)\put(1500,0){\usebox{\GreyCircle}}\multiput(0,300)(1800,0){2}{\usebox{\segm}}\end{picture}}
\savebox{\dm}{\begin{picture}(8700,2400)\put(0,900){\usebox{\GreyCircle}}\multiput(300,1200)(5400,0){2}{\usebox{\segm}}\put(2100,1200){\usebox{\susp}}\put(7500,0){\usebox{\bifurc}}\end{picture}}
\savebox{\shortdm}{\begin{picture}(6900,2400)\put(0,900){\usebox{\GreyCircle}}\multiput(300,1200)(3600,0){2}{\usebox{\segm}}\put(2100,1200){\usebox{\shortsusp}}\put(5700,0){\usebox{\bifurc}}\end{picture}}
\savebox{\atwo}{\begin{picture}(2400,600)\put(300,300){\usebox{\segm}}\multiput(300,300)(1800,0){2}{\circle{600}}\multiput(300,300)(25,25){13}{\circle*{70}}\multiput(600,600)(300,0){4}{\multiput(0,0)(25,-25){7}{\circle*{70}}}\multiput(750,450)(300,0){4}{\multiput(0,0)(25,25){7}{\circle*{70}}}\multiput(2100,300)(-25,25){13}{\circle*{70}}\end{picture}}
\savebox{\tosw}{\begin{picture}(300,300)\multiput(-550,-550)(10,30){16}{\line(0,1){30}}\multiput(-550,-550)(30,10){16}{\line(1,0){30}}\end{picture}}
\savebox{\tose}{\begin{picture}(300,300)\multiput(550,-550)(-10,30){16}{\line(0,1){30}}\multiput(550,-550)(-30,10){16}{\line(-1,0){30}}\end{picture}}
\savebox{\tonw}{\begin{picture}(300,300)\multiput(-550,550)(10,-30){16}{\line(0,-1){30}}\multiput(-550,550)(30,-10){16}{\line(1,0){30}}\end{picture}}
\savebox{\tone}{\begin{picture}(300,300)\multiput(550,550)(-10,-30){16}{\line(0,-1){30}}\multiput(550,550)(-30,-10){16}{\line(-1,0){30}}\end{picture}}
\savebox{\toe}{\begin{picture}(400,400)\multiput(0,0)(20,10){21}{\line(1,0){20}}\multiput(0,400)(20,-10){21}{\line(1,0){20}}\end{picture}}
\savebox{\tow}{\begin{picture}(400,400)\multiput(400,0)(-20,10){21}{\line(-1,0){20}}\multiput(400,400)(-20,-10){21}{\line(-1,0){20}}\end{picture}}
\savebox{\longam}{\begin{picture}(7800,600)\put(300,300){\circle*{150}}\put(2100,300){\circle*{150}}\put(5700,300){\circle*{150}}\put(7500,300){\circle*{150}}\put(300,300){\circle{600}}\put(7500,300){\circle{600}}\multiput(300,300)(25,25){13}{\circle*{70}}\multiput(600,600)(3900,0){2}{\multiput(0,0)(300,0){9}{\multiput(0,0)(25,-25){7}{\circle*{70}}}\multiput(150,-150)(300,0){9}{\multiput(0,0)(25,25){7}{\circle*{70}}}}\multiput(7500,300)(-25,25){13}{\circle*{70}}\thicklines\put(300,300){\line( 1, 0){2250}}\put(7500,300){\line(-1, 0){2250}}\multiput(2850,300)(600,0){4}{\line( 1, 0){300}}\end{picture}}
\savebox{\mediumam}{\begin{picture}(6000,600)\multiput(300,300)(3600,0){2}{\usebox{\segm}}\put(2100,300){\usebox{\shortsusp}}\multiput(300,300)(1800,0){4}{\circle*{150}}\multiput(300,300)(5400,0){2}{\circle{600}}\multiput(300,300)(25,25){13}{\circle*{70}}\multiput(600,600)(3000,0){2}{\multiput(0,0)(300,0){6}{\multiput(0,0)(25,-25){7}{\circle*{70}}}\multiput(150,-150)(300,0){6}{\multiput(0,0)(25,25){7}{\circle*{70}}}}\multiput(5700,300)(-25,25){13}{\circle*{70}}\multiput(2700,600)(300,0){3}{\circle*{70}}\end{picture}}
\savebox{\shortam}{\begin{picture}(4200,600)\put(300,300){\usebox{\susp}}\multiput(300,300)(3600,0){2}{\circle{600}}\multiput(300,300)(25,25){13}{\circle*{70}}\multiput(600,600)(2100,0){2}{\multiput(0,0)(300,0){3}{\multiput(0,0)(25,-25){7}{\circle*{70}}}\multiput(150,-150)(300,0){3}{\multiput(0,0)(25,25){7}{\circle*{70}}}}\multiput(3900,300)(-25,25){13}{\circle*{70}}\multiput(1800,600)(300,0){3}{\circle*{70}}\end{picture}}
\savebox{\plusaaoneone}{\begin{picture}(1950,3000)\put(0,300){\usebox{\bifurc}}\multiput(1200,300)(0,2400){2}{\circle{600}}\multiput(1500,300)(0,2400){2}{\line(1,0){450}}\put(1950,300){\line(0,1){2400}}\end{picture}}
%\savebox{\plusdm}{\begin{picture}(6600,2400)\put(0,1200){\usebox{\segm}}\put(1800,1200){\usebox{\susp}}\put(5400,0){\usebox{\bifurc}}\put(1500,900){\usebox{\GreyCircle}}\end{picture}}
\savebox{\plusdm}{\begin{picture}(8400,2400)\put(0,1200){\usebox{\segm}}\put(1500,0){\usebox{\shortdm}}\end{picture}}
\savebox{\vsegm}{\begin{picture}(0,1800)\multiput(0,0)(0,1800){2}{\circle*{150}}\thicklines\put(0,0){\line(0,1){1800}}\end{picture}}

\savebox{\GreyCircleTwo}{\begin{picture}(600,600)(300,300)\put(600,600){\circle{600}}\put(30,25){\multiput(500,350)(150,0){2}{\circle*{70}}\multiput(425,425)(150,0){3}{\circle*{70}}\multiput(350,500)(150,0){4}{\circle*{70}}\multiput(425,575)(150,0){3}{\circle*{70}}\multiput(350,650)(150,0){4}{\circle*{70}}\multiput(425,725)(150,0){3}{\circle*{70}}\multiput(500,800)(150,0){2}{\circle*{70}}}\put(420,1020){\tiny2}\end{picture}}

\savebox{\rightbisegm}{\begin{picture}(1800,0)\multiput(0,0)(1800,0){2}{\circle*{150}}\thicklines\multiput(0,-60)(0,150){2}{\line(1,0){1800}}\multiput(1050,0)(-25,25){10}{\circle*{50}}\multiput(1050,0)(-25,-25){10}{\circle*{50}}\end{picture}}

%\savebox{\bsecondtwo}{\begin{picture}(2400,750)\put(300,300){\usebox{\rightbisegm}}\multiput(300,300)(1800,0){2}{\circle{600}}\multiput(300,300)(25,25){18}{\circle*{70}}\multiput(750,750)(300,0){3}{\multiput(0,0)(25,-25){7}{\circle*{70}}}\multiput(900,600)(300,0){3}{\multiput(0,0)(25,25){7}{\circle*{70}}}\multiput(2100,300)(-25,25){18}{\circle*{70}}\end{picture}}
\savebox{\bsecondtwo}{\begin{picture}(2400,600)\put(300,300){\usebox{\rightbisegm}}\put(0,0){\usebox{\GreyCircle}}\put(2100,300){\circle{600}}\end{picture}}

%\savebox{\shortbprimem}{\begin{picture}(7500,750)\multiput(300,300)(3600,0){2}{\usebox{\segm}}\put(2100,300){\usebox{\shortsusp}}\put(5700,300){\usebox{\rightbisegm}}\put(0,0){\usebox{\GreyCircle}}\multiput(300,300)(25,25){18}{\circle*{70}}\multiput(750,750)(300,0){6}{\multiput(0,0)(25,-25){7}{\circle*{70}}}\multiput(900,600)(300,0){5}{\multiput(0,0)(25,25){7}{\circle*{70}}}\multiput(7050,750)(-300,0){12}{\multiput(0,0)(-25,-25){7}{\circle*{70}}}\multiput(6900,600)(-300,0){11}{\multiput(0,0)(-25,25){7}{\circle*{70}}}\multiput(7500,300)(-25,25){18}{\circle*{70}}\multiput(2700,600)(300,0){3}{\circle*{70}}\end{picture}}
\savebox{\shortbprimem}{\begin{picture}(7500,1200)\multiput(300,300)(3600,0){2}{\usebox{\segm}}\put(2100,300){\usebox{\shortsusp}}\put(5700,300){\usebox{\rightbisegm}}\put(0,0){\usebox{\GreyCircleTwo}}\end{picture}}

%\savebox{\btwo}{\begin{picture}(2100,750)\put(300,300){\usebox{\rightbisegm}}\put(300,300){\circle{600}}\multiput(300,300)(25,25){18}{\circle*{70}}\multiput(750,750)(300,0){3}{\multiput(0,0)(25,-25){7}{\circle*{70}}}\multiput(900,600)(300,0){3}{\multiput(0,0)(25,25){7}{\circle*{70}}}\multiput(2100,300)(-25,25){18}{\circle*{70}}\end{picture}}
\savebox{\btwo}{\begin{picture}(2100,600)\put(300,300){\usebox{\rightbisegm}}\put(0,0){\usebox{\GreyCircle}}\end{picture}}

%\savebox{\shortbm}{\begin{picture}(7500,750)\multiput(300,300)(3600,0){2}{\usebox{\segm}}\put(2100,300){\usebox{\shortsusp}}\put(5700,300){\usebox{\rightbisegm}}\put(300,300){\circle{600}}\multiput(300,300)(25,25){18}{\circle*{70}}\multiput(750,750)(300,0){6}{\multiput(0,0)(25,-25){7}{\circle*{70}}}\multiput(900,600)(300,0){5}{\multiput(0,0)(25,25){7}{\circle*{70}}}\multiput(7050,750)(-300,0){12}{\multiput(0,0)(-25,-25){7}{\circle*{70}}}\multiput(6900,600)(-300,0){11}{\multiput(0,0)(-25,25){7}{\circle*{70}}}\multiput(7500,300)(-25,25){18}{\circle*{70}}\multiput(2700,600)(300,0){3}{\circle*{70}}\end{picture}}
\savebox{\shortbm}{\begin{picture}(7500,600)\multiput(300,300)(3600,0){2}{\usebox{\segm}}\put(2100,300){\usebox{\shortsusp}}\put(5700,300){\usebox{\rightbisegm}}\put(0,0){\usebox{\GreyCircle}}\end{picture}}

\savebox{\leftbisegm}{\begin{picture}(1800,0)\multiput(0,0)(1800,0){2}{\circle*{150}}\thicklines\multiput(0,-60)(0,150){2}{\line(1,0){1800}}\multiput(750,0)(25,25){10}{\circle*{50}}\multiput(750,0)(25,-25){10}{\circle*{50}}\end{picture}}

%\savebox{\shortcm}{\begin{picture}(9000,750)\multiput(0,300)(1800,0){2}{\usebox{\segm}}\put(3600,300){\usebox{\shortsusp}}\put(5400,300){\usebox{\segm}}\put(7200,300){\usebox{\leftbisegm}}\put(1800,300){\circle{600}}\multiput(0,300)(25,25){18}{\circle*{70}}\multiput(450,750)(300,0){12}{\multiput(0,0)(25,-25){7}{\circle*{70}}}\multiput(600,600)(300,0){11}{\multiput(0,0)(25,25){7}{\circle*{70}}}\multiput(8550,750)(-300,0){12}{\multiput(0,0)(-25,-25){7}{\circle*{70}}}\multiput(8400,600)(-300,0){11}{\multiput(0,0)(-25,25){7}{\circle*{70}}}\multiput(9000,300)(-25,25){18}{\circle*{70}}\multiput(4200,600)(300,0){3}{\circle*{70}}\end{picture}}
\savebox{\shortcm}{\begin{picture}(9000,600)\multiput(0,300)(1800,0){2}{\usebox{\segm}}\put(3600,300){\usebox{\shortsusp}}\put(5400,300){\usebox{\segm}}\put(7200,300){\usebox{\leftbisegm}}\put(1500,0){\usebox{\GreyCircle}}\end{picture}}

%\savebox{\cprimetwo}{\begin{picture}(2100,750)\put(0,300){\usebox{\leftbisegm}}\put(1500,0){\usebox{\GreyCircle}}\put(-300,0){\multiput(300,300)(25,25){18}{\circle*{70}}\multiput(750,750)(300,0){3}{\multiput(0,0)(25,-25){7}{\circle*{70}}}\multiput(900,600)(300,0){3}{\multiput(0,0)(25,25){7}{\circle*{70}}}\multiput(2100,300)(-25,25){18}{\circle*{70}}}\end{picture}}
\savebox{\cprimetwo}{\begin{picture}(2100,1200)\put(0,300){\usebox{\leftbisegm}}\put(1500,0){\usebox{\GreyCircleTwo}}\end{picture}}

%\savebox{\shortbsecondm}{\begin{picture}(7800,750)\multiput(300,300)(3600,0){2}{\usebox{\segm}}\put(2100,300){\usebox{\shortsusp}}\put(5700,300){\usebox{\rightbisegm}}\multiput(300,300)(7200,0){2}{\circle{600}}\multiput(300,300)(25,25){18}{\circle*{70}}\multiput(750,750)(300,0){6}{\multiput(0,0)(25,-25){7}{\circle*{70}}}\multiput(900,600)(300,0){5}{\multiput(0,0)(25,25){7}{\circle*{70}}}\multiput(7050,750)(-300,0){12}{\multiput(0,0)(-25,-25){7}{\circle*{70}}}\multiput(6900,600)(-300,0){11}{\multiput(0,0)(-25,25){7}{\circle*{70}}}\multiput(7500,300)(-25,25){18}{\circle*{70}}\multiput(2700,600)(300,0){3}{\circle*{70}}\end{picture}}
\savebox{\shortbsecondm}{\begin{picture}(7800,600)\multiput(300,300)(3600,0){2}{\usebox{\segm}}\put(2100,300){\usebox{\shortsusp}}\put(5700,300){\usebox{\rightbisegm}}\put(0,0){\usebox{\GreyCircle}}\put(7500,300){\circle{600}}\end{picture}}

%\savebox{\shortcsecondm}{\begin{picture}(9300,750)\put(300,300){\circle{600}}\put(300,0){\usebox{\shortcm}}\end{picture}}
\savebox{\shortcsecondm}{\begin{picture}(9300,600)\put(300,300){\circle{600}}\put(300,0){\usebox{\shortcm}}\end{picture}}

%\savebox{\bthirdthree}{\begin{picture}(3900,750)\put(0,300){\usebox{\segm}}\put(1800,300){\usebox{\rightbisegm}}\put(3600,300){\circle{600}}\multiput(0,300)(25,25){18}{\circle*{70}}\multiput(450,750)(300,0){9}{\multiput(0,0)(25,-25){7}{\circle*{70}}}\multiput(600,600)(300,0){9}{\multiput(0,0)(25,25){7}{\circle*{70}}}\multiput(3600,300)(-25,25){18}{\circle*{70}}\end{picture}}
\savebox{\bthirdthree}{\begin{picture}(3900,600)\put(0,300){\usebox{\segm}}\put(1800,300){\usebox{\rightbisegm}}\put(3300,0){\usebox{\GreyCircle}}\end{picture}}

%\savebox{\ffour}{\begin{picture}(5700,750)\multiput(0,300)(3600,0){2}{\usebox{\segm}}\put(1800,300){\usebox{\rightbisegm}}\put(5400,300){\circle{600}}\multiput(0,300)(25,25){18}{\circle*{70}}\multiput(450,750)(300,0){15}{\multiput(0,0)(25,-25){7}{\circle*{70}}}\multiput(600,600)(300,0){15}{\multiput(0,0)(25,25){7}{\circle*{70}}}\multiput(5400,300)(-25,25){18}{\circle*{70}}\end{picture}}
\savebox{\ffour}{\begin{picture}(5700,600)\multiput(0,300)(3600,0){2}{\usebox{\segm}}\put(1800,300){\usebox{\rightbisegm}}\put(5100,0){\usebox{\GreyCircle}}\end{picture}}

\savebox{\lefttrisegm}{\begin{picture}(1800,0)\multiput(0,0)(1800,0){2}{\circle*{150}}\thicklines\multiput(0,-120)(0,135){3}{\line(1,0){1800}}\multiput(750,0)(25,25){10}{\circle*{50}}\multiput(750,0)(25,-25){10}{\circle*{50}}\end{picture}}

%\savebox{\gtwo}{\begin{picture}(2100,750)\put(300,300){\usebox{\lefttrisegm}}\put(300,300){\circle{600}}\multiput(300,300)(25,25){18}{\circle*{70}}\multiput(750,750)(300,0){3}{\multiput(0,0)(25,-25){7}{\circle*{70}}}\multiput(900,600)(300,0){3}{\multiput(0,0)(25,25){7}{\circle*{70}}}\multiput(2100,300)(-25,25){18}{\circle*{70}}\end{picture}}
\savebox{\gtwo}{\begin{picture}(2100,600)\put(300,300){\usebox{\lefttrisegm}}\put(0,0){\usebox{\GreyCircle}}\end{picture}}

%\savebox{\gprimetwo}{\begin{picture}(2100,750)\put(300,300){\usebox{\lefttrisegm}}\put(0,0){\usebox{\GreyCircle}}\multiput(300,300)(25,25){18}{\circle*{70}}\multiput(750,750)(300,0){3}{\multiput(0,0)(25,-25){7}{\circle*{70}}}\multiput(900,600)(300,0){3}{\multiput(0,0)(25,25){7}{\circle*{70}}}\multiput(2100,300)(-25,25){18}{\circle*{70}}\end{picture}}
\savebox{\gprimetwo}{\begin{picture}(2100,1200)\put(300,300){\usebox{\lefttrisegm}}\put(0,0){\usebox{\GreyCircleTwo}}\end{picture}}

%\savebox{\gsecondtwo}{\begin{picture}(2400,750)\put(300,300){\usebox{\lefttrisegm}}\multiput(300,300)(1800,0){2}{\circle{600}}\multiput(300,300)(25,25){18}{\circle*{70}}\multiput(750,750)(300,0){3}{\multiput(0,0)(25,-25){7}{\circle*{70}}}\multiput(900,600)(300,0){3}{\multiput(0,0)(25,25){7}{\circle*{70}}}\multiput(2100,300)(-25,25){18}{\circle*{70}}\end{picture}}
\savebox{\gsecondtwo}{\begin{picture}(2400,600)\put(300,300){\usebox{\lefttrisegm}}\put(1800,0){\usebox{\GreyCircle}}\put(300,300){\circle{600}}\end{picture}}

\savebox{\plusbm}{\begin{picture}(9000,600)\put(0,300){\usebox{\segm}}\put(1500,0){\usebox{\shortbm}}\end{picture}}

\savebox{\plusbprimem}{\begin{picture}(9000,1200)\put(0,300){\usebox{\segm}}\put(1500,0){\usebox{\shortbprimem}}\end{picture}}

\savebox{\pluscsecondm}{\begin{picture}(9000,600)\put(0,0){\usebox{\shortcm}}\end{picture}}

\savebox{\dynkinf}{\begin{picture}(5400,0)\multiput(0,0)(3600,0){2}{\usebox{\segm}}\put(1800,0){\usebox{\rightbisegm}}\end{picture}}

%%%%%%%%%% minimal cases with a 1-projective element in the support of
%%%%%%%%%% a non-simple spherical root

\newcommand{\projcsecondm}{\begin{picture}(9300,1800)\put(0,0){\usebox{\aone}}\put(300,600){\usebox{\pluscsecondm}}\end{picture}}

\newcommand{\projgsecondtwo}{\begin{picture}(2400,1800)\put(0,0){\usebox{\aone}}\put(300,900){\usebox{\lefttrisegm}}\put(1800,600){\usebox{\GreyCircle}}\end{picture}}

\newcommand{\projffourone}{\begin{picture}(6000,1800)\multiput(300,900)(3600,0){2}{\usebox{\segm}}\put(2100,900){\usebox{\rightbisegm}}\multiput(0,600)(3600,0){2}{\usebox{\GreyCircle}}\put(5400,0){\usebox{\aone}}\end{picture}}

\newcommand{\projffourtwo}{\begin{picture}(6000,2250)\multiput(300,900)(3600,0){2}{\usebox{\segm}}\put(2100,900){\usebox{\rightbisegm}}\multiput(0,0)(3600,0){2}{\usebox{\aone}}\put(5400,0){\usebox{\aone}}\put(1800,600){\usebox{\GreyCircle}}\multiput(300,1800)(5400,0){2}{\line(0,1){450}}\put(300,2250){\line(1,0){5400}}\put(700,1300){\usebox{\toe}}\end{picture}}

\newcommand{\projffourthree}{\begin{picture}(6000,1800)\multiput(300,900)(3600,0){2}{\usebox{\segm}}\put(2100,900){\usebox{\rightbisegm}}\put(0,600){\usebox{\atwo}}\put(1800,600){\usebox{\GreyCircle}}\multiput(3600,0)(1800,0){2}{\usebox{\aone}}\put(4900,1300){\usebox{\tow}}\end{picture}}

\newcommand{\projbcprime}{\begin{picture}(11400,1800)\put(0,600){\multiput(0,0)(1800,0){2}{\usebox{\atwo}}\put(7200,0){\usebox{\atwo}}\put(3900,300){\usebox{\susp}}\multiput(3900,300)(25,25){13}{\circle*{70}}\multiput(4200,600)(2400,0){2}{\multiput(0,0)(300,0){2}{\multiput(0,0)(25,-25){7}{\circle*{70}}}\multiput(150,-150)(300,0){2}{\multiput(0,0)(25,25){7}{\circle*{70}}}}\multiput(7500,300)(-25,25){13}{\circle*{70}}}\put(9000,600){\usebox{\btwo}}\put(10800,0){\usebox{\aone}}\end{picture}}

%%%%%%%%%% primitive cases

\newcommand{\diagrambbpp}{\begin{picture}(9600,3300)\multiput(0,0)(0,2700){2}{\multiput(300,300)(5400,0){2}{\usebox{\segm}}\put(7500,300){\usebox{\rightbisegm}}\put(2100,300){\usebox{\susp}}\multiput(300,300)(1800,0){2}{\circle{600}}\multiput(5700,300)(1800,0){3}{\circle{600}}}\multiput(300,600)(1800,0){2}{\line(0,1){2100}}\multiput(5700,600)(1800,0){3}{\line(0,1){2100}}\end{picture}}

\newcommand{\diagrambopplusq}{\begin{picture}(14700,1350)\multiput(300,900)(5400,0){2}{\usebox{\segm}}\put(2100,900){\usebox{\susp}}\multiput(0,0)(1800,0){2}{\usebox{\aprime}}\put(5400,0){\usebox{\aprime}}\put(7200,600){\usebox{\shortbprimem}}\end{picture}}

\newcommand{\diagramacastpplusbq}{\begin{picture}(18300,750)\multiput(0,0)(1800,0){2}{\usebox{\atwo}}\put(7200,0){\usebox{\atwo}}\put(9300,300){\usebox{\segm}}\put(3900,300){\usebox{\susp}}\multiput(3900,300)(25,25){13}{\circle*{70}}\multiput(4200,600)(2400,0){2}{\multiput(0,0)(300,0){2}{\multiput(0,0)(25,-25){7}{\circle*{70}}}\multiput(150,-150)(300,0){2}{\multiput(0,0)(25,25){7}{\circle*{70}}}}\multiput(7500,300)(-25,25){13}{\circle*{70}}\put(10800,0){\usebox{\shortbm}}\end{picture}}

\newcommand{\diagramacastpplusbprimeq}{\begin{picture}(18300,750)\multiput(0,0)(1800,0){2}{\usebox{\atwo}}\put(7200,0){\usebox{\atwo}}\put(9300,300){\usebox{\segm}}\put(3900,300){\usebox{\susp}}\multiput(3900,300)(25,25){13}{\circle*{70}}\multiput(4200,600)(2400,0){2}{\multiput(0,0)(300,0){2}{\multiput(0,0)(25,-25){7}{\circle*{70}}}\multiput(150,-150)(300,0){2}{\multiput(0,0)(25,25){7}{\circle*{70}}}}\multiput(7500,300)(-25,25){13}{\circle*{70}}\put(10800,0){\usebox{\shortbprimem}}\end{picture}}

\newcommand{\diagramccpp}{\begin{picture}(9600,3300)\multiput(0,0)(0,2700){2}{\multiput(300,300)(5400,0){2}{\usebox{\segm}}\put(7500,300){\usebox{\leftbisegm}}\put(2100,300){\usebox{\susp}}\multiput(300,300)(1800,0){2}{\circle{600}}\multiput(5700,300)(1800,0){3}{\circle{600}}}\multiput(300,600)(1800,0){2}{\line(0,1){2100}}\multiput(5700,600)(1800,0){3}{\line(0,1){2100}}\end{picture}}

\newcommand{\diagramcoastn}{\begin{picture}(9600,1800)\multiput(300,900)(5400,0){2}{\usebox{\segm}}\put(2100,900){\usebox{\susp}}\multiput(0,0)(1800,0){2}{\usebox{\aprime}}\multiput(5400,0)(1800,0){2}{\usebox{\aprime}}\put(7500,900){\usebox{\leftbisegm}}\put(9000,0){\usebox{\aone}}\put(8500,1300){\usebox{\tow}}\end{picture}}

\newcommand{\diagramcon}{\begin{picture}(9600,1200)\multiput(300,900)(5400,0){2}{\usebox{\segm}}\put(2100,900){\usebox{\susp}}\multiput(0,0)(1800,0){2}{\usebox{\aprime}}\multiput(5400,0)(1800,0){3}{\usebox{\aprime}}\put(7500,900){\usebox{\leftbisegm}}\end{picture}}

\newcommand{\diagramddpp}{\begin{picture}(9900,5700)\put(0,1800){\multiput(300,1200)(900,-2700){2}{\multiput(0,1200)(5400,0){2}{\usebox{\segm}}\put(1800,1200){\usebox{\susp}}\put(7200,0){\usebox{\bifurc}}\multiput(0,1200)(5400,0){2}{\multiput(0,0)(1800,0){2}{\circle{600}}}\multiput(8400,0)(0,2400){2}{\circle{600}}}\put(400,900){\multiput(0,1200)(5400,0){2}{\multiput(0,0)(1800,0){2}{\line(1,-3){700}}}\multiput(8400,0)(0,2400){2}{\line(1,-3){700}}}}\end{picture}}

\newcommand{\diagramdopplusq}{\begin{picture}(14100,2400)\multiput(0,300)(5400,0){2}{\put(0,0){\usebox{\aprime}}\put(300,900){\usebox{\segm}}}\put(1800,300){\usebox{\aprime}}\put(7200,0){\usebox{\shortdm}}\put(2100,1200){\usebox{\susp}}\end{picture}}

\newcommand{\diagramdon}{\begin{picture}(9000,3300)\multiput(0,1200)(5400,0){2}{\multiput(0,0)(1800,0){2}{\usebox{\aprime}}}\multiput(300,2100)(5400,0){2}{\usebox{\segm}}\put(2100,2100){\usebox{\susp}}\put(7500,900){\usebox{\bifurc}}\multiput(8400,0)(0,2400){2}{\usebox{\aprime}}\end{picture}}

\newcommand{\diagramdcneven}{\begin{picture}(14100,3300)\multiput(0,1800)(7200,0){2}{\usebox{\dthree}}\put(3600,2100){\usebox{\susp}}\put(10800,2100){\usebox{\segm}}\put(12600,900){\usebox{\bifurc}}\put(12300,1800){\usebox{\GreyCircle}}\put(13500,0){\usebox{\aone}}\put(13800,1500){\usebox{\tonw}}\end{picture}}

\newcommand{\diagramdcprimeneven}{\begin{picture}(14100,3300)\multiput(0,1800)(7200,0){2}{\usebox{\dthree}}\put(3600,2100){\usebox{\susp}}\put(10800,2100){\usebox{\segm}}\put(12600,900){\usebox{\bifurc}}\put(12300,1800){\usebox{\GreyCircle}}\put(13500,0){\usebox{\aprime}}\end{picture}}

\newcommand{\diagramdcnodd}{\begin{picture}(12300,3000)\put(0,1200){\usebox{\dthree}}\put(3600,1500){\usebox{\susp}}\put(7200,1200){\usebox{\dthree}}\put(10800,300){\usebox{\bifurc}}\multiput(12000,300)(0,2400){2}{\circle{600}}\thicklines\multiput(12000,300)(0,2000){2}{\line(0,1){400}}\multiput(12000,700)(-200,200){4}{\line(-1,0){200}}\multiput(11800,700)(-200,200){4}{\line(0,1){200}}\multiput(12000,2300)(-200,-200){4}{\line(-1,0){200}}\multiput(11800,2300)(-200,-200){4}{\line(0,-1){200}}\end{picture}}

\newcommand{\diagramdcastn}{\begin{picture}(10800,3000)\multiput(0,1200)(1800,0){2}{\usebox{\atwo}}\put(7200,1200){\usebox{\atwo}}\put(3900,1500){\usebox{\susp}}\put(9300,300){\usebox{\bifurc}}\multiput(10500,300)(0,2400){2}{\circle{600}}\multiput(3900,1500)(25,25){13}{\circle*{70}}\multiput(4200,1800)(2400,0){2}{\multiput(0,0)(300,0){2}{\multiput(0,0)(25,-25){7}{\circle*{70}}}\multiput(150,-150)(300,0){2}{\multiput(0,0)(25,25){7}{\circle*{70}}}}\multiput(7500,1500)(-25,25){13}{\circle*{70}}\thicklines\put(9300,1500){\line(0,1){400}}\multiput(9300,1900)(200,200){4}{\line(1,0){200}}\multiput(9500,1900)(200,200){4}{\line(0,1){200}}\put(10500,2700){\line(-1,0){400}}\put(9300,1500){\line(1,0){400}}\multiput(9700,1500)(200,-200){4}{\line(0,-1){200}}\multiput(9700,1300)(200,-200){4}{\line(1,0){200}}\put(10500,300){\line(0,1){400}}\end{picture}}

\newcommand{\diagrameepp}{\begin{picture}(12450,5550)\put(1050,0){\multiput(0,0)(-1050,3150){2}{\multiput(300,2100)(1800,0){4}{\usebox{\segm}}\put(3900,300){\usebox{\vsegm}}\multiput(7500,2100)(1800,0){2}{\usebox{\shortsusp}}\multiput(300,2100)(1800,0){7}{\circle{600}}\put(3900,300){\circle{600}}}\put(-100,300){\multiput(300,2100)(1800,0){5}{\line(-1,3){850}}\put(3900,300){\line(-1,3){850}}\multiput(9300,2100)(1800,0){2}{\multiput(0,-50)(-50,150){18}{\line(0,1){100}}}}}\end{picture}}

\newcommand{\diagrameon}{\begin{picture}(11400,3000)\multiput(300,2700)(1800,0){4}{\usebox{\segm}}\put(3900,900){\usebox{\vsegm}}\multiput(0,1800)(1800,0){7}{\usebox{\aprime}}\put(3600,0){\usebox{\aprime}}\multiput(7500,2700)(1800,0){2}{\usebox{\shortsusp}}\end{picture}}

\newcommand{\diagrameasix}{\begin{picture}(7800,3600)\multiput(300,2700)(1800,0){4}{\usebox{\segm}}\put(3900,900){\usebox{\vsegm}}\multiput(3600,0)(0,1800){2}{\usebox{\aprime}}\multiput(300,2700)(5400,0){2}{\multiput(0,0)(1800,0){2}{\circle{600}}}\multiput(300,3000)(7200,0){2}{\line(0,1){600}}\put(300,3600){\line(1,0){7200}}\multiput(2100,3000)(3600,0){2}{\line(0,1){300}}\put(2100,3300){\line(1,0){3600}}\end{picture}}

\newcommand{\diagramecseven}{\begin{picture}(9300,2100)\multiput(300,1800)(1800,0){5}{\usebox{\segm}}\put(3900,0){\usebox{\vsegm}}\multiput(3600,1500)(3600,0){2}{\usebox{\GreyCircle}}\multiput(0,900)(1800,0){2}{\usebox{\aprime}}\end{picture}}

\newcommand{\diagramefseven}{\begin{picture}(9600,2700)\multiput(300,1800)(1800,0){5}{\usebox{\segm}}\put(3900,0){\usebox{\vsegm}}\multiput(0,1500)(7200,0){2}{\usebox{\GreyCircle}}\put(9000,900){\usebox{\aone}}\put(8500,2200){\usebox{\tow}}\end{picture}}

\newcommand{\diagramefprimen}{\begin{picture}(11400,2100)\multiput(300,1800)(1800,0){5}{\usebox{\segm}}\put(3900,0){\usebox{\vsegm}}\multiput(0,1500)(7200,0){2}{\usebox{\GreyCircle}}\multiput(9000,900)(1800,0){2}{\usebox{\aprime}}\put(9300,1800){\usebox{\shortsusp}}\end{picture}}

\newcommand{\diagramecastn}{\begin{picture}(11400,2400)\multiput(300,2100)(1800,0){4}{\usebox{\segm}}\put(3900,300){\usebox{\vsegm}}\multiput(7500,2100)(1800,0){2}{\usebox{\shortsusp}}\multiput(0,1800)(1800,0){6}{\multiput(300,300)(1800,0){2}{\circle{600}}\multiput(300,300)(25,25){13}{\circle*{70}}\multiput(600,600)(300,0){4}{\multiput(0,0)(25,-25){7}{\circle*{70}}}\multiput(750,450)(300,0){4}{\multiput(0,0)(25,25){7}{\circle*{70}}}\multiput(2100,300)(-25,25){13}{\circle*{70}}}\put(3900,300){\circle{600}}\put(3600,0){\multiput(300,300)(25,25){13}{\circle*{70}}\multiput(600,600)(0,300){4}{\multiput(0,0)(-25,25){7}{\circle*{70}}}\multiput(450,750)(0,300){4}{\multiput(0,0)(25,25){7}{\circle*{70}}}\multiput(300,2100)(25,-25){13}{\circle*{70}}}\end{picture}}

\newcommand{\diagramefsix}{\begin{picture}(7800,2100)\multiput(300,1800)(1800,0){4}{\usebox{\segm}}\put(3900,0){\usebox{\vsegm}}\multiput(0,1500)(7200,0){2}{\usebox{\GreyCircle}}\end{picture}}

\newcommand{\diagramefsixplusatwo}{\begin{picture}(11400,2100)\multiput(300,1800)(1800,0){4}{\usebox{\segm}}\put(3900,0){\usebox{\vsegm}}\multiput(0,1500)(7200,0){2}{\usebox{\GreyCircle}}\put(7500,1800){\usebox{\segm}}\put(9000,1500){\usebox{\atwo}}\end{picture}}

\newcommand{\diagramacfiveplusatwo}{\begin{picture}(9300,2100)\put(0,1500){\usebox{\atwo}}\put(2100,0){\multiput(0,1800)(1800,0){4}{\usebox{\segm}}\put(1800,0){\usebox{\vsegm}}\multiput(1500,1500)(3600,0){2}{\usebox{\GreyCircle}}}\end{picture}}

\newcommand{\diagramfffourfour}{\begin{picture}(6000,3300)\multiput(0,0)(0,2700){2}{\multiput(300,300)(3600,0){2}{\usebox{\segm}}\put(2100,300){\usebox{\rightbisegm}}\multiput(300,300)(1800,0){4}{\circle{600}}}\multiput(300,600)(1800,0){4}{\line(0,1){2100}}\end{picture}}

\newcommand{\diagramfofour}{\begin{picture}(6000,900)\multiput(300,900)(3600,0){2}{\usebox{\segm}}\put(2100,900){\usebox{\rightbisegm}}\multiput(0,0)(1800,0){4}{\usebox{\aprime}}\end{picture}}

\newcommand{\diagramccpplusq}{\begin{picture}(19800,750)\put(0,0){\usebox{\dthree}}\put(3600,300){\usebox{\susp}}\put(7200,0){\usebox{\dthree}}\put(10800,0){\usebox{\shortcm}}\end{picture}}

\newcommand{\diagramccprimepplustwo}{\begin{picture}(12900,750)\put(0,0){\usebox{\dthree}}\put(3600,300){\usebox{\susp}}\put(7200,0){\usebox{\dthree}}\put(10800,0){\usebox{\cprimetwo}}\end{picture}}

\newcommand{\diagrambcastn}{\begin{picture}(11400,750)\multiput(0,0)(1800,0){2}{\usebox{\atwo}}\put(7200,0){\usebox{\atwo}}\put(9000,0){\usebox{\bsecondtwo}}\put(3900,300){\usebox{\susp}}\multiput(3900,300)(25,25){13}{\circle*{70}}\multiput(4200,600)(2400,0){2}{\multiput(0,0)(300,0){2}{\multiput(0,0)(25,-25){7}{\circle*{70}}}\multiput(150,-150)(300,0){2}{\multiput(0,0)(25,25){7}{\circle*{70}}}}\multiput(7500,300)(-25,25){13}{\circle*{70}}\end{picture}}

\newcommand{\diagrambcn}{\begin{picture}(11400,1800)\put(0,600){\multiput(0,0)(1800,0){2}{\usebox{\atwo}}\put(7200,0){\usebox{\atwo}}\put(3900,300){\usebox{\susp}}\multiput(3900,300)(25,25){13}{\circle*{70}}\multiput(4200,600)(2400,0){2}{\multiput(0,0)(300,0){2}{\multiput(0,0)(25,-25){7}{\circle*{70}}}\multiput(150,-150)(300,0){2}{\multiput(0,0)(25,25){7}{\circle*{70}}}}\multiput(7500,300)(-25,25){13}{\circle*{70}}}\put(9000,600){\usebox{\btwo}}\put(10800,0){\usebox{\aone}}\put(10300,1300){\usebox{\tow}}\end{picture}}

\newcommand{\diagrambcprimen}{\begin{picture}(11400,1350)\put(0,600){\multiput(0,0)(1800,0){2}{\usebox{\atwo}}\put(7200,0){\usebox{\atwo}}\put(3900,300){\usebox{\susp}}\multiput(3900,300)(25,25){13}{\circle*{70}}\multiput(4200,600)(2400,0){2}{\multiput(0,0)(300,0){2}{\multiput(0,0)(25,-25){7}{\circle*{70}}}\multiput(150,-150)(300,0){2}{\multiput(0,0)(25,25){7}{\circle*{70}}}}\multiput(7500,300)(-25,25){13}{\circle*{70}}}\put(9000,600){\usebox{\btwo}}\put(10800,0){\usebox{\aprime}}\end{picture}}

\newcommand{\diagramdsn}{\begin{picture}(7200,3300)\multiput(300,1500)(3600,0){2}{\usebox{\segm}}\put(2100,1500){\usebox{\shortsusp}}\put(5700,300){\usebox{\bifurc}}\put(300,1500){\circle{600}}\multiput(6900,300)(0,2400){2}{\circle{600}}\multiput(300,1500)(25,25){13}{\circle*{70}}\put(600,1800){\multiput(0,0)(300,0){6}{\multiput(0,0)(25,-25){7}{\circle*{70}}}\multiput(150,-150)(300,0){6}{\multiput(0,0)(25,25){7}{\circle*{70}}}}\put(3600,1800){\multiput(0,0)(300,0){7}{\multiput(0,0)(25,-25){7}{\circle*{70}}}\multiput(150,-150)(300,0){6}{\multiput(0,0)(25,25){7}{\circle*{70}}}}\multiput(300,1500)(25,-25){13}{\circle*{70}}\put(600,1200){\multiput(0,0)(300,0){6}{\multiput(0,0)(25,25){7}{\circle*{70}}}\multiput(150,150)(300,0){6}{\multiput(0,0)(25,-25){7}{\circle*{70}}}}\put(3600,1200){\multiput(0,0)(300,0){7}{\multiput(0,0)(25,25){7}{\circle*{70}}}\multiput(150,150)(300,0){6}{\multiput(0,0)(25,-25){7}{\circle*{70}}}}\multiput(2700,1200)(0,600){2}{\multiput(0,0)(300,0){3}{\circle*{70}}}\thicklines\put(6900,2700){\line(-1,0){400}}\multiput(6500,2700)(-200,-200){5}{\line(0,-1){200}}\multiput(6500,2500)(-200,-200){4}{\line(-1,0){200}}\multiput(5700,1700)(-30,-10){5}{\line(-1,0){30}}\put(6900,300){\line(-1,0){400}}\multiput(6500,300)(-200,200){5}{\line(0,1){200}}\multiput(6500,500)(-200,200){4}{\line(-1,0){200}}\multiput(5700,1300)(-30,10){5}{\line(-1,0){30}}\end{picture}}

\newcommand{\diagramdsastfour}{\begin{picture}(3600,3000)\put(300,1500){\usebox{\segm}}\put(2100,300){\usebox{\bifurc}}\put(300,1500){\circle{600}}\multiput(3300,300)(0,2400){2}{\circle{600}}\multiput(300,1500)(25,25){13}{\circle*{70}}\multiput(600,1800)(300,0){5}{\multiput(0,0)(25,-25){7}{\circle*{70}}}\multiput(750,1650)(300,0){4}{\multiput(0,0)(25,25){7}{\circle*{70}}}\multiput(300,1500)(25,-25){13}{\circle*{70}}\multiput(600,1200)(300,0){5}{\multiput(0,0)(25,25){7}{\circle*{70}}}\multiput(750,1350)(300,0){4}{\multiput(0,0)(25,-25){7}{\circle*{70}}}\thicklines\put(3300,2700){\line(-1,0){400}}\multiput(2900,2700)(-200,-200){5}{\line(0,-1){200}}\multiput(2900,2500)(-200,-200){4}{\line(-1,0){200}}\multiput(2100,1700)(-30,-10){5}{\line(-1,0){30}}\put(3300,300){\line(-1,0){400}}\multiput(2900,300)(-200,200){5}{\line(0,1){200}}\multiput(2900,500)(-200,200){4}{\line(-1,0){200}}\multiput(2100,1300)(-30,10){5}{\line(-1,0){30}}\multiput(3300,300)(0,2000){2}{\line(0,1){400}}\multiput(3300,700)(-200,200){4}{\line(-1,0){200}}\multiput(3100,700)(-200,200){4}{\line(0,1){200}}\multiput(3300,2300)(-200,-200){4}{\line(-1,0){200}}\multiput(3100,2300)(-200,-200){4}{\line(0,-1){200}}\end{picture}}

\newcommand{\diagramaatwotwoplusatwo}{\begin{picture}(7800,3000)
\put(0,-600){\multiput(300,2700)(1800,0){4}{\usebox{\segm}}\put(3900,900){\usebox{\vsegm}}\multiput(300,2700)(5400,0){2}{\multiput(0,0)(1800,0){2}{\circle{600}}}\multiput(300,3000)(7200,0){2}{\line(0,1){600}}\put(300,3600){\line(1,0){7200}}\multiput(2100,3000)(3600,0){2}{\line(0,1){300}}\put(2100,3300){\line(1,0){3600}}}
\put(3600,0){\multiput(300,300)(0,1800){2}{\circle{600}}\multiput(300,300)(25,25){13}{\circle*{70}}\multiput(600,600)(0,300){4}{\multiput(0,0)(-25,25){7}{\circle*{70}}}\multiput(450,750)(0,300){4}{\multiput(0,0)(25,25){7}{\circle*{70}}}\multiput(300,2100)(25,-25){13}{\circle*{70}}}
\end{picture}}

\newcommand{\diagramaotwoplusatwo}{\begin{picture}(6000,1200)\multiput(300,900)(3600,0){2}{\usebox{\segm}}\put(2100,900){\usebox{\rightbisegm}}\multiput(3600,0)(1800,0){2}{\usebox{\aprime}}\put(0,600){\usebox{\atwo}}\end{picture}}

\newcommand{\diagramfcastfour}{\begin{picture}(6000,750)\multiput(0,0)(3600,0){2}{\usebox{\atwo}}\put(1800,0){\usebox{\btwo}}\end{picture}}

\newcommand{\diagrambxthreethree}{\begin{picture}(10950,3300)
\put(0,900){
\multiput(300,900)(1800,0){2}{\usebox{\segm}}
\put(6600,900){\usebox{\segm}}
\put(8400,900){\usebox{\rightbisegm}}
\multiput(0,0)(6300,0){2}{\multiput(0,0)(1800,0){3}{\usebox{\aone}}}
}
\multiput(300,2700)(3600,0){2}{\line(0,1){600}}
\put(8400,2700){\line(0,1){600}}
\put(300,3300){\line(1,0){8100}}
\multiput(2100,2700)(4500,0){2}{\line(0,1){300}}
\put(2100,3000){\line(1,0){1700}}
\put(4000,3000){\line(1,0){2600}}
\put(300,0){\line(0,1){900}}
\put(300,0){\line(1,0){10650}}
\put(10950,0){\line(0,1){2400}}
\put(10950,2400){\line(-1,0){450}}
\multiput(2100,300)(6300,0){2}{\line(0,1){600}}
\put(2100,300){\line(1,0){6300}}
\multiput(3900,600)(6300,0){2}{\line(0,1){300}}
\put(3900,600){\line(1,0){4400}}
\put(8500,600){\line(1,0){1700}}
\put(650,2200){\usebox{\toe}}
\multiput(3150,2200)(6300,0){2}{\usebox{\tow}}
\multiput(7050,2200)(1800,0){2}{\usebox{\toe}}
\end{picture}}

\newcommand{\diagramacastn}{\begin{picture}(11400,600)\multiput(0,0)(7200,0){2}{\multiput(0,0)(1800,0){2}{\usebox{\atwo}}}\put(3900,300){\usebox{\susp}}\multiput(3900,300)(25,25){13}{\circle*{70}}\multiput(4200,600)(2400,0){2}{\multiput(0,0)(300,0){2}{\multiput(0,0)(25,-25){7}{\circle*{70}}}\multiput(150,-150)(300,0){2}{\multiput(0,0)(25,25){7}{\circle*{70}}}}\multiput(7500,300)(-25,25){13}{\circle*{70}}\end{picture}}

%% file: figures.tex
\newcommand{\stronebthreeoneone}{
\begin{picture}(34800,30150)(0,-30150)
\multiput(0,0)(6000,-4800){2}{
\put(3000,-2700){\vector(0,-1){5700}}
\put(6900,-900){\vector(1,0){4200}}
\multiput(0,-1800)(1800,0){3}{\usebox{\aone}}
\multiput(1300,-500)(1800,0){2}{\usebox{\tow}}
\put(12000,-1800){\usebox{\aone}}
\put(13800,-1200){\usebox{\GreyCircle}}
\put(15900,-900){\circle{600}}
\put(0,-10500){\usebox{\atwo}}
\put(3600,-11100){\usebox{\aone}}
\put(3100,-9800){\usebox{\tow}}
\put(12000,-10500){\usebox{\GreyCircle}}
\put(15900,-10200){\circle{600}}
\multiput(0,0)(0,-9300){2}{
\multiput(0,0)(12000,0){2}{
\put(300,-900){\usebox{\dynkinf}}
}}}
\multiput(0,0)(0,-9300){2}{
\put(24900,-5700){\vector(1,0){3000}}
\put(29100,-5700){\usebox{\dynkinf}}
\put(34500,-5700){\circle{600}}
\multiput(0,0)(12000,0){2}{
\put(5100,-2400){\vector(1,-1){1800}}
\put(5400,-1800){\usebox{\aone}}
\put(11700,-5700){\circle{600}}
}}
\put(18300,-24000){\usebox{\dynkinf}}
\multiput(18000,-24900)(5400,0){2}{\usebox{\aone}}
\put(19800,-24300){\usebox{\GreyCircle}}
\multiput(18300,-23100)(5400,0){2}{\line(0,1){450}}
\put(18300,-22650){\line(1,0){5400}}
\put(18700,-23600){\usebox{\toe}}
\put(12300,-19500){\usebox{\dynkinf}}
\put(12900,-15000){\vector(1,0){4200}}
\multiput(31800,-7500)(-10800,0){2}{\vector(0,-1){5700}}
\put(21000,-22200){\vector(0,1){5400}}
\put(15000,-17700){\line(0,1){2550}}
\put(15000,-14850){\vector(0,1){2850}}
\put(15000,-2700){\line(0,-1){2850}}
\put(15000,-5850){\vector(0,-1){2550}}
\put(6900,-10200){\line(1,0){1950}}
\put(9150,-10200){\vector(1,0){1950}}
\put(28800,-6600){\usebox{\aone}}
\put(30600,-6000){\usebox{\GreyCircle}}
\put(28800,-15300){\usebox{\GreyCircle}}
\put(12000,-19800){\usebox{\atwo}}
\put(13800,-19800){\usebox{\GreyCircle}}
\multiput(15600,-20400)(1800,0){2}{\usebox{\aone}}
\multiput(15100,-19100)(1800,0){2}{\usebox{\tow}}
\put(21900,-24000){\circle{600}}
\put(22900,-23600){\usebox{\tow}}
\put(24000,-20400){
\put(300,900){\usebox{\dynkinf}}
\multiput(0,600)(3600,0){2}{\usebox{\GreyCircle}}
\put(5700,900){\circle{600}}
}
\put(24000,-30150){
\put(300,900){\usebox{\dynkinf}}
\multiput(0,0)(3600,0){2}{\usebox{\aone}}
\put(5400,0){\usebox{\aone}}
\multiput(300,1800)(5400,0){2}{\line(0,1){450}}
\put(300,2250){\line(1,0){5400}}
\put(700,1300){\usebox{\toe}}
\put(1800,600){\usebox{\GreyCircle}}
}
\put(28800,-24900){
\put(300,900){\usebox{\dynkinf}}
\multiput(0,600)(3600,0){2}{\usebox{\GreyCircle}}
\put(5400,0){\usebox{\aone}}
}
\put(27000,-27000){\vector(0,1){5700}}
\put(24900,-27300){\vector(-1,1){1800}}
\put(28500,-18000){\vector(1,1){1800}}
\put(30300,-22500){\vector(-1,1){1800}}
\put(18900,-19500){\line(1,0){1950}}
\put(21150,-19500){\vector(1,0){1950}}
\end{picture}
}

\newcommand{\strtwoffourtwoone}{
\begin{picture}(33000,16350)
%\multiput(0,0)(33000,0){2}{\line(0,1){16350}}\multiput(0,0)(0,16350){2}{\line(1,0){33000}}
\put(29100,5700){\vector(-2,-1){8400}}
\multiput(4200,5850)(600,-300){14}{\multiput(0,0)(30,-15){10}{\line(1,0){30}}}
\put(12600,1650){\vector(2,-1){300}}
\multiput(13500,5700)(300,-600){5}{\multiput(0,0)(15,-30){10}{\line(0,-1){30}}}
\put(15000,2700){\vector(1,-2){150}}
\multiput(19500,5700)(-300,-600){5}{\multiput(0,0)(-15,-30){10}{\line(0,-1){30}}}
\put(18000,2700){\vector(-1,-2){150}}
\put(21000,13200){\vector(0,-1){3000}}
\put(0,6600){
\put(300,1350){\usebox{\dynkinf}}
\multiput(0,450)(1800,0){4}{\usebox{\aone}}
\multiput(300,2250)(5400,0){2}{\line(0,1){450}}
\put(300,2700){\line(1,0){5400}}
\multiput(2100,0)(3600,0){2}{\line(0,1){450}}
\put(2100,0){\line(1,0){3600}}
\put(700,1750){\usebox{\toe}}
\put(3100,1750){\usebox{\tow}}
}
\put(9000,6600){
\put(300,1350){\usebox{\dynkinf}}
\multiput(0,450)(1800,0){4}{\usebox{\aone}}
\multiput(300,2250)(5400,0){2}{\line(0,1){450}}
\put(300,2700){\line(1,0){5400}}
\multiput(300,0)(3600,0){2}{\line(0,1){450}}
\put(300,0){\line(1,0){3600}}
\multiput(1300,1750)(1800,0){3}{\usebox{\tow}}
}
\put(27000,6600){
\put(300,1350){\usebox{\dynkinf}}
\multiput(0,450)(3600,0){2}{\usebox{\aone}}
\put(5400,450){\usebox{\aone}}
\multiput(300,2250)(5400,0){2}{\line(0,1){450}}
\put(300,2700){\line(1,0){5400}}
\multiput(300,0)(3600,0){2}{\line(0,1){450}}
\put(300,0){\line(1,0){3600}}
\put(700,1750){\usebox{\toe}}
\put(4900,1750){\usebox{\tow}}
\put(1800,1050){\usebox{\GreyCircle}}
}
\put(18000,13650){
\put(300,1350){\usebox{\dynkinf}}
\multiput(0,450)(3600,0){2}{\usebox{\aone}}
\put(5400,450){\usebox{\aone}}
\multiput(300,2250)(5400,0){2}{\line(0,1){450}}
\put(300,2700){\line(1,0){5400}}
\put(700,1750){\usebox{\toe}}
\put(1800,1050){\usebox{\GreyCircle}}
}
\put(18000,6600){
\put(300,1350){\usebox{\dynkinf}}
\multiput(0,450)(5400,0){2}{\usebox{\aone}}
\put(3900,1350){\circle{600}}
\multiput(300,2250)(5400,0){2}{\line(0,1){450}}
\put(300,2700){\line(1,0){5400}}
\put(700,1750){\usebox{\toe}}
\put(4900,1750){\usebox{\tow}}
\put(1800,1050){\usebox{\GreyCircle}}
}
\put(13500,-450){
\put(300,1350){\usebox{\dynkinf}}
\multiput(300,1350)(5400,0){2}{\circle{600}}
\put(3900,1350){\circle{600}}
\multiput(300,450)(5400,0){2}{\line(0,1){600}}
\put(300,450){\line(1,0){5400}}
\put(1800,1050){\usebox{\GreyCircle}}
}
\end{picture}
}

\newcommand{\strtwoffouronefour}{
\begin{picture}(24000,8250)
\put(18000,6000){
\put(300,900){\usebox{\dynkinf}}
\put(0,600){\usebox{\atwo}}
\put(1800,600){\usebox{\GreyCircle}}
\multiput(3600,0)(1800,0){2}{\usebox{\aone}}
\multiput(3100,1300)(1800,0){2}{\usebox{\tow}}
}
\put(9000,6000){
\put(300,900){\usebox{\dynkinf}}
\multiput(0,0)(3600,0){2}{\usebox{\aone}}
\put(5400,0){\usebox{\aone}}
\multiput(300,1800)(5400,0){2}{\line(0,1){450}}
\put(300,2250){\line(1,0){5400}}
\put(700,1300){\usebox{\toe}}
\put(1800,600){\usebox{\GreyCircle}}
}
\put(0,6000){
\put(300,900){\usebox{\dynkinf}}
\multiput(0,600)(3600,0){2}{\usebox{\GreyCircle}}
\put(5400,0){\usebox{\aone}}
}
\put(9000,0){
\put(300,300){\usebox{\dynkinf}}
\multiput(0,0)(3600,0){2}{\usebox{\GreyCircle}}
\put(5700,300){\circle{600}}
}
\put(6000,4800){\vector(1,-1){3000}}
\multiput(12000,5100)(0,-600){5}{\line(0,-1){300}}
\put(12000,2100){\vector(0,-1){300}}
\multiput(18000,4800)(-400,-400){7}{\multiput(0,0)(-20,-20){10}{\line(-1,0){30}}}
\put(15200,2000){\vector(-1,-1){200}}
\end{picture}
}

\newcommand{\stroneffouroneone}{
\begin{picture}(33000,15300)
%\multiput(0,0)(33000,0){2}{\line(0,1){15300}}\multiput(0,0)(0,15300){2}{\line(1,0){33000}}
\put(13500,0){
\put(300,300){\usebox{\dynkinf}}
\put(5400,0){\usebox{\GreyCircle}}
}
\put(0,6450){
\put(300,900){\usebox{\dynkinf}}
\put(3600,600){\usebox{\GreyCircle}}
\put(5400,0){\usebox{\aone}}
\put(4900,1300){\usebox{\tow}}
}
\put(4500,13500){
\put(300,900){\usebox{\dynkinf}}
\put(0,600){\usebox{\atwo}}
\put(1800,600){\usebox{\GreyCircle}}
\multiput(3600,0)(1800,0){2}{\usebox{\aone}}
\put(4250,1300){\usebox{\toe}}
\put(4950,1300){\usebox{\tow}}
}
\put(9000,6450){
\put(300,900){\usebox{\dynkinf}}
\multiput(0,600)(3600,0){2}{\usebox{\atwo}}
\put(1800,600){\usebox{\GreyCircle}}
}
\put(27000,6450){
\put(300,900){\usebox{\dynkinf}}
\multiput(0,600)(3600,0){2}{\usebox{\GreyCircle}}
\put(5400,0){\usebox{\aone}}
}
\put(18000,6000){
\put(300,1350){\usebox{\dynkinf}}
\multiput(0,450)(3600,0){2}{\usebox{\aone}}
\put(5400,450){\usebox{\aone}}
\multiput(300,2250)(5400,0){2}{\line(0,1){450}}
\put(300,2700){\line(1,0){5400}}
\multiput(300,0)(3600,0){2}{\line(0,1){450}}
\put(300,0){\line(1,0){3600}}
\put(700,1750){\usebox{\toe}}
\put(4900,1750){\usebox{\tow}}
\put(1800,1050){\usebox{\GreyCircle}}
}
\multiput(3900,5850)(600,-300){14}{\multiput(0,0)(30,-15){10}{\line(1,0){30}}}
\put(12300,1650){\vector(2,-1){300}}
\multiput(29100,5850)(-600,-300){14}{\multiput(0,0)(-30,-15){10}{\line(-1,0){30}}}
\put(20700,1650){\vector(-2,-1){300}}
\multiput(0,0)(-4500,7200){2}{
\multiput(13500,5400)(300,-600){5}{\multiput(0,0)(15,-30){10}{\line(0,-1){30}}}
\put(15000,2400){\vector(1,-2){150}}
}
\multiput(0,0)(-13500,7200){2}{
\multiput(19500,5400)(-300,-600){5}{\multiput(0,0)(-15,-30){10}{\line(0,-1){30}}}
\put(18000,2400){\vector(-1,-2){150}}
}
\end{picture}
}

\newcommand{\stronecthreeoneone}{
\begin{picture}(30000,33000)(-12000,-16200)
%\multiput(-12000,-16200)(30000,0){2}{\line(0,1){33000}}\multiput(-12000,-16200)(0,33000){2}{\line(1,0){30000}}
\put(0,450){
\put(300,900){\usebox{\dynkinf}}
\put(300,900){\circle{600}}
\put(3600,600){\usebox{\GreyCircle}}
}
\put(-6000,-4800){
\put(300,900){\usebox{\dynkinf}}
\put(300,900){\circle{600}}
\multiput(2100,900)(3600,0){2}{\circle{600}}
\put(3600,0){\usebox{\aone}}
\multiput(2100,-450)(3600,0){2}{\line(0,1){1050}}
\put(2100,-450){\line(1,0){3600}}
\put(3100,1300){\usebox{\tow}}
}
\put(0,-10500){
\put(300,900){\usebox{\dynkinf}}
\put(300,900){\circle{600}}
\multiput(1800,0)(1800,0){3}{\usebox{\aone}}
\multiput(2100,1800)(3600,0){2}{\line(0,1){450}}
\put(2100,2250){\line(1,0){3600}}
\multiput(2450,1300)(1800,0){2}{\usebox{\toe}}
\put(4950,1300){\usebox{\tow}}
}
\put(-10800,7500){
\put(300,900){\usebox{\dynkinf}}
\put(300,900){\circle{600}}
\multiput(1800,0)(1800,0){3}{\usebox{\aone}}
\multiput(2100,1800)(3600,0){2}{\line(0,1){450}}
\put(2100,2250){\line(1,0){3600}}
\put(2450,1300){\usebox{\toe}}
\multiput(3150,1300)(1800,0){2}{\usebox{\tow}}
}
\put(0,7500){
\put(300,900){\usebox{\dynkinf}}
\put(300,900){\circle{600}}
\put(1800,600){\usebox{\GreyCircle}}
\put(3600,600){\usebox{\atwo}}
}
\put(10800,7500){
\put(300,900){\usebox{\dynkinf}}
\put(300,900){\circle{600}}
\put(1800,600){\usebox{\GreyCircle}}
\multiput(3600,0)(1800,0){2}{\usebox{\aone}}
\put(4900,1300){\usebox{\tow}}
}
\put(6000,-4800){
\put(300,900){\usebox{\dynkinf}}
\multiput(300,900)(5400,0){2}{\circle{600}}
\put(3600,600){\usebox{\GreyCircle}}
}
\put(0,7050){
\put(-10800,7500){
\put(300,900){\usebox{\dynkinf}}
\multiput(0,0)(1800,0){4}{\usebox{\aone}}
\multiput(2100,1800)(3600,0){2}{\line(0,1){450}}
\put(2100,2250){\line(1,0){3600}}
\put(2450,1300){\usebox{\toe}}
\multiput(3150,1300)(1800,0){2}{\usebox{\tow}}
}
\put(0,7500){
\put(300,900){\usebox{\dynkinf}}
\put(0,0){\usebox{\aone}}
\put(1800,600){\usebox{\GreyCircle}}
\put(3600,600){\usebox{\atwo}}
}
\put(10800,7500){
\put(300,900){\usebox{\dynkinf}}
\put(0,0){\usebox{\aone}}
\put(1800,600){\usebox{\GreyCircle}}
\multiput(3600,0)(1800,0){2}{\usebox{\aone}}
\put(4900,1300){\usebox{\tow}}
}
}
\put(-6000,-16200){
\put(300,900){\usebox{\dynkinf}}
\multiput(0,0)(1800,0){4}{\usebox{\aone}}
\multiput(2100,1800)(3600,0){2}{\line(0,1){450}}
\put(2100,2250){\line(1,0){3600}}
\multiput(2450,1300)(1800,0){2}{\usebox{\toe}}
\put(4950,1300){\usebox{\tow}}
}
\put(6000,-16200){
\put(300,900){\usebox{\dynkinf}}
\put(0,600){\usebox{\atwo}}
\put(1800,600){\usebox{\GreyCircle}}
\multiput(3600,0)(1800,0){2}{\usebox{\aone}}
\put(4250,1300){\usebox{\toe}}
\put(4950,1300){\usebox{\tow}}
}
\put(-12000,450){
\put(300,900){\usebox{\dynkinf}}
\multiput(0,0)(1800,0){4}{\usebox{\aone}}
\multiput(300,1800)(5400,0){2}{\line(0,1){450}}
\put(300,2250){\line(1,0){5400}}
\multiput(2100,-450)(3600,0){2}{\line(0,1){450}}
\put(2100,-450){\line(1,0){3600}}
\put(700,1300){\usebox{\toe}}
\put(3100,1300){\usebox{\tow}}
}
\put(-12000,-10500){
\put(300,900){\usebox{\dynkinf}}
\multiput(0,0)(1800,0){4}{\usebox{\aone}}
\multiput(300,1800)(5400,0){2}{\line(0,1){450}}
\put(300,2250){\line(1,0){5400}}
\multiput(2100,-450)(3600,0){2}{\line(0,1){450}}
\put(2100,-450){\line(1,0){3600}}
\multiput(700,1300)(3600,0){2}{\usebox{\toe}}
}
\put(12000,450){
\put(300,900){\usebox{\dynkinf}}
\put(300,900){\circle{600}}
\put(3600,600){\usebox{\GreyCircle}}
\put(5400,0){\usebox{\aprime}}
}
\put(12000,-10500){
\put(300,900){\usebox{\dynkinf}}
\put(300,900){\circle{600}}
\put(3600,600){\usebox{\GreyCircle}}
\put(5400,0){\usebox{\aone}}
}
\multiput(-3900,8400)(0,7050){2}{\vector(1,0){3000}}
\multiput(9900,8400)(0,7050){2}{\vector(-1,0){3000}}
\put(3000,6600){\vector(0,-1){3000}}
\multiput(-7800,13650)(10800,0){3}{\vector(0,-1){3000}}
\put(9000,-13050){\vector(0,1){6900}}
\multiput(-900,-2400)(0,-10950){2}{\vector(1,1){1800}}
\multiput(-6900,-7650)(12000,0){2}{\vector(1,1){1800}}
\put(-6900,-900){\vector(1,-1){1800}}
\multiput(900,-7650)(12000,0){2}{\vector(-1,1){1800}}
\put(6900,-2400){\vector(-1,1){1800}}
%\put(12900,-900){\vector(-1,-1){1800}}
\put(-5700,6300){\vector(1,-3){2700}}
\put(11700,6300){\vector(-1,-3){2700}}
\put(11100,1350){\line(-1,0){900}}
\put(9900,1350){\vector(-1,0){3000}}
\end{picture}
}

\newcommand{\stronebthreeonethree}{
\begin{picture}(12600,10800)
%\multiput(0,0)(12600,0){2}{\line(0,1){10800}}\multiput(0,0)(0,10800){2}{\line(1,0){12600}}
\put(0,4200){
\put(300,900){\usebox{\segm}}
\put(2100,900){\usebox{\rightbisegm}}
\put(0,600){\usebox{\atwo}}
\put(3600,0){\usebox{\aprime}}
}
\put(4200,-600){
\put(300,900){\usebox{\segm}}
\put(2100,900){\usebox{\rightbisegm}}
\put(0,600){\usebox{\GreyCircleTwo}}
}
\put(4200,9000){
\put(300,900){\usebox{\segm}}
\put(2100,900){\usebox{\rightbisegm}}
\multiput(0,0)(1800,0){2}{\usebox{\aone}}
\put(1300,1300){\usebox{\tow}}
\put(3600,0){\usebox{\aprime}}
}
\put(8400,4200){
\put(300,900){\usebox{\segm}}
\put(2100,900){\usebox{\rightbisegm}}
\put(0,0){\usebox{\aone}}
\put(1800,600){\usebox{\GreyCircleTwo}}
}
\multiput(9300,3600)(-4200,4800){2}{
\multiput(0,0)(-400,-400){4}{\multiput(0,0)(-20,-20){10}{\line(-1,0){30}}}
\put(-1600,-1600){\vector(-1,-1){200}}
}
\multiput(3300,3600)(4200,4800){2}{
\multiput(0,0)(400,-400){4}{\multiput(0,0)(20,-20){10}{\line(1,0){30}}}
\put(1600,-1600){\vector(1,-1){200}}
}
\end{picture}
}

\newcommand{\strthreebthreefivefour}{
\begin{picture}(15000,8400)
%\multiput(0,0)(15000,0){2}{\line(0,1){8400}}\multiput(0,0)(0,8400){2}{\line(1,0){15000}}
\put(4500,0){
\put(300,900){\usebox{\dynkinf}}
\put(0,600){\usebox{\atwo}}
\put(1800,600){\usebox{\GreyCircle}}
\put(3600,0){\usebox{\aone}}
\put(5700,900){\circle{600}}
}
\put(9000,6600){
\put(300,900){\usebox{\dynkinf}}
\put(0,600){\usebox{\atwo}}
\put(1800,600){\usebox{\GreyCircle}}
\multiput(3600,0)(1800,0){2}{\usebox{\aone}}
\put(4300,1300){\usebox{\toe}}
}
\put(0,6600){
\put(300,900){\usebox{\dynkinf}}
\put(0,600){\usebox{\atwo}}
\put(1800,600){\usebox{\GreyCircle}}
\multiput(3600,0)(1800,0){2}{\usebox{\aone}}
}
\multiput(10650,6000)(-300,-600){5}{\multiput(0,0)(-15,-30){10}{\line(0,-1){30}}}
\put(9150,3000){\vector(-1,-2){150}}
\multiput(4350,6000)(300,-600){5}{\multiput(0,0)(15,-30){10}{\line(0,-1){30}}}
\put(5850,3000){\vector(1,-2){150}}
\end{picture}
}

\newcommand{\stronebthreeonefour}{
\begin{picture}(11400,8400)
%\multiput(0,0)(11400,0){2}{\line(0,1){8400}}\multiput(0,0)(0,8400){2}{\line(1,0){11400}}
\put(0,6150){
\put(300,900){\usebox{\segm}}
\put(2100,900){\usebox{\rightbisegm}}
\multiput(0,0)(1800,0){3}{\usebox{\aone}}
\multiput(300,1800)(1800,0){2}{\line(0,1){450}}
\put(300,2250){\line(1,0){1800}}
\multiput(300,0)(3600,0){2}{\line(0,-1){450}}
\put(300,-450){\line(1,0){3600}}
\put(2500,1300){\usebox{\toe}}
}
\put(3600,-600){
\put(300,900){\usebox{\segm}}
\put(2100,900){\usebox{\rightbisegm}}
\put(3600,600){\usebox{\GreyCircle}}
}
\put(7200,6150){
\put(300,900){\usebox{\segm}}
\put(2100,900){\usebox{\rightbisegm}}
\put(0,600){\usebox{\atwo}}
\put(1800,600){\usebox{\GreyCircle}}
\put(3600,0){\usebox{\aone}}
}
\multiput(8400,5100)(-300,-600){5}{\multiput(0,0)(-15,-30){10}{\line(0,-1){30}}}
\put(6900,2100){\vector(-1,-2){150}}
\multiput(3000,5100)(300,-600){5}{\multiput(0,0)(15,-30){10}{\line(0,-1){30}}}
\put(4500,2100){\vector(1,-2){150}}
\end{picture}
}

\newcommand{\orbitsone}{
\begin{picture}(30000,3600)(0,-3600)\multiput(0,0)(0,-3600){2}{\line(1,0){30000}}
\put(600,-1800){
\put(0,0){\usebox{\dynkinf}}
\put(0,0){\circle{600}}
}
\put(8400,-1800){
\put(0,0){\usebox{\dynkinf}}
\put(-300,-300){\usebox{\GreyCircleTwo}}
\put(5400,0){\circle{600}}
}
\put(16200,-1800){
\put(0,0){\usebox{\dynkinf}}
\put(-300,-300){\usebox{\atwo}}
\multiput(3300,-900)(1800,0){2}{\usebox{\aprime}}
}
\end{picture}
}

\newcommand{\orbitsoneone}{
\begin{picture}(24000,1800)
\put(300,900){
\put(0,0){\usebox{\dynkinf}}
\put(0,0){\circle{600}}
}
\put(9300,900){
\put(0,0){\usebox{\dynkinf}}
\put(-300,-300){\usebox{\GreyCircleTwo}}
\put(5400,0){\circle{600}}
}
\put(18300,900){
\put(0,0){\usebox{\dynkinf}}
\put(-300,-300){\usebox{\atwo}}
\multiput(3300,-900)(1800,0){2}{\usebox{\aprime}}
}
\end{picture}
}

\newcommand{\orbitstwo}{
\begin{picture}(30000,10800)(0,-10800)\multiput(0,0)(0,-10800){2}{\line(1,0){30000}}
\put(600,-1800){
\put(0,0){\usebox{\dynkinf}}
\put(1800,0){\circle{600}}
}
\put(8400,-1800){
\put(0,0){\usebox{\dynkinf}}
\multiput(0,0)(5400,0){2}{\circle{600}}
\put(1500,-300){\usebox{\GreyCircleTwo}}
}
\put(16200,-1800){
\put(0,0){\usebox{\dynkinf}}
\multiput(0,0)(5400,0){2}{\circle{600}}
\multiput(0,-900)(5400,0){2}{\line(0,1){600}}
\put(0,-900){\line(1,0){5400}}
\put(1500,-300){\usebox{\GreyCircle}}
\put(3600,0){\circle{600}}
}
\put(24000,-1800){
\put(0,0){\usebox{\dynkinf}}
\put(-300,-900){\usebox{\aprime}}
\put(1500,-300){\usebox{\GreyCircleTwo}}
\put(5400,0){\circle{600}}
}
\put(600,-5400){
\put(0,0){\usebox{\dynkinf}}
\put(-300,-300){\usebox{\atwo}}
\put(1500,-300){\usebox{\GreyCircle}}
\multiput(3600,0)(1800,0){2}{\circle{600}}
}
\put(8400,-5400){
\put(0,0){\usebox{\dynkinf}}
\multiput(-300,-300)(3600,0){2}{\usebox{\atwo}}
\put(1500,-300){\usebox{\GreyCircle}}
}
\put(16200,-5400){
\put(0,0){\usebox{\dynkinf}}
\put(-300,-300){\usebox{\atwo}}
\put(1500,-300){\usebox{\GreyCircle}}
\put(3300,-900){\usebox{\aprime}}
\put(5400,0){\circle{600}}
}
\put(24000,-5400){
\put(0,0){\usebox{\dynkinf}}
\put(0,0){\circle{600}}
\put(1500,-900){\usebox{\aone}}
\multiput(3300,-900)(1800,0){2}{\usebox{\aprime}}
}
\put(600,-9000){
\put(0,0){\usebox{\dynkinf}}
\multiput(-300,-900)(1800,0){4}{\usebox{\aone}}
\multiput(0,900)(5400,0){2}{\line(0,1){450}}
\put(0,1350){\line(1,0){5400}}
\multiput(0,-900)(3600,0){2}{\line(0,-1){450}}
\put(0,-1350){\line(1,0){3600}}
\multiput(2800,400)(1800,0){2}{\usebox{\tow}}
}
\put(8400,-9000){
\put(0,0){\usebox{\dynkinf}}
\put(-300,-300){\usebox{\atwo}}
\put(1500,-300){\usebox{\GreyCircle}}
\multiput(3300,-900)(1800,0){2}{\usebox{\aone}}
\multiput(2800,400)(1800,0){2}{\usebox{\tow}}
}
\end{picture}
}

\newcommand{\orbitstwoone}{
\begin{picture}(33000,4800)
\put(300,3900){
\put(0,0){\usebox{\dynkinf}}
\put(1800,0){\circle{600}}
}
\put(9300,3900){
\put(0,0){\usebox{\dynkinf}}
\multiput(0,0)(5400,0){2}{\circle{600}}
\put(1500,-300){\usebox{\GreyCircleTwo}}
}
\put(18300,3900){
\put(0,0){\usebox{\dynkinf}}
\multiput(0,0)(5400,0){2}{\circle{600}}
\multiput(0,-900)(5400,0){2}{\line(0,1){600}}
\put(0,-900){\line(1,0){5400}}
\put(1500,-300){\usebox{\GreyCircle}}
\put(3600,0){\circle{600}}
}
\put(27300,3900){
\put(0,0){\usebox{\dynkinf}}
\put(-300,-900){\usebox{\aprime}}
\put(1500,-300){\usebox{\GreyCircleTwo}}
\put(5400,0){\circle{600}}
}
\put(300,900){
\put(0,0){\usebox{\dynkinf}}
\put(-300,-300){\usebox{\atwo}}
\put(1500,-300){\usebox{\GreyCircle}}
\multiput(3600,0)(1800,0){2}{\circle{600}}
}
\put(9300,900){
\put(0,0){\usebox{\dynkinf}}
\multiput(-300,-300)(3600,0){2}{\usebox{\atwo}}
\put(1500,-300){\usebox{\GreyCircle}}
}
\put(18300,900){
\put(0,0){\usebox{\dynkinf}}
\put(-300,-300){\usebox{\atwo}}
\put(1500,-300){\usebox{\GreyCircle}}
\put(3300,-900){\usebox{\aprime}}
\put(5400,0){\circle{600}}
}
\put(27300,900){
\put(0,0){\usebox{\dynkinf}}
\put(-300,-300){\usebox{\atwo}}
\put(1500,-300){\usebox{\GreyCircle}}
\multiput(3300,-900)(1800,0){2}{\usebox{\aone}}
\multiput(2800,400)(1800,0){2}{\usebox{\tow}}
}
\end{picture}
}

\newcommand{\orbitstwotwo}{
\begin{picture}(15000,2700)
\put(300,1350){
\put(0,0){\usebox{\dynkinf}}
\put(0,0){\circle{600}}
\put(1500,-900){\usebox{\aone}}
\multiput(3300,-900)(1800,0){2}{\usebox{\aprime}}
}
\put(9300,1350){
\put(0,0){\usebox{\dynkinf}}
\multiput(-300,-900)(1800,0){4}{\usebox{\aone}}
\multiput(0,900)(5400,0){2}{\line(0,1){450}}
\put(0,1350){\line(1,0){5400}}
\multiput(0,-900)(3600,0){2}{\line(0,-1){450}}
\put(0,-1350){\line(1,0){3600}}
\multiput(2800,400)(1800,0){2}{\usebox{\tow}}
}
\end{picture}
}

\newcommand{\orbitsthree}{
\begin{picture}(30000,7200)(0,-7200)\multiput(0,0)(0,-7200){2}{\line(1,0){30000}}
\put(600,-1800){
\put(0,0){\usebox{\dynkinf}}
\put(3600,0){\circle{600}}
}
\put(8400,-1800){
\put(0,0){\usebox{\dynkinf}}
\put(5400,0){\circle{600}}
\put(3300,-300){\usebox{\GreyCircle}}
}
\put(16200,-1800){
\put(0,0){\usebox{\dynkinf}}
\put(0,0){\circle{600}}
\put(3300,-300){\usebox{\GreyCircle}}
}
\put(24000,-1800){
\put(0,0){\usebox{\dynkinf}}
\put(5100,-900){\usebox{\aone}}
\put(4600,400){\usebox{\tow}}
\put(3300,-300){\usebox{\GreyCircle}}
}
\put(600,-5400){
\put(0,0){\usebox{\dynkinf}}
\put(5400,0){\circle{600}}
\multiput(-300,-300)(3600,0){2}{\usebox{\GreyCircle}}
}
\put(8400,-5400){
\put(0,0){\usebox{\dynkinf}}
\put(0,0){\circle{600}}
\multiput(1800,0)(3600,0){2}{\circle{600}}
\multiput(1800,-300)(3600,0){2}{\line(0,-1){1050}}
\put(1800,-1350){\line(1,0){3600}}
\put(3300,-900){\usebox{\aone}}
}
\put(16200,-5400){
\put(0,0){\usebox{\dynkinf}}
\multiput(-300,-900)(1800,0){4}{\usebox{\aone}}
\multiput(0,900)(5400,0){2}{\line(0,1){450}}
\put(1800,900){\line(0,1){450}}
\put(0,1350){\line(1,0){5400}}
\multiput(0,-900)(3600,0){2}{\line(0,-1){450}}
\put(0,-1350){\line(1,0){3600}}
\put(2200,400){\usebox{\toe}}
\put(4600,400){\usebox{\tow}}
}
\put(24000,-5400){
\put(0,0){\usebox{\dynkinf}}
\put(-300,-300){\usebox{\atwo}}
\put(1500,-300){\usebox{\GreyCircle}}
\multiput(3300,-900)(1800,0){2}{\usebox{\aone}}
\put(4600,400){\usebox{\tow}}
}
\end{picture}
}

\newcommand{\orbitsthreeone}{
\begin{picture}(24000,4200)
\put(300,3300){
\put(0,0){\usebox{\dynkinf}}
\put(3600,0){\circle{600}}
}
\put(9300,3300){
\put(0,0){\usebox{\dynkinf}}
\put(5400,0){\circle{600}}
\put(3300,-300){\usebox{\GreyCircle}}
}
\put(18300,3300){
\put(0,0){\usebox{\dynkinf}}
\put(0,0){\circle{600}}
\put(3300,-300){\usebox{\GreyCircle}}
}
\put(300,900){
\put(0,0){\usebox{\dynkinf}}
\put(5100,-900){\usebox{\aone}}
\put(4600,400){\usebox{\tow}}
\put(3300,-300){\usebox{\GreyCircle}}
}
\put(9300,900){
\put(0,0){\usebox{\dynkinf}}
\put(5400,0){\circle{600}}
\multiput(-300,-300)(3600,0){2}{\usebox{\GreyCircle}}
}
\end{picture}
}

\newcommand{\orbitsthreetwo}{
\begin{picture}(24000,2700)
\put(300,1350){
\put(0,0){\usebox{\dynkinf}}
\put(0,0){\circle{600}}
\multiput(1800,0)(3600,0){2}{\circle{600}}
\multiput(1800,-300)(3600,0){2}{\line(0,-1){1050}}
\put(1800,-1350){\line(1,0){3600}}
\put(3300,-900){\usebox{\aone}}
}
\put(9300,1350){
\put(0,0){\usebox{\dynkinf}}
\multiput(-300,-900)(1800,0){4}{\usebox{\aone}}
\multiput(0,900)(5400,0){2}{\line(0,1){450}}
\put(1800,900){\line(0,1){450}}
\put(0,1350){\line(1,0){5400}}
\multiput(0,-900)(3600,0){2}{\line(0,-1){450}}
\put(0,-1350){\line(1,0){3600}}
\put(2200,400){\usebox{\toe}}
\put(4600,400){\usebox{\tow}}
}
\put(18300,1350){
\put(0,0){\usebox{\dynkinf}}
\put(-300,-300){\usebox{\atwo}}
\put(1500,-300){\usebox{\GreyCircle}}
\multiput(3300,-900)(1800,0){2}{\usebox{\aone}}
\put(4600,400){\usebox{\tow}}
}
\end{picture}
}

\newcommand{\orbitsfour}{
\begin{picture}(30000,3600)(0,-3600)\multiput(0,0)(0,-3600){2}{\line(1,0){30000}}
\put(600,-1800){
\put(0,0){\usebox{\dynkinf}}
\put(5400,0){\circle{600}}
}
\put(8400,-1800){
\put(0,0){\usebox{\dynkinf}}
\put(5100,-300){\usebox{\GreyCircle}}
}
\put(16200,-1800){
\put(0,0){\usebox{\dynkinf}}
\put(3300,-300){\usebox{\GreyCircle}}
\put(5100,-900){\usebox{\aone}}
}
\end{picture}
}

\newcommand{\orbitsfourone}{
\begin{picture}(15000,1200)
\put(300,600){
\put(0,0){\usebox{\dynkinf}}
\put(5400,0){\circle{600}}
}
\put(9300,600){
\put(0,0){\usebox{\dynkinf}}
\put(5100,-300){\usebox{\GreyCircle}}
}
\end{picture}
}

\newcommand{\orbitsfourtwo}{
\begin{picture}(6000,1800)
\put(300,900){
\put(0,0){\usebox{\dynkinf}}
\put(3300,-300){\usebox{\GreyCircle}}
\put(5100,-900){\usebox{\aone}}
}
\end{picture}
}

\newcommand{\bthreeone}{
\begin{picture}(16800,16200)
%\multiput(0,0)(0,16200){2}{\line(1,0){16800}}\multiput(0,0)(16800,0){2}{\line(0,1){16200}}
\put(4200,14400){
\put(300,900){\usebox{\segm}}
\put(2100,900){\usebox{\rightbisegm}}
\multiput(0,0)(1800,0){3}{\usebox{\aone}}
\multiput(1300,1300)(1800,0){2}{\usebox{\tow}}
}
\put(0,9600){
\put(300,900){\usebox{\segm}}
\put(2100,900){\usebox{\rightbisegm}}
\put(0,600){\usebox{\atwo}}
\put(3600,0){\usebox{\aone}}
\put(3100,1300){\usebox{\tow}}
}
\put(8400,9600){
\put(300,900){\usebox{\segm}}
\put(2100,900){\usebox{\rightbisegm}}
\put(0,0){\usebox{\aone}}
\put(1800,600){\usebox{\GreyCircle}}
\put(3900,900){\circle{600}}
}
\put(12600,4800){
\put(300,900){\usebox{\segm}}
\put(2100,900){\usebox{\rightbisegm}}
\put(0,0){\usebox{\aone}}
\put(1800,600){\usebox{\GreyCircle}}
}
\put(4200,4800){
\put(300,900){\usebox{\segm}}
\put(2100,900){\usebox{\rightbisegm}}
\put(0,600){\usebox{\GreyCircle}}
\put(3900,900){\circle{600}}
}
\put(8400,0){
\put(300,900){\usebox{\segm}}
\put(2100,900){\usebox{\rightbisegm}}
\put(0,600){\usebox{\GreyCircle}}
}
\multiput(0,0)(4200,-4800){3}{
\multiput(5100,13800)(-400,-400){4}{\multiput(0,0)(-20,-20){10}{\line(-1,0){30}}}
\put(3500,12200){\vector(-1,-1){200}}
}
\multiput(0,0)(-4200,-4800){2}{
\multiput(0,0)(4200,-4800){2}{
\multiput(7500,13800)(400,-400){4}{\multiput(0,0)(20,-20){10}{\line(1,0){30}}}
\put(9100,12200){\vector(1,-1){200}}
}
}
\end{picture}
}

\newcommand{\stronglysolvabletwo}{
\begin{picture}(24000,2400)
\put(300,900){\usebox{\dynkinf}}
\multiput(0,0)(1800,0){4}{\usebox{\aone}}
\put(-300,2700){
\multiput(600,-900)(3600,0){2}{\line(0,1){600}}
\multiput(2400,-900)(3600,0){2}{\line(0,1){300}}
\put(600,-300){\line(1,0){3600}}
\put(2400,-600){\line(1,0){1700}}
\put(4300,-600){\line(1,0){1700}}
\put(-1800,0){
\put(2800,-1400){\usebox{\toe}}
\put(5200,-1400){\usebox{\tow}}
}
}
\put(9000,0){
\put(300,900){\usebox{\dynkinf}}
\multiput(0,0)(3600,0){2}{\usebox{\aone}}
\multiput(2100,900)(3600,0){2}{\circle{600}}
\multiput(300,1800)(3600,0){2}{\line(0,1){450}}
\put(300,2250){\line(1,0){3600}}
}
\put(18000,600){
\put(300,300){\usebox{\dynkinf}}
\multiput(300,300)(1800,0){4}{\circle{600}}
}
\multiput(6600,900)(9000,0){2}{\vector(1,0){1800}}
\end{picture}
}

\newcommand{\stronglysolvablethree}{
\begin{picture}(33000,2250)
\put(300,900){\usebox{\dynkinf}}
\multiput(0,0)(1800,0){4}{\usebox{\aone}}
\put(-300,2700){
\multiput(2400,-900)(3600,0){2}{\line(0,1){450}}
\put(2400,-450){\line(1,0){3600}}
\put(2800,-1400){\usebox{\toe}}
\put(5200,-1400){\usebox{\tow}}
}
\put(9000,0){
\put(300,900){\usebox{\dynkinf}}
\put(300,900){\circle{600}}
\multiput(1800,0)(1800,0){3}{\usebox{\aone}}
\put(-300,2700){
\multiput(2400,-900)(3600,0){2}{\line(0,1){450}}
\put(2400,-450){\line(1,0){3600}}
\put(2800,-1400){\usebox{\toe}}
\put(5200,-1400){\usebox{\tow}}
}
}
\put(18000,0){
\put(300,900){\usebox{\dynkinf}}
\multiput(300,900)(3600,0){2}{\circle{600}}
\multiput(1800,0)(3600,0){2}{\usebox{\aone}}
\put(-300,2700){
\multiput(2400,-900)(3600,0){2}{\line(0,1){450}}
\put(2400,-450){\line(1,0){3600}}
}
}
\put(27000,600){
\put(300,300){\usebox{\dynkinf}}
\multiput(300,300)(1800,0){4}{\circle{600}}
}
\multiput(6600,900)(9000,0){3}{\vector(1,0){1800}}
\end{picture}
}

\newcommand{\sevenoneone}{
\begin{picture}(18600,14100)
%\multiput(0,0)(18600,0){2}{\line(0,1){14100}}\multiput(0,0)(0,14100){2}{\line(1,0){18600}}
\put(7200,0){
\put(300,300){\usebox{\rightbisegm}}
\put(2100,300){\usebox{\segm}}
\put(1800,0){\usebox{\GreyCircle}}
}
\put(14400,12450){
\put(300,300){\usebox{\rightbisegm}}
\put(2100,300){\usebox{\segm}}
\put(0,0){\usebox{\GreyCircle}}
\multiput(1800,-600)(1800,0){2}{\usebox{\aone}}
\put(3100,700){\usebox{\tow}}
}
\put(7200,6450){
\put(300,300){\usebox{\rightbisegm}}
\put(2100,300){\usebox{\segm}}
\put(0,0){\usebox{\GreyCircle}}
\put(1800,0){\usebox{\atwo}}
}
\put(0,12450){
\put(300,300){\usebox{\rightbisegm}}
\put(2100,300){\usebox{\segm}}
\multiput(0,-600)(1800,0){3}{\usebox{\aone}}
\multiput(300,1200)(3600,0){2}{\line(0,1){450}}
\put(300,1650){\line(1,0){3600}}
\multiput(1350,700)(1800,0){2}{\usebox{\tow}}
\put(750,700){\usebox{\toe}}
}
\put(0,6450){
\put(300,300){\usebox{\rightbisegm}}
\put(2100,300){\usebox{\segm}}
\multiput(300,300)(3600,0){2}{\circle{600}}
\multiput(300,0)(3600,0){2}{\line(0,-1){1050}}
\put(300,-1050){\line(1,0){3600}}
\put(1800,-600){\usebox{\aone}}
\put(1300,700){\usebox{\tow}}
}
\put(7200,12450){
\put(300,300){\usebox{\rightbisegm}}
\put(2100,300){\usebox{\segm}}
\multiput(0,-600)(1800,0){3}{\usebox{\aone}}
\multiput(300,1200)(3600,0){2}{\line(0,1){450}}
\put(300,1650){\line(1,0){3600}}
\put(3150,700){\usebox{\tow}}
\multiput(750,700)(1800,0){2}{\usebox{\toe}}
}
\put(14400,6450){
\put(300,300){\usebox{\rightbisegm}}
\put(2100,300){\usebox{\segm}}
\put(1800,0){\usebox{\GreyCircle}}
\put(3900,300){\circle{600}}
}
\multiput(4200,4800)(0,6450){2}{
\multiput(0,0)(400,-400){7}{\multiput(0,0)(20,-20){10}{\line(1,0){30}}}
\put(2800,-2800){\vector(1,-1){200}}
}
\put(11400,11250){
\multiput(0,0)(400,-400){7}{\multiput(0,0)(20,-20){10}{\line(1,0){30}}}
\put(2800,-2800){\vector(1,-1){200}}
}
\put(7200,11250){
\multiput(0,0)(-400,-400){7}{\multiput(0,0)(-20,-20){10}{\line(-1,0){30}}}
\put(-2800,-2800){\vector(-1,-1){200}}
}
\multiput(14400,11250)(0,-6450){2}{
\multiput(0,0)(-400,-400){7}{\multiput(0,0)(-20,-20){10}{\line(-1,0){30}}}
\put(-2800,-2800){\vector(-1,-1){200}}
}
\put(9100,4800){
\multiput(0,0)(0,-600){4}{\line(0,-1){300}}
\put(0,-2400){\vector(0,-1){300}}
}
\multiput(2100,11250)(14400,0){2}{
\multiput(0,0)(0,-600){4}{\line(0,-1){300}}
\put(0,-2400){\vector(0,-1){300}}
}
\end{picture}
}
 
\newcommand{\constantone}{
\begin{picture}(27600,2700)
\put(0,0){
\put(300,1350){\usebox{\dynkinf}}
\multiput(0,450)(3600,0){2}{\usebox{\aone}}
\put(5400,450){\usebox{\aone}}
\multiput(300,2250)(5400,0){2}{\line(0,1){450}}
\put(300,2700){\line(1,0){5400}}
\multiput(300,0)(3600,0){2}{\line(0,1){450}}
\put(300,0){\line(1,0){3600}}
\put(700,1750){\usebox{\toe}}
\put(4900,1750){\usebox{\tow}}
\put(1800,1050){\usebox{\GreyCircle}}
}
\put(10800,0){
\put(300,1350){\usebox{\dynkinf}}
\multiput(300,1350)(5400,0){2}{\circle{600}}
\put(3900,1350){\circle{600}}
\multiput(300,450)(5400,0){2}{\line(0,1){600}}
\put(300,450){\line(1,0){5400}}
\put(1800,1050){\usebox{\GreyCircle}}
}
\put(21600,1050){
\put(300,300){\usebox{\dynkinf}}
\put(3900,300){\circle{600}}
}
\multiput(6900,1350)(10800,0){2}{\vector(1,0){3000}}
\end{picture}
}

\newcommand{\constanttwo}{
\begin{picture}(27600,1800)
\put(0,600){
\put(300,300){\usebox{\dynkinf}}
\multiput(0,0)(3600,0){2}{\usebox{\GreyCircle}}
\put(5400,-600){\usebox{\aone}}
}
\put(10800,600){
\put(300,300){\usebox{\dynkinf}}
\multiput(0,0)(3600,0){2}{\usebox{\GreyCircle}}
\put(5700,300){\circle{600}}
}
\put(21600,600){
\put(300,300){\usebox{\dynkinf}}
\put(0,0){\usebox{\GreyCircle}}
\put(5700,300){\circle{600}}
}
\multiput(6900,900)(10800,0){2}{\vector(1,0){3000}}
\end{picture}
}

\newcommand{\constantthree}{
\begin{picture}(16800,1800)
\put(0,600){
\put(300,300){\usebox{\dynkinf}}
\put(3600,0){\usebox{\GreyCircle}}
\put(5400,-600){\usebox{\aone}}
\put(4900,700){\usebox{\tow}}
}
\put(10800,600){
\put(300,300){\usebox{\dynkinf}}
\put(5700,300){\circle{600}}
}
\put(6900,900){\vector(1,0){3000}}
\end{picture}
}

\newcommand{\constantfour}{
\begin{picture}(26400,8850)
\multiput(300,7500)(10800,0){2}{\usebox{\dynkinf}}
\put(0,6600){
\multiput(0,0)(3600,0){2}{\usebox{\aone}}
\put(5400,0){\usebox{\aone}}
\multiput(300,1800)(5400,0){2}{\line(0,1){450}}
\put(300,2250){\line(1,0){5400}}
\put(700,1300){\usebox{\toe}}
\put(1800,600){\usebox{\GreyCircle}}
}
\put(10800,6600){
\multiput(0,0)(5400,0){2}{\usebox{\aone}}
\put(3900,900){\circle{600}}
\multiput(300,1800)(5400,0){2}{\line(0,1){450}}
\put(300,2250){\line(1,0){5400}}
\put(700,1300){\usebox{\toe}}
\put(4900,1300){\usebox{\tow}}
\put(1800,600){\usebox{\GreyCircle}}
}
\put(20400,3300){
\put(300,900){\usebox{\dynkinf}}
\put(0,600){\usebox{\GreyCircle}}
\multiput(3900,900)(1800,0){2}{\circle{600}}
}
\put(6900,7500){\vector(1,0){3000}}
\multiput(300,900)(10800,0){2}{\usebox{\dynkinf}}
\put(0,0){
\put(0,600){\usebox{\atwo}}
\put(1800,600){\usebox{\GreyCircle}}
\multiput(3600,0)(1800,0){2}{\usebox{\aone}}
\multiput(3100,1300)(1800,0){2}{\usebox{\tow}}
}
\put(10800,0){
\put(0,600){\usebox{\GreyCircle}}
\put(3900,900){\circle{600}}
\put(5400,0){\usebox{\aone}}
}
\put(6900,900){\vector(1,0){3000}}
\put(17700,1800){\vector(1,1){1800}}
\put(17700,6600){\vector(1,-1){1800}}
\end{picture}
}

\newcommand{\constantfive}{
\begin{picture}(16800,7200)
\put(0,5400){
\put(300,900){\usebox{\dynkinf}}
\put(0,600){\usebox{\atwo}}
\put(1800,600){\usebox{\GreyCircle}}
\multiput(3600,0)(1800,0){2}{\usebox{\aone}}
\put(4250,1300){\usebox{\toe}}
\put(4950,1300){\usebox{\tow}}
}
\multiput(6900,6300)(0,-6000){2}{\vector(1,0){3000}}
\put(10800,5400){
\put(300,900){\usebox{\dynkinf}}
\multiput(300,900)(5400,0){2}{\circle{600}}
\put(3600,600){\usebox{\GreyCircle}}
}
\put(0,-600){
\put(300,900){\usebox{\dynkinf}}
\multiput(0,600)(3600,0){2}{\usebox{\atwo}}
\put(1800,600){\usebox{\GreyCircle}}
}
\put(10800,-600){
\put(300,900){\usebox{\dynkinf}}
\put(300,900){\circle{600}}
\put(3600,600){\usebox{\GreyCircle}}
}
\multiput(3000,4850)(10800,0){2}{
\multiput(0,0)(0,-600){5}{\line(0,-1){300}}
\put(0,-3000){\vector(0,-1){300}}
}
\end{picture}
}

\newcommand{\constantsix}{
\begin{picture}(15600,10100)(0,2250)
\put(0,9650){
\put(300,1350){\usebox{\dynkinf}}
\multiput(0,450)(1800,0){4}{\usebox{\aone}}
\multiput(300,2250)(5400,0){2}{\line(0,1){450}}
\put(300,2700){\line(1,0){5400}}
\multiput(2100,0)(3600,0){2}{\line(0,1){450}}
\put(2100,0){\line(1,0){3600}}
\put(700,1750){\usebox{\toe}}
\put(3100,1750){\usebox{\tow}}
}
\put(0,2550){
\put(300,1350){\usebox{\dynkinf}}
\multiput(0,450)(1800,0){4}{\usebox{\aone}}
\multiput(300,2250)(5400,0){2}{\line(0,1){450}}
\put(300,2700){\line(1,0){5400}}
\multiput(2100,0)(3600,0){2}{\line(0,1){450}}
\put(2100,0){\line(1,0){3600}}
\multiput(700,1750)(3600,0){2}{\usebox{\toe}}
}
\put(9600,6150){
\put(300,1350){\usebox{\dynkinf}}
\put(3600,450){\usebox{\aone}}
\put(300,1350){\circle{600}}
\multiput(2100,1350)(3600,0){2}{\circle{600}}
\multiput(2100,0)(3600,0){2}{\line(0,1){1050}}
\put(2100,0){\line(1,0){3600}}
\put(3100,1750){\usebox{\tow}}
}
%\put(9600,0){
%\put(300,300){\usebox{\dynkinf}}
%\put(0,0){\usebox{\atwo}}
%\put(3900,300){\circle{600}}
%}
\put(6900,10200){\vector(1,-1){1800}}
\put(6900,4800){\vector(1,1){1800}}
\end{picture}
}

\newcommand{\constantseven}{
\begin{picture}(16800,2700)
\put(0,0){
\put(300,1350){\usebox{\dynkinf}}
\multiput(0,450)(1800,0){4}{\usebox{\aone}}
\multiput(300,2250)(5400,0){2}{\line(0,1){450}}
\put(300,2700){\line(1,0){5400}}
\multiput(2100,0)(3600,0){2}{\line(0,1){450}}
\put(2100,0){\line(1,0){3600}}
\multiput(700,1750)(3600,0){2}{\usebox{\toe}}
}
\put(10800,1050){
\put(300,300){\usebox{\dynkinf}}
\put(0,0){\usebox{\atwo}}
\put(3900,300){\circle{600}}
}
\put(6900,1350){\vector(1,0){3000}}
\end{picture}
}

\newcommand{\constanteight}{
\begin{picture}(16800,2700)
\put(0,0){
\put(300,1350){\usebox{\dynkinf}}
\multiput(0,450)(1800,0){4}{\usebox{\aone}}
\multiput(300,2250)(5400,0){2}{\line(0,1){450}}
\put(300,2700){\line(1,0){5400}}
\multiput(300,0)(3600,0){2}{\line(0,1){450}}
\put(300,0){\line(1,0){3600}}
\multiput(1300,1750)(1800,0){3}{\usebox{\tow}}
}
\put(10800,1050){
\put(300,300){\usebox{\dynkinf}}
\put(3600,0){\usebox{\atwo}}
\put(2100,300){\circle{600}}
}
\put(6900,1350){\vector(1,0){3000}}
\end{picture}
}

\newcommand{\increasingone}{
\begin{picture}(18000,13350)
\put(6000,10500){
\put(300,1350){\usebox{\dynkinf}}
\multiput(0,450)(1800,0){4}{\usebox{\aone}}
\multiput(300,2250)(3600,0){2}{\line(0,1){600}}
\multiput(2100,2250)(3600,0){2}{\line(0,1){300}}
\put(300,2850){\line(1,0){3600}}
\put(2100,2550){\line(1,0){1700}}
\put(4000,2550){\line(1,0){1700}}
\multiput(300,0)(5400,0){2}{\line(0,1){450}}
\put(300,0){\line(1,0){5400}}
\put(1300,1750){\usebox{\tow}}
\put(4300,1750){\usebox{\toe}}
}
\put(0,4800){
\put(300,1350){\usebox{\dynkinf}}
\multiput(0,450)(3600,0){2}{\usebox{\aone}}
\put(5400,450){\usebox{\aone}}
\put(2100,1350){\circle{600}}
\multiput(300,2250)(3600,0){2}{\line(0,1){450}}
\put(300,2700){\line(1,0){3600}}
\multiput(300,0)(5400,0){2}{\line(0,1){450}}
\put(300,0){\line(1,0){5400}}
\put(4300,1750){\usebox{\toe}}
}
\put(12000,4800){
\put(300,1350){\usebox{\dynkinf}}
\put(0,450){\usebox{\aone}}
\multiput(1800,450)(3600,0){2}{\usebox{\aone}}
\put(3900,1350){\circle{600}}
\multiput(2100,2250)(3600,0){2}{\line(0,1){450}}
\put(2100,2700){\line(1,0){3600}}
\multiput(300,0)(5400,0){2}{\line(0,1){450}}
\put(300,0){\line(1,0){5400}}
\put(1300,1750){\usebox{\tow}}
}
\put(6000,0){
\put(300,900){\usebox{\dynkinf}}
\multiput(300,900)(1800,0){4}{\circle{600}}
\multiput(300,0)(5400,0){2}{\line(0,1){600}}
\put(300,0){\line(1,0){5400}}
}
\put(9000,9600){\vector(0,-1){6900}}
\put(5100,4200){\vector(1,-1){1800}}
\put(12900,4200){\vector(-1,-1){1800}}
\end{picture}
}

\newcommand{\increasingtwo}{
\begin{picture}(18000,13200)
\put(6000,10500){
\put(300,1350){\usebox{\dynkinf}}
\multiput(0,450)(1800,0){4}{\usebox{\aone}}
\multiput(300,2250)(1800,0){2}{\line(0,1){450}}
\put(5700,2250){\line(0,1){450}}
\put(300,2700){\line(1,0){5400}}
\multiput(300,0)(3600,0){2}{\line(0,1){450}}
\put(300,0){\line(1,0){3600}}
\put(2500,1750){\usebox{\toe}}
\put(4900,1750){\usebox{\tow}}
}
\put(6000,0){
\put(300,900){\usebox{\dynkinf}}
\multiput(300,900)(1800,0){4}{\circle{600}}
\multiput(300,0)(3600,0){2}{\line(0,1){600}}
\put(300,0){\line(1,0){3600}}
}
\put(0,4800){
\put(300,1350){\usebox{\dynkinf}}
\multiput(0,450)(3600,0){2}{\usebox{\aone}}
\put(5400,450){\usebox{\aone}}
\put(2100,1350){\circle{600}}
\multiput(300,2250)(5400,0){2}{\line(0,1){450}}
\put(300,2700){\line(1,0){5400}}
\multiput(300,0)(3600,0){2}{\line(0,1){450}}
\put(300,0){\line(1,0){3600}}
\put(4900,1750){\usebox{\tow}}
}
\put(12000,4800){
\put(300,1350){\usebox{\dynkinf}}
\multiput(0,450)(1800,0){3}{\usebox{\aone}}
\put(5700,1350){\circle{600}}
\multiput(300,2250)(1800,0){2}{\line(0,1){450}}
\put(300,2700){\line(1,0){1800}}
\multiput(300,0)(3600,0){2}{\line(0,1){450}}
\put(300,0){\line(1,0){3600}}
\put(2500,1750){\usebox{\toe}}
}
\put(9000,9600){\vector(0,-1){6900}}
\put(5100,4200){\vector(1,-1){1800}}
\put(12900,4200){\vector(-1,-1){1800}}
\end{picture}
}

\newcommand{\increasingthree}{
\begin{picture}(18000,11850)
\put(6000,9000){
\put(300,1350){\usebox{\dynkinf}}
\multiput(0,450)(1800,0){4}{\usebox{\aone}}
\multiput(300,2250)(3600,0){2}{\line(0,1){600}}
\multiput(2100,2250)(3600,0){2}{\line(0,1){300}}
\put(300,2850){\line(1,0){3600}}
\put(2100,2550){\line(1,0){1700}}
\put(4000,2550){\line(1,0){1700}}
\multiput(0,0)(1800,0){2}{
\put(3150,1750){\usebox{\tow}}
\put(650,1750){\usebox{\toe}}
}
}
\put(0,3750){
\put(300,1350){\usebox{\dynkinf}}
\multiput(1800,450)(1800,0){3}{\usebox{\aone}}
\put(300,1350){\circle{600}}
\multiput(2100,2250)(3600,0){2}{\line(0,1){450}}
\put(2100,2700){\line(1,0){3600}}
\multiput(0,0)(1800,0){2}{
\put(3150,1750){\usebox{\tow}}
}
\put(2450,1750){\usebox{\toe}}
}
\put(12000,3750){
\put(300,1350){\usebox{\dynkinf}}
\multiput(0,450)(1800,0){3}{\usebox{\aone}}
\put(5700,1350){\circle{600}}
\multiput(300,2250)(3600,0){2}{\line(0,1){450}}
\put(300,2700){\line(1,0){3600}}
\multiput(0,0)(1800,0){2}{
\put(650,1750){\usebox{\toe}}
}
\put(3150,1750){\usebox{\tow}}
}
\put(6000,-300){
\put(300,900){\usebox{\dynkinf}}
\multiput(300,900)(1800,0){4}{\circle{600}}
\put(1800,600){\usebox{\GreyCircle}}
}
\put(9000,8550){\vector(0,-1){6300}}
\put(5100,3600){\vector(1,-1){1800}}
\put(12900,3600){\vector(-1,-1){1800}}
\end{picture}
}

\newcommand{\increasingfour}{
\begin{picture}(23400,9600)
\put(7800,8700){
\multiput(300,300)(1800,0){3}{\usebox{\segm}}
\put(5700,300){\usebox{\rightbisegm}}
\multiput(0,0)(1800,0){3}{\usebox{\atwo}}
\put(5400,0){\usebox{\GreyCircle}}
\put(7500,300){\circle{600}}
}
\put(0,4500){
\multiput(300,300)(1800,0){3}{\usebox{\segm}}
\put(5700,300){\usebox{\rightbisegm}}
\multiput(1800,0)(1800,0){2}{\usebox{\atwo}}
\put(5400,0){\usebox{\GreyCircle}}
\multiput(300,300)(7200,0){2}{\circle{600}}
}
\put(15600,4500){
\multiput(300,300)(1800,0){3}{\usebox{\segm}}
\put(5700,300){\usebox{\rightbisegm}}
\multiput(0,0)(1800,0){3}{\usebox{\atwo}}
\put(7500,300){\circle{600}}
}
\put(7800,300){
\multiput(300,300)(1800,0){3}{\usebox{\segm}}
\put(5700,300){\usebox{\rightbisegm}}
\put(3600,0){\usebox{\GreyCircle}}
\multiput(300,300)(7200,0){2}{\circle{600}}
}
\put(11700,7500){\vector(0,-1){5400}}
\put(6900,3600){\vector(1,-1){1800}}
\put(16500,3600){\vector(-1,-1){1800}}
\end{picture}
}

\newcommand{\thebigone}{
%\begin{picture}(90000,30000)
%\multiput(0,0)(4500,0){21}{\line(0,1){25500}}
%\multiput(0,0)(0,25500){2}{\line(1,0){90000}}
\put(0,11850){
\put(300,900){\usebox{\dynkinf}}
\multiput(300,900)(5400,0){2}{\circle{600}}
\put(3600,600){\usebox{\GreyCircle}}
}
\put(9000,11850){
\put(300,900){\usebox{\dynkinf}}
\put(5400,0){\usebox{\aone}}
\put(4900,1300){\usebox{\tow}}
\put(3600,600){\usebox{\GreyCircle}}
}
\put(18000,11850){
\put(300,900){\usebox{\dynkinf}}
\multiput(0,600)(3600,0){2}{\usebox{\atwo}}
\put(1800,600){\usebox{\GreyCircle}}
}
\put(27000,11850){
\put(300,900){\usebox{\dynkinf}}
\put(3600,0){\usebox{\aone}}
\put(3100,1300){\usebox{\tow}}
\put(300,900){\circle{600}}
\multiput(2100,900)(3600,0){2}{\circle{600}}
\multiput(2100,600)(3600,0){2}{\line(0,-1){1050}}
\put(2100,-450){\line(1,0){3600}}
}
\put(36000,15600){
\put(0,-450){
\put(300,1350){\usebox{\dynkinf}}
\multiput(0,450)(3600,0){2}{\usebox{\aone}}
\put(5400,450){\usebox{\aone}}
\multiput(300,2250)(5400,0){2}{\line(0,1){450}}
\put(300,2700){\line(1,0){5400}}
\multiput(300,0)(3600,0){2}{\line(0,1){450}}
\put(300,0){\line(1,0){3600}}
\put(700,1750){\usebox{\toe}}
\put(4900,1750){\usebox{\tow}}
\put(1800,1050){\usebox{\GreyCircle}}
}
}
\put(48000,15600){
\put(0,-450){
\put(300,1350){\usebox{\dynkinf}}
\multiput(0,450)(1800,0){4}{\usebox{\aone}}
\multiput(300,2250)(5400,0){2}{\line(0,1){450}}
\put(300,2700){\line(1,0){5400}}
\multiput(300,0)(3600,0){2}{\line(0,1){450}}
\put(300,0){\line(1,0){3600}}
\multiput(1300,1750)(1800,0){3}{\usebox{\tow}}
}
}
\put(57000,15600){
\put(300,900){\usebox{\dynkinf}}
\put(5400,0){\usebox{\aone}}
\multiput(0,600)(3600,0){2}{\usebox{\GreyCircle}}
}
\put(69000,15600){
\put(300,900){\usebox{\dynkinf}}
\multiput(0,0)(5400,0){2}{\usebox{\aone}}
\put(3900,900){\circle{600}}
\multiput(300,1800)(5400,0){2}{\line(0,1){450}}
\put(300,2250){\line(1,0){5400}}
\put(700,1300){\usebox{\toe}}
\put(4900,1300){\usebox{\tow}}
\put(1800,600){\usebox{\GreyCircle}}
}
\put(81000,15600){
\put(300,900){\usebox{\dynkinf}}
\put(5400,0){\usebox{\aone}}
\put(0,600){\usebox{\GreyCircle}}
\put(3900,900){\circle{600}}
}
\put(9000,19350){
\put(300,900){\usebox{\dynkinf}}
\put(0,600){\usebox{\atwo}}
\put(1800,600){\usebox{\GreyCircle}}
\multiput(3600,0)(1800,0){2}{\usebox{\aone}}
\put(4250,1300){\usebox{\toe}}
\put(4950,1300){\usebox{\tow}}
}
\put(22500,18900){
\put(300,1350){\usebox{\dynkinf}}
\multiput(0,450)(1800,0){4}{\usebox{\aone}}
\multiput(300,2250)(5400,0){2}{\line(0,1){450}}
\put(300,2700){\line(1,0){5400}}
\multiput(2100,0)(3600,0){2}{\line(0,1){450}}
\put(2100,0){\line(1,0){3600}}
\multiput(700,1750)(3600,0){2}{\usebox{\toe}}
}
\put(36000,22650){
\put(300,1350){\usebox{\dynkinf}}
\multiput(0,450)(1800,0){4}{\usebox{\aone}}
\multiput(300,2250)(5400,0){2}{\line(0,1){450}}
\put(300,2700){\line(1,0){5400}}
\multiput(2100,0)(3600,0){2}{\line(0,1){450}}
\put(2100,0){\line(1,0){3600}}
\put(700,1750){\usebox{\toe}}
\put(3100,1750){\usebox{\tow}}
}
\put(63000,23100){
\put(300,900){\usebox{\dynkinf}}
\multiput(0,0)(3600,0){2}{\usebox{\aone}}
\put(5400,0){\usebox{\aone}}
\multiput(300,1800)(5400,0){2}{\line(0,1){450}}
\put(300,2250){\line(1,0){5400}}
\put(700,1300){\usebox{\toe}}
\put(1800,600){\usebox{\GreyCircle}}
}
\put(75000,23100){
\put(300,900){\usebox{\dynkinf}}
\put(0,600){\usebox{\atwo}}
\put(1800,600){\usebox{\GreyCircle}}
\multiput(3600,0)(1800,0){2}{\usebox{\aone}}
\multiput(3100,1300)(1800,0){2}{\usebox{\tow}}
}
\put(9000,4350){
\put(300,900){\usebox{\dynkinf}}
\put(300,900){\circle{600}}
\put(3600,600){\usebox{\GreyCircle}}
}
\put(31500,4350){
\put(300,900){\usebox{\dynkinf}}
\put(5400,600){\usebox{\GreyCircle}}
}
\put(51000,7650){
\put(300,1350){\usebox{\dynkinf}}
\multiput(300,1350)(5400,0){2}{\circle{600}}
\put(3900,1350){\circle{600}}
\multiput(300,450)(5400,0){2}{\line(0,1){600}}
\put(300,450){\line(1,0){5400}}
\put(1800,1050){\usebox{\GreyCircle}}
}
\put(69000,8100){
\put(300,900){\usebox{\dynkinf}}
\put(5700,900){\circle{600}}
\multiput(0,600)(3600,0){2}{\usebox{\GreyCircle}}
}
\put(81000,8100){
\put(300,900){\usebox{\dynkinf}}
\multiput(3900,900)(1800,0){2}{\circle{600}}
\put(0,600){\usebox{\GreyCircle}}
}
\put(75000,600){
\put(300,900){\usebox{\dynkinf}}
\put(5700,900){\circle{600}}
\put(0,600){\usebox{\GreyCircle}}
}
\put(9600,18600){\vector(-1,-1){4200}}
\put(26700,18300){\vector(1,-2){2100}}
\put(36500,22050){\vector(-1,-2){4000}}
\put(18600,11100){\vector(-1,-1){4200}}
\put(41700,14550){\vector(2,-1){9000}}
\put(63750,15000){\vector(1,-1){4500}}
\put(67900,22650){\vector(1,-2){2100}}
\put(79900,22200){\vector(1,-2){2100}}
\put(84000,14850){\vector(0,-1){4200}}
\put(73900,7350){\vector(1,-2){2100}}
\put(73750,14650){\vector(2,-1){8400}}
\put(12000,18750){
\multiput(0,0)(0,-600){7}{\line(0,-1){300}}
\put(0,-4200){\vector(0,-1){300}}
}
\put(14400,18600){
\multiput(0,0)(400,-400){10}{\multiput(0,0)(20,-20){10}{\line(1,0){30}}}
\put(4000,-4000){\vector(1,-1){200}}
}
\put(5400,11100){
\multiput(0,0)(400,-400){10}{\multiput(0,0)(20,-20){10}{\line(1,0){30}}}
\put(4000,-4000){\vector(1,-1){200}}
}
\put(15900,11200){
\multiput(0,0)(600,-200){25}{\multiput(0,0)(30,-10){10}{\line(1,0){30}}}
\put(15000,-5000){\vector(3,-1){300}}
}
\put(23250,11400){
\multiput(0,0)(600,-300){16}{\multiput(0,0)(30,-15){10}{\line(1,0){30}}}
\put(9600,-4800){\vector(2,-1){300}}
}
\put(26400,11550){
\multiput(0,0)(-600,-300){17}{\multiput(0,0)(-30,-15){10}{\line(-1,0){30}}}
\put(-10200,-5100){\vector(-2,-1){300}}
}
\put(38800,14450){
\multiput(0,0)(-300,-600){13}{\multiput(0,0)(-15,-30){10}{\line(0,-1){30}}}
\put(-3750,-7500){\vector(-1,-2){300}}
}
\put(40950,21600){
\multiput(0,0)(400,-400){28}{\multiput(0,0)(20,-20){10}{\line(1,0){30}}}
\put(11200,-11200){\vector(1,-1){200}}
}
\put(52150,14700){
\multiput(0,0)(300,-600){7}{\multiput(0,0)(15,-30){10}{\line(0,-1){30}}}
\put(2100,-4200){\vector(1,-2){150}}
}
\put(56400,15450){
\multiput(0,0)(-600,-300){30}{\multiput(0,0)(-30,-15){10}{\line(-1,0){30}}}
\put(-18000,-9000){\vector(-2,-1){300}}
}
\put(64800,22500){
\multiput(0,0)(300,-600){20}{\multiput(0,0)(15,-30){10}{\line(0,-1){30}}}
\put(6000,-12000){\vector(1,-2){150}}
}
\put(70500,14800){
\multiput(0,0)(-600,-200){21}{\multiput(0,0)(-30,-10){10}{\line(-1,0){30}}}
\put(-12600,-4200){\vector(-3,-1){300}}
}
\put(79200,22500){
\multiput(0,0)(-300,-600){20}{\multiput(0,0)(-15,-30){10}{\line(0,-1){30}}}
\put(-6000,-12000){\vector(-1,-2){150}}
}
\put(82050,7500){
\multiput(0,0)(-300,-600){7}{\multiput(0,0)(-15,-30){10}{\line(0,-1){30}}}
\put(-2100,-4200){\vector(-1,-2){150}}
}
%\end{picture}
}

%% file: tables.1.tex
\bigskip
\begin{tablef} 
Rank 0 spherical systems.
\begin{center}
\begin{picture}(30000,14400)(0,-14400)%\multiput(0,0)(0,-14400){2}{\line(1,0){30000}}
\multiput(0,-1800)(0,-3600){4}{\multiput(600,0)(7800,0){4}{\put(0,0){\usebox{\dynkinf}}}}
\put(600,-1800){\multiput(0,0)(1800,0){4}{\circle{600}}}
\put(8400,-1800){\multiput(0,0)(1800,0){3}{\circle{600}}}
\put(16200,-1800){\multiput(0,0)(1800,0){2}{\circle{600}}\put(5400,0){\circle{600}}}
\put(24000,-1800){\multiput(3600,0)(1800,0){2}{\circle{600}}\put(0,0){\circle{600}}}
\put(600,-5400){\multiput(1800,0)(1800,0){3}{\circle{600}}}
\put(8400,-5400){\multiput(0,0)(1800,0){2}{\circle{600}}}
\put(16200,-5400){\multiput(0,0)(3600,0){2}{\circle{600}}}
\put(24000,-5400){\multiput(0,0)(5400,0){2}{\circle{600}}}
\put(600,-9000){\multiput(1800,0)(1800,0){2}{\circle{600}}}
\put(8400,-9000){\multiput(1800,0)(3600,0){2}{\circle{600}}}
\put(16200,-9000){\multiput(3600,0)(1800,0){2}{\circle{600}}}
\put(24000,-9000){\put(0,0){\circle{600}}}
\put(600,-12600){\put(1800,0){\circle{600}}}
\put(8400,-12600){\put(3600,0){\circle{600}}}
\put(16200,-12600){\put(5400,0){\circle{600}}}
\put(24000,-12600){}
\end{picture}
\end{center}
\end{tablef}

\bigskip
\begin{tablef}\label{tr1a1}
Rank 1 spherical systems with $\mathrm{supp}(\Sigma)$ of type $\mathsf A_1$.
\begin{center}
\begin{picture}(30000,21600)(0,-21600)%\multiput(0,0)(0,-21600){2}{\line(1,0){30000}}
\multiput(0,-1800)(0,-3600){6}{\multiput(600,0)(7800,0){4}{\put(0,0){\usebox{\dynkinf}}}}
\put(600,-1800){\put(-300,-900){\usebox{\aone}}\put(1800,0){\circle{600}}}
\put(8400,-1800){\put(-300,-900){\usebox{\aone}}\multiput(1800,0)(1800,0){2}{\circle{600}}}
\put(16200,-1800){\put(-300,-900){\usebox{\aone}}\multiput(1800,0)(3600,0){2}{\circle{600}}}
\put(24000,-1800){\put(-300,-900){\usebox{\aone}}\multiput(1800,0)(1800,0){3}{\circle{600}}}
\put(600,-5400){\put(1500,-900){\usebox{\aone}}\multiput(0,0)(3600,0){2}{\circle{600}}}
\put(8400,-5400){\put(1500,-900){\usebox{\aone}}\multiput(0,0)(3600,0){2}{\circle{600}}\put(5400,0){\circle{600}}}
\put(16200,-5400){\put(3300,-900){\usebox{\aone}}\multiput(1800,0)(3600,0){2}{\circle{600}}}
\put(24000,-5400){\put(3300,-900){\usebox{\aone}}\multiput(0,0)(1800,0){2}{\circle{600}}\put(5400,0){\circle{600}}}
\put(600,-9000){\put(5100,-900){\usebox{\aone}}\put(3600,0){\circle{600}}}
\put(8400,-9000){\put(5100,-900){\usebox{\aone}}\multiput(0,0)(3600,0){2}{\circle{600}}}
\put(16200,-9000){\put(5100,-900){\usebox{\aone}}\multiput(1800,0)(1800,0){2}{\circle{600}}}
\put(24000,-9000){\put(5100,-900){\usebox{\aone}}\multiput(0,0)(1800,0){3}{\circle{600}}}
\put(600,-12600){\put(-300,-900){\usebox{\aprime}}\put(1800,0){\circle{600}}}
\put(8400,-12600){\put(-300,-900){\usebox{\aprime}}\multiput(1800,0)(1800,0){2}{\circle{600}}}
\put(16200,-12600){\put(-300,-900){\usebox{\aprime}}\multiput(1800,0)(3600,0){2}{\circle{600}}}
\put(24000,-12600){\put(-300,-900){\usebox{\aprime}}\multiput(1800,0)(1800,0){3}{\circle{600}}}
\put(600,-16200){\put(1500,-900){\usebox{\aprime}}\multiput(0,0)(3600,0){2}{\circle{600}}}
\put(8400,-16200){\put(1500,-900){\usebox{\aprime}}\multiput(0,0)(3600,0){2}{\circle{600}}\put(5400,0){\circle{600}}}
\put(16200,-16200){\put(3300,-900){\usebox{\aprime}}\multiput(1800,0)(3600,0){2}{\circle{600}}}
\put(24000,-16200){\put(3300,-900){\usebox{\aprime}}\multiput(0,0)(1800,0){2}{\circle{600}}\put(5400,0){\circle{600}}}
\put(600,-19800){\put(5100,-900){\usebox{\aprime}}\put(3600,0){\circle{600}}}
\put(8400,-19800){\put(5100,-900){\usebox{\aprime}}\multiput(0,0)(3600,0){2}{\circle{600}}}
\put(16200,-19800){\put(5100,-900){\usebox{\aprime}}\multiput(1800,0)(1800,0){2}{\circle{600}}}
\put(24000,-19800){\put(5100,-900){\usebox{\aprime}}\multiput(0,0)(1800,0){3}{\circle{600}}}
\end{picture}
\end{center}
\end{tablef}

\bigskip
\begin{tablef}\label{tr1a1a1}
Rank 1 spherical systems with $\mathrm{supp}(\Sigma)$ of type $\mathsf A_1\times\mathsf A_1$.
\begin{center}
\begin{picture}(30000,3600)(0,-3600)%\multiput(0,0)(0,-3600){2}{\line(1,0){30000}}
\multiput(600,-1800)(7800,0){3}{\put(0,0){\usebox{\dynkinf}}\multiput(0,0)(1800,0){4}{\circle{600}}}
\put(600,-1800){\multiput(0,-300)(3600,0){2}{\line(0,-1){600}}\put(0,-900){\line(1,0){3600}}}
\put(8400,-1800){\multiput(0,-300)(5400,0){2}{\line(0,-1){600}}\put(0,-900){\line(1,0){5400}}}
\put(16200,-1800){\multiput(1800,-300)(3600,0){2}{\line(0,-1){600}}\put(1800,-900){\line(1,0){3600}}}
\end{picture}
\end{center}
\end{tablef}

\bigskip
\begin{tablef}\label{tr1a2}
Rank 1 spherical systems with $\mathrm{supp}(\Sigma)$ of type $\mathsf A_2$.
\begin{center}
\begin{picture}(30000,3600)(0,-3600)%\multiput(0,0)(0,-3600){2}{\line(1,0){30000}}
\multiput(600,-1800)(7800,0){4}{\put(0,0){\usebox{\dynkinf}}}
\put(600,-1800){\put(-300,-300){\usebox{\atwo}}\put(3600,0){\circle{600}}}
\put(8400,-1800){\put(-300,-300){\usebox{\atwo}}\multiput(3600,0)(1800,0){2}{\circle{600}}}
\put(16200,-1800){\put(3300,-300){\usebox{\atwo}}\put(1800,0){\circle{600}}}
\put(24000,-1800){\put(3300,-300){\usebox{\atwo}}\multiput(0,0)(1800,0){2}{\circle{600}}}
\end{picture}
\end{center}
\end{tablef}

\bigskip
\begin{tablef}\label{tr1b2}
Rank 1 spherical systems with $\mathrm{supp}(\Sigma)$ of type $\mathsf B_2$.
\begin{center}
\begin{picture}(30000,3600)(0,-3600)%\multiput(0,0)(0,-3600){2}{\line(1,0){30000}}
\multiput(600,-1800)(7800,0){3}{\put(0,0){\usebox{\dynkinf}}\multiput(0,0)(5400,0){2}{\circle{600}}}
\put(600,-1800){\put(1500,-300){\usebox{\GreyCircle}}}
\put(8400,-1800){\put(1500,-300){\usebox{\GreyCircle}}\put(3600,0){\circle{600}}}
\put(16200,-1800){\put(1500,-300){\usebox{\GreyCircleTwo}}}
\end{picture}
\end{center}
\end{tablef}

\bigskip
\begin{tablef}\label{tr1b3}
Rank 1 spherical systems with $\mathrm{supp}(\Sigma)$ of type $\mathsf B_3$.
\begin{center}
\begin{picture}(30000,3600)(0,-3600)%\multiput(0,0)(0,-3600){2}{\line(1,0){30000}}
\multiput(600,-1800)(7800,0){4}{\put(0,0){\usebox{\dynkinf}}\put(5400,0){\circle{600}}}
\put(600,-1800){\put(-300,-300){\usebox{\GreyCircle}}}
\put(8400,-1800){\put(-300,-300){\usebox{\GreyCircle}}\put(3600,0){\circle{600}}}
\put(16200,-1800){\put(-300,-300){\usebox{\GreyCircleTwo}}}
\put(24000,-1800){\put(3300,-300){\usebox{\GreyCircle}}}
\end{picture}
\end{center}
\end{tablef}

\bigskip
\begin{tablef}\label{tr1c3}
Rank 1 spherical systems with $\mathrm{supp}(\Sigma)$ of type $\mathsf C_3$.
\begin{center}
\begin{picture}(30000,3600)(0,-3600)%\multiput(0,0)(0,-3600){2}{\line(1,0){30000}}
\multiput(600,-1800)(7800,0){2}{\put(0,0){\usebox{\dynkinf}}\put(0,0){\circle{600}}}
\put(600,-1800){\put(3300,-300){\usebox{\GreyCircle}}}
\put(8400,-1800){\put(3300,-300){\usebox{\GreyCircle}}\put(5400,0){\circle{600}}}
\end{picture}
\end{center}
\end{tablef}

\bigskip
\begin{tablef}\label{tr1f4}
Rank 1 spherical systems with $\mathrm{supp}(\Sigma)$ of type $\mathsf F_4$.
\begin{center}
\begin{picture}(30000,3600)(0,-3600)%\multiput(0,0)(0,-3600){2}{\line(1,0){30000}}
\put(600,-1800){\put(0,0){\usebox{\dynkinf}}\put(5100,-300){\usebox{\GreyCircle}}}
\end{picture}
\end{center}
\end{tablef}

\clearpage

\bigskip
\begin{tablef}
Rank 2 spherical systems with $\mathrm{supp}(\Sigma)$ of type $\mathsf A_1\times\mathsf A_1$.
\begin{center}
\begin{picture}(30000,14400)(0,-14400)%\multiput(0,0)(0,-14400){2}{\line(1,0){30000}}
\multiput(0,-1800)(0,-3600){3}{\multiput(600,0)(7800,0){4}{\put(0,0){\usebox{\dynkinf}}}}
\put(600,-1800){\multiput(-300,-900)(3600,0){2}{\usebox{\aone}}\multiput(1800,0)(3600,0){2}{\circle{600}}\multiput(0,900)(3600,0){2}{\line(0,1){450}}\put(0,1350){\line(1,0){3600}}}
\put(8400,-1800){\multiput(-300,-900)(3600,0){2}{\usebox{\aone}}\multiput(1800,0)(3600,0){2}{\circle{600}}}
\put(16200,-1800){\put(-300,-900){\usebox{\aone}}\put(3300,-900){\usebox{\aprime}}\multiput(1800,0)(3600,0){2}{\circle{600}}}
\put(24000,-1800){\put(-300,-900){\usebox{\aprime}}\put(3300,-900){\usebox{\aone}}\multiput(1800,0)(3600,0){2}{\circle{600}}}
\put(600,-5400){\multiput(-300,-900)(3600,0){2}{\usebox{\aprime}}\multiput(1800,0)(3600,0){2}{\circle{600}}}
\put(8400,-5400){\multiput(-300,-900)(5400,0){2}{\usebox{\aone}}\multiput(1800,0)(1800,0){2}{\circle{600}}\multiput(0,900)(5400,0){2}{\line(0,1){450}}\put(0,1350){\line(1,0){5400}}}
\put(16200,-5400){\multiput(-300,-900)(5400,0){2}{\usebox{\aone}}\multiput(1800,0)(1800,0){2}{\circle{600}}}
\put(24000,-5400){\put(-300,-900){\usebox{\aone}}\put(5100,-900){\usebox{\aprime}}\multiput(1800,0)(1800,0){2}{\circle{600}}}
\put(600,-9000){\put(-300,-900){\usebox{\aprime}}\put(5100,-900){\usebox{\aone}}\multiput(1800,0)(1800,0){2}{\circle{600}}}
\put(8400,-9000){\multiput(-300,-900)(5400,0){2}{\usebox{\aprime}}\multiput(1800,0)(1800,0){2}{\circle{600}}}
\put(16200,-9000){\multiput(1500,-900)(3600,0){2}{\usebox{\aone}}\multiput(0,0)(3600,0){2}{\circle{600}}\multiput(1800,900)(3600,0){2}{\line(0,1){450}}\put(1800,1350){\line(1,0){3600}}}
\put(24000,-9000){\multiput(1500,-900)(3600,0){2}{\usebox{\aone}}\multiput(0,0)(3600,0){2}{\circle{600}}}
\multiput(600,-12600)(7800,0){3}{\put(0,0){\usebox{\dynkinf}}}
\put(600,-12600){\put(1500,-900){\usebox{\aone}}\put(5100,-900){\usebox{\aprime}}\multiput(0,0)(3600,0){2}{\circle{600}}}
\put(8400,-12600){\put(1500,-900){\usebox{\aprime}}\put(5100,-900){\usebox{\aone}}\multiput(0,0)(3600,0){2}{\circle{600}}}
\put(16200,-12600){\multiput(1500,-900)(3600,0){2}{\usebox{\aprime}}\multiput(0,0)(3600,0){2}{\circle{600}}}
\end{picture}
\end{center}
\end{tablef}

\bigskip
\begin{tablef}
Rank 2 spherical systems with $\mathrm{supp}(\Sigma)$ of type $\mathsf A_2$.
\begin{center}
\begin{picture}(30000,10800)(0,-10800)%\multiput(0,0)(0,-10800){2}{\line(1,0){30000}}
\multiput(0,-1800)(0,-3600){3}{\multiput(600,0)(7800,0){4}{\put(0,0){\usebox{\dynkinf}}}}
\put(600,-1800){\multiput(-300,-900)(1800,0){2}{\usebox{\aone}}\put(3600,0){\circle{600}}\multiput(0,900)(1800,0){2}{\line(0,1){450}}\put(0,1350){\line(1,0){1800}}}
\put(8400,-1800){\multiput(-300,-900)(1800,0){2}{\usebox{\aone}}\multiput(3600,0)(1800,0){2}{\circle{600}}\multiput(0,900)(1800,0){2}{\line(0,1){450}}\put(0,1350){\line(1,0){1800}}}
\put(16200,-1800){\multiput(-300,-900)(1800,0){2}{\usebox{\aone}}\put(3600,0){\circle{600}}}
\put(24000,-1800){\multiput(-300,-900)(1800,0){2}{\usebox{\aone}}\multiput(3600,0)(1800,0){2}{\circle{600}}}
\put(600,-5400){\multiput(-300,-900)(1800,0){2}{\usebox{\aprime}}\put(3600,0){\circle{600}}}
\put(8400,-5400){\multiput(-300,-900)(1800,0){2}{\usebox{\aprime}}\multiput(3600,0)(1800,0){2}{\circle{600}}}
\put(16200,-5400){\multiput(3300,-900)(1800,0){2}{\usebox{\aone}}\put(1800,0){\circle{600}}\multiput(3600,900)(1800,0){2}{\line(0,1){450}}\put(3600,1350){\line(1,0){1800}}}
\put(24000,-5400){\multiput(3300,-900)(1800,0){2}{\usebox{\aone}}\multiput(0,0)(1800,0){2}{\circle{600}}\multiput(3600,900)(1800,0){2}{\line(0,1){450}}\put(3600,1350){\line(1,0){1800}}}
\put(600,-9000){\multiput(3300,-900)(1800,0){2}{\usebox{\aone}}\put(1800,0){\circle{600}}}
\put(8400,-9000){\multiput(3300,-900)(1800,0){2}{\usebox{\aone}}\multiput(0,0)(1800,0){2}{\circle{600}}}
\put(16200,-9000){\multiput(3300,-900)(1800,0){2}{\usebox{\aprime}}\put(1800,0){\circle{600}}}
\put(24000,-9000){\multiput(3300,-900)(1800,0){2}{\usebox{\aprime}}\multiput(0,0)(1800,0){2}{\circle{600}}}
\end{picture}
\end{center}
\end{tablef}

\bigskip
\begin{tablef}
Rank 2 spherical systems with $\mathrm{supp}(\Sigma)$ of type $\mathsf B_2$.
\begin{center}
\begin{picture}(30000,7200)(0,-7200)%\multiput(0,0)(0,-7200){2}{\line(1,0){30000}}
\multiput(0,-1800)(0,-3600){2}{\multiput(600,0)(7800,0){4}{\put(0,0){\usebox{\dynkinf}}\multiput(0,0)(5400,0){2}{\circle{600}}}}
\put(600,-1800){\multiput(1500,-900)(1800,0){2}{\usebox{\aone}}\multiput(1800,900)(1800,0){2}{\line(0,1){450}}\put(1800,1350){\line(1,0){1800}}}
\put(8400,-1800){\multiput(1500,-900)(1800,0){2}{\usebox{\aone}}}
\put(16200,-1800){\multiput(1500,-900)(1800,0){2}{\usebox{\aone}}\put(2800,400){\usebox{\tow}}}
\put(24000,-1800){\put(1500,-300){\usebox{\GreyCircle}}\put(3300,-900){\usebox{\aone}}}
\put(600,-5400){\put(1500,-300){\usebox{\GreyCircle}}\put(3300,-900){\usebox{\aprime}}}
\put(8400,-5400){\put(1500,-900){\usebox{\aone}}\put(3300,-900){\usebox{\aprime}}}
\put(16200,-5400){\put(1500,-900){\usebox{\aone}}\put(3300,-900){\usebox{\aprime}}\put(2200,400){\usebox{\toe}}}
\put(24000,-5400){\multiput(1500,-900)(1800,0){2}{\usebox{\aprime}}}
\end{picture}
\end{center}
\end{tablef}

\bigskip
\begin{tablef}
Rank 2 spherical systems with $\mathrm{supp}(\Sigma)$ of type $\mathsf A_2\times\mathsf A_1$.
\begin{center}
\begin{picture}(30000,3600)(0,-3600)%\multiput(0,0)(0,-3600){2}{\line(1,0){30000}}
\multiput(600,-1800)(7800,0){4}{\put(0,0){\usebox{\dynkinf}}}
\put(600,-1800){\put(-300,-300){\usebox{\atwo}}\put(3600,0){\circle{600}}\put(5100,-900){\usebox{\aone}}}
\put(8400,-1800){\put(-300,-300){\usebox{\atwo}}\put(3600,0){\circle{600}}\put(5100,-900){\usebox{\aprime}}}
\put(16200,-1800){\put(-300,-900){\usebox{\aone}}\put(1800,0){\circle{600}}\put(3300,-300){\usebox{\atwo}}}
\put(24000,-1800){\put(-300,-900){\usebox{\aprime}}\put(1800,0){\circle{600}}\put(3300,-300){\usebox{\atwo}}}
\end{picture}
\end{center}
\end{tablef}

\bigskip
\begin{tablef}\label{tr2b3}
Rank 2 spherical systems with $\mathrm{supp}(\Sigma)$ of type $\mathsf B_3$.
\begin{center}
\begin{picture}(30000,10800)(0,-10800)%\multiput(0,0)(0,-10800){2}{\line(1,0){30000}}
\multiput(0,-1800)(0,-3600){2}{\multiput(600,0)(7800,0){4}{\put(0,0){\usebox{\dynkinf}}\put(5400,0){\circle{600}}}}
\put(600,-1800){\put(-300,-900){\usebox{\aone}}\put(1500,-300){\usebox{\GreyCircleTwo}}\put(5400,0){\circle{600}}}
\put(8400,-1800){\put(-300,-900){\usebox{\aone}}\put(1500,-300){\usebox{\GreyCircleTwo}}\put(5400,0){\circle{600}}\put(400,400){\usebox{\toe}}}
\put(16200,-1800){\put(-300,-900){\usebox{\aprime}}\put(1500,-300){\usebox{\GreyCircleTwo}}\put(5400,0){\circle{600}}}
\put(24000,-1800){\put(-300,-900){\usebox{\aone}}\put(1500,-300){\usebox{\GreyCircle}}\put(5400,0){\circle{600}}}
\put(600,-5400){\put(-300,-900){\usebox{\aone}}\put(1500,-300){\usebox{\GreyCircle}}\multiput(3600,0)(1800,0){2}{\circle{600}}}
\put(8400,-5400){\put(-300,-300){\usebox{\atwo}}\put(3300,-900){\usebox{\aone}}}
\put(16200,-5400){\put(-300,-300){\usebox{\atwo}}\put(3300,-900){\usebox{\aone}}\put(2800,400){\usebox{\tow}}}
\put(24000,-5400){\put(-300,-300){\usebox{\atwo}}\put(3300,-900){\usebox{\aprime}}}
\put(0,-9000){\put(600,0){\put(0,0){\usebox{\dynkinf}}\put(5400,0){\circle{600}}}}
\put(600,-9000){\put(-300,-300){\usebox{\atwo}}\put(1500,-300){\usebox{\GreyCircle}}\put(3600,0){\circle{600}}}
\end{picture}
\end{center}
\end{tablef}

\bigskip
\begin{tablef}\label{tr2c3}
Rank 2 spherical systems with $\mathrm{supp}(\Sigma)$ of type $\mathsf C_3$.
\begin{center}
\begin{picture}(30000,7200)(0,-7200)%\multiput(0,0)(0,-7200){2}{\line(1,0){30000}}
\multiput(600,-1800)(7800,0){4}{\put(0,0){\usebox{\dynkinf}}\put(0,0){\circle{600}}}
\put(600,-1800){\multiput(1800,0)(3600,0){2}{\circle{600}}\put(3300,-900){\usebox{\aone}}\multiput(1800,-300)(3600,0){2}{\line(0,-1){900}}\put(1800,-1200){\line(1,0){3600}}}
\put(8400,-1800){\multiput(1800,0)(3600,0){2}{\circle{600}}\put(3300,-900){\usebox{\aone}}\multiput(1800,-300)(3600,0){2}{\line(0,-1){900}}\put(1800,-1200){\line(1,0){3600}}\put(4000,400){\usebox{\toe}}}
\put(16200,-1800){\put(3300,-300){\usebox{\GreyCircle}}\put(5100,-900){\usebox{\aone}}}
\put(24000,-1800){\put(3300,-300){\usebox{\GreyCircle}}\put(5100,-900){\usebox{\aprime}}}
\multiput(600,-5400)(7800,0){3}{\put(0,0){\usebox{\dynkinf}}\put(0,0){\circle{600}}}
\put(600,-5400){\put(1500,-900){\usebox{\aone}}\put(3300,-300){\usebox{\atwo}}}
\put(8400,-5400){\put(1500,-300){\usebox{\GreyCircle}}\put(3600,0){\circle{600}}\put(5100,-900){\usebox{\aone}}}
\put(16200,-5400){\put(1500,-300){\usebox{\GreyCircle}}\put(3300,-300){\usebox{\atwo}}}
\end{picture}
\end{center}
\end{tablef}

\bigskip
\begin{tablef}\label{tr2f4}
Rank 2 spherical systems with $\mathrm{supp}(\Sigma)$ of type $\mathsf F_4$.
\begin{center}
\begin{picture}(30000,7200)(0,-7200)%\multiput(0,0)(0,-7200){2}{\line(1,0){30000}}
\multiput(600,-1800)(7800,0){4}{\put(0,0){\usebox{\dynkinf}}}
\put(600,-1800){\put(3300,-300){\usebox{\GreyCircle}}\put(5100,-900){\usebox{\aone}}}
\put(8400,-1800){\put(3300,-300){\usebox{\GreyCircle}}\put(5100,-900){\usebox{\aone}}\put(4600,400){\usebox{\tow}}}
\put(16200,-1800){\put(-300,-300){\usebox{\GreyCircle}}\put(3600,0){\circle{600}}\put(5100,-900){\usebox{\aone}}}
\put(24000,-1800){\multiput(-300,-300)(3600,0){2}{\usebox{\GreyCircle}}\put(5400,0){\circle{600}}}
\multiput(600,-5400)(7800,0){2}{\put(0,0){\usebox{\dynkinf}}}
\put(600,-5400){\multiput(0,0)(1800,0){4}{\circle{600}}\put(1500,-300){\usebox{\GreyCircle}}\multiput(0,-300)(5400,0){2}{\line(0,-1){600}}\put(0,-900){\line(1,0){5400}}}
\put(8400,-5400){\multiput(-300,-300)(3600,0){2}{\usebox{\atwo}}}
\end{picture}
\end{center}
\end{tablef}

\clearpage

\bigskip
\begin{tablef}
Rank 3 spherical systems with $\mathrm{supp}(\Sigma)$ of type $\mathsf A_2\times\mathsf A_1$.
\begin{center}
\begin{picture}(30000,21600)(0,-21600)%\multiput(0,0)(0,-21600){2}{\line(1,0){30000}}
\multiput(0,-1800)(0,-3600){2}{\multiput(600,0)(7800,0){4}{\put(0,0){\usebox{\dynkinf}}\put(3600,0){\circle{600}}}}
\multiput(600,-1800)(7800,0){4}{\multiput(-300,-900)(1800,0){2}{\usebox{\aone}}\put(5100,-900){\usebox{\aone}}}
\multiput(600,-900)(1800,0){2}{\line(0,1){450}}\put(6000,-900){\line(0,1){450}}\put(600,-450){\line(1,0){5400}}
\put(7800,0){\multiput(600,-900)(5400,0){2}{\line(0,1){450}}\put(600,-450){\line(1,0){5400}}\multiput(2400,-2700)(3600,0){2}{\line(0,-1){450}}\put(2400,-3150){\line(1,0){3600}}\put(1000,-1400){\usebox{\toe}}}
\put(15600,0){\multiput(600,-900)(1800,0){2}{\line(0,1){450}}\put(600,-450){\line(1,0){1800}}}
\put(23400,0){\multiput(600,-900)(5400,0){2}{\line(0,1){450}}\put(600,-450){\line(1,0){5400}}}
\multiput(600,-5400)(7800,0){2}{\multiput(-300,-900)(1800,0){2}{\usebox{\aone}}\put(5100,-900){\usebox{\aone}}}
\put(0,-3600){\multiput(2400,-900)(3600,0){2}{\line(0,1){450}}\put(2400,-450){\line(1,0){3600}}}
\put(7800,-3600){}
\multiput(16200,-5400)(7800,0){2}{\multiput(-300,-900)(1800,0){2}{\usebox{\aone}}\put(5100,-900){\usebox{\aprime}}}
\put(15600,-3600){\multiput(600,-900)(1800,0){2}{\line(0,1){450}}\put(600,-450){\line(1,0){1800}}}
\put(23400,-3600){}
\multiput(600,-9000)(7800,0){2}{\put(0,0){\usebox{\dynkinf}}\put(3600,0){\circle{600}}}
\put(600,-9000){\multiput(-300,-900)(1800,0){2}{\usebox{\aprime}}\put(5100,-900){\usebox{\aone}}}
\put(8400,-9000){\multiput(-300,-900)(1800,0){2}{\usebox{\aprime}}\put(5100,-900){\usebox{\aprime}}}
\put(0,-10800){
\multiput(0,-1800)(0,-3600){2}{\multiput(600,0)(7800,0){4}{\put(0,0){\usebox{\dynkinf}}\put(1800,0){\circle{600}}}}
\multiput(600,-1800)(7800,0){4}{\put(-300,-900){\usebox{\aone}}\multiput(3300,-900)(1800,0){2}{\usebox{\aone}}}
\put(600,-900){\line(0,1){450}}\multiput(4200,-900)(1800,0){2}{\line(0,1){450}}\put(600,-450){\line(1,0){5400}}
\put(7800,0){\multiput(600,-900)(5400,0){2}{\line(0,1){450}}\put(600,-450){\line(1,0){5400}}\multiput(600,-2700)(3600,0){2}{\line(0,-1){450}}\put(600,-3150){\line(1,0){3600}}\put(5200,-1400){\usebox{\tow}}}
\put(15600,0){\multiput(600,-900)(3600,0){2}{\line(0,1){450}}\put(600,-450){\line(1,0){3600}}}
\put(23400,0){\multiput(600,-900)(5400,0){2}{\line(0,1){450}}\put(600,-450){\line(1,0){5400}}}
\multiput(600,-5400)(7800,0){2}{\put(-300,-900){\usebox{\aone}}\multiput(3300,-900)(1800,0){2}{\usebox{\aone}}}
\put(0,-3600){\multiput(4200,-900)(1800,0){2}{\line(0,1){450}}\put(4200,-450){\line(1,0){1800}}}
\put(7800,-3600){}
\multiput(16200,-5400)(7800,0){2}{\put(-300,-900){\usebox{\aprime}}\multiput(3300,-900)(1800,0){2}{\usebox{\aone}}}
\put(15600,-3600){\multiput(4200,-900)(1800,0){2}{\line(0,1){450}}\put(4200,-450){\line(1,0){1800}}}
\put(23400,-3600){}
\multiput(600,-9000)(7800,0){2}{\put(0,0){\usebox{\dynkinf}}\put(1800,0){\circle{600}}}
\put(600,-9000){\put(-300,-900){\usebox{\aone}}\multiput(3300,-900)(1800,0){2}{\usebox{\aprime}}}
\put(8400,-9000){\put(-300,-900){\usebox{\aprime}}\multiput(3300,-900)(1800,0){2}{\usebox{\aprime}}}
}
\end{picture}
\end{center}
\end{tablef}

\bigskip
\begin{tablef}\label{tr3b3}
Rank 3 spherical systems with $\mathrm{supp}(\Sigma)$ of type $\mathsf B_3$. 
\begin{center}
\begin{picture}(30000,21600)(0,-21600)%\multiput(0,0)(0,-21600){2}{\line(1,0){30000}}
\multiput(0,-1800)(0,-3600){6}{\multiput(600,0)(7800,0){4}{\put(0,0){\usebox{\dynkinf}}\put(5400,0){\circle{600}}}}
\multiput(0,-1800)(0,-3600){3}{\multiput(600,0)(7800,0){4}{\multiput(-300,-900)(1800,0){3}{\usebox{\aone}}}}
\multiput(600,-900)(1800,0){3}{\line(0,1){450}}\put(600,-450){\line(1,0){3600}}
\put(7800,0){\multiput(600,-900)(1800,0){2}{\line(0,1){450}}\put(600,-450){\line(1,0){1800}}\multiput(600,-2700)(3600,0){2}{\line(0,-1){450}}\put(600,-3150){\line(1,0){3600}}\put(2800,-1400){\usebox{\toe}}}
\put(15600,0){\multiput(600,-900)(1800,0){2}{\line(0,1){450}}\put(600,-450){\line(1,0){1800}}}
\put(23400,0){\multiput(600,-900)(1800,0){2}{\line(0,1){450}}\put(600,-450){\line(1,0){1800}}\put(3400,-1400){\usebox{\tow}}}
\put(0,-3600){\multiput(600,-900)(3600,0){2}{\line(0,1){450}}\put(600,-450){\line(1,0){3600}}}
\put(7800,-3600){\multiput(600,-900)(3600,0){2}{\line(0,1){450}}\put(600,-450){\line(1,0){3600}}\put(1600,-1400){\usebox{\tow}}}
\put(15600,-3600){\multiput(600,-900)(3600,0){2}{\line(0,1){450}}\put(600,-450){\line(1,0){3600}}\put(1000,-1400){\usebox{\toe}}\put(3400,-1400){\usebox{\tow}}}
\put(23400,-3600){\multiput(600,-900)(3600,0){2}{\line(0,1){450}}\put(600,-450){\line(1,0){3600}}\put(950,-1400){\usebox{\toe}}\put(1650,-1400){\usebox{\tow}}\put(3450,-1400){\usebox{\tow}}}
\put(0,-7200){\multiput(2400,-900)(1800,0){2}{\line(0,1){450}}\put(2400,-450){\line(1,0){1800}}}
\put(7800,-7200){}
\put(15600,-7200){\put(3400,-1400){\usebox{\tow}}}
\put(23400,-7200){\put(1600,-1400){\usebox{\tow}}}
\put(600,-12600){\multiput(-300,-900)(1800,0){3}{\usebox{\aone}}}
\put(0,-10800){\put(1600,-1400){\usebox{\tow}}\put(3400,-1400){\usebox{\tow}}}
\multiput(8400,-12600)(7800,0){3}{\multiput(-300,-900)(1800,0){2}{\usebox{\aone}}\put(3300,-900){\usebox{\aprime}}}
\put(7800,-10800){\multiput(600,-900)(1800,0){2}{\line(0,1){450}}\put(600,-450){\line(1,0){1800}}}
\put(15600,-10800){}
\put(23400,-10800){\put(2800,-1400){\usebox{\toe}}}
\put(600,-16200){\multiput(-300,-900)(1800,0){2}{\usebox{\aone}}\put(3300,-900){\usebox{\aprime}}}
\put(0,-14400){\put(1600,-1400){\usebox{\tow}}}
\multiput(8400,-16200)(7800,0){2}{\multiput(-300,-900)(3600,0){2}{\usebox{\aone}}\put(1500,-300){\usebox{\GreyCircle}}}
\put(7800,-14400){\multiput(600,-900)(3600,0){2}{\line(0,1){450}}\put(600,-450){\line(1,0){3600}}}
\put(15600,-14400){}
\put(24000,-16200){\put(-300,-300){\usebox{\atwo}}\put(3300,-900){\usebox{\aone}}\put(1500,-300){\usebox{\GreyCircle}}}
\put(23400,-14400){}
\put(600,-19800){\put(-300,-300){\usebox{\atwo}}\put(3300,-900){\usebox{\aone}}\put(1500,-300){\usebox{\GreyCircle}}}
\put(0,-18000){\put(3400,-1400){\usebox{\tow}}}
\put(8400,-19800){\put(-300,-900){\usebox{\aone}}\put(1500,-300){\usebox{\GreyCircle}}\put(3300,-900){\usebox{\aprime}}}
\put(16200,-19800){\put(-300,-300){\usebox{\atwo}}\put(3300,-900){\usebox{\aprime}}\put(1500,-300){\usebox{\GreyCircle}}}
\put(24000,-19800){\multiput(-300,-900)(1800,0){3}{\usebox{\aprime}}}
\end{picture}
\end{center}
\end{tablef}

\clearpage

\bigskip
\begin{tablef}
Rank 3 spherical systems with $\mathrm{supp}(\Sigma)$ of type $\mathsf C_3$.
\begin{center}
\begin{picture}(30000,18000)(0,-18000)%\multiput(0,0)(0,-18000){2}{\line(1,0){30000}}
\multiput(0,-1800)(0,-3600){4}{\multiput(600,0)(7800,0){4}{\put(0,0){\usebox{\dynkinf}}\put(0,0){\circle{600}}}}
\multiput(0,-1800)(0,-3600){3}{\multiput(600,0)(7800,0){4}{\multiput(1500,-900)(1800,0){3}{\usebox{\aone}}}}
\multiput(2400,-900)(1800,0){3}{\line(0,1){450}}\put(2400,-450){\line(1,0){3600}}
\put(7800,0){\multiput(2400,-900)(1800,0){2}{\line(0,1){450}}\put(2400,-450){\line(1,0){1800}}}
\put(15600,0){\multiput(2400,-900)(3600,0){2}{\line(0,1){450}}\put(2400,-450){\line(1,0){3600}}}
\put(23400,0){\multiput(2400,-900)(3600,0){2}{\line(0,1){450}}\put(2400,-450){\line(1,0){3600}}\put(4600,-1400){\usebox{\toe}}}
\put(0,-3600){\multiput(2400,-900)(3600,0){2}{\line(0,1){450}}\put(2400,-450){\line(1,0){3600}}\put(3400,-1400){\usebox{\tow}}}
\put(7800,-3600){\multiput(2400,-900)(3600,0){2}{\line(0,1){450}}\put(2400,-450){\line(1,0){3600}}\put(2800,-1400){\usebox{\toe}}\put(5200,-1400){\usebox{\tow}}}
\put(15600,-3600){\multiput(2400,-900)(3600,0){2}{\line(0,1){450}}\put(2400,-450){\line(1,0){3600}}\put(4550,-1400){\usebox{\toe}}\put(2750,-1400){\usebox{\toe}}\put(5250,-1400){\usebox{\tow}}}
\put(23400,-3600){\multiput(2400,-900)(3600,0){2}{\line(0,1){450}}\put(2400,-450){\line(1,0){3600}}\put(3450,-1400){\usebox{\tow}}\put(2750,-1400){\usebox{\toe}}\put(5250,-1400){\usebox{\tow}}}
\put(0,-7200){\multiput(4200,-900)(1800,0){2}{\line(0,1){450}}\put(4200,-450){\line(1,0){1800}}}
\put(7800,-7200){}
\put(15600,-7200){\put(4600,-1400){\usebox{\toe}}}
\put(23400,-7200){\put(3400,-1400){\usebox{\tow}}}
\multiput(600,-12600)(7800,0){3}{\put(1500,-300){\usebox{\GreyCircle}}\multiput(3300,-900)(1800,0){2}{\usebox{\aone}}}
\put(0,-10800){\multiput(4200,-900)(1800,0){2}{\line(0,1){450}}\put(4200,-450){\line(1,0){1800}}}
\put(7800,-10800){}
\put(15600,-10800){\put(5200,-1400){\usebox{\tow}}}
\put(24000,-12600){\put(1500,-900){\usebox{\aone}}\multiput(3300,-900)(1800,0){2}{\usebox{\aprime}}}
\put(23400,-10800){}
\multiput(600,-16200)(7800,0){2}{\put(0,0){\usebox{\dynkinf}}\put(0,0){\circle{600}}}
\put(600,-16200){\put(1500,-900){\usebox{\aone}}\multiput(3300,-900)(1800,0){2}{\usebox{\aprime}}}
\put(0,-14400){\put(2800,-1400){\usebox{\toe}}}
\put(8400,-16200){\multiput(1500,-900)(1800,0){3}{\usebox{\aprime}}}
\end{picture}
\end{center}
\end{tablef}

\bigskip
\begin{tablef}\label{tr3f4}
Rank 3 spherical systems with $\mathrm{supp}(\Sigma)$ of type $\mathsf F_4$.
\begin{center}
\begin{picture}(30000,14400)(0,-14400)%\multiput(0,0)(0,-14400){2}{\line(1,0){30000}}
\multiput(0,-1800)(0,-3600){3}{\multiput(600,0)(7800,0){4}{\put(0,0){\usebox{\dynkinf}}}}
\multiput(600,-1800)(7800,0){3}{\multiput(-300,-900)(1800,0){2}{\usebox{\aone}}\put(3300,-300){\usebox{\atwo}}}
\multiput(600,-900)(1800,0){2}{\line(0,1){450}}\put(600,-450){\line(1,0){1800}}
\put(7800,0){}
\put(15600,0){\put(1600,-1400){\usebox{\tow}}}
\put(24000,-1800){\multiput(-300,-900)(5400,0){2}{\usebox{\aone}}\put(1500,-300){\usebox{\GreyCircle}}\put(3600,0){\circle{600}}}
\put(23400,0){\multiput(600,-900)(5400,0){2}{\line(0,1){450}}\put(600,-450){\line(1,0){5400}}}
\multiput(600,-5400)(7800,0){2}{\multiput(-300,-900)(5400,0){2}{\usebox{\aone}}\put(1500,-300){\usebox{\GreyCircle}}\put(3600,0){\circle{600}}}
\put(0,-3600){\multiput(600,-900)(5400,0){2}{\line(0,1){450}}\put(600,-450){\line(1,0){5400}}\put(1000,-1400){\usebox{\toe}}\put(5200,-1400){\usebox{\tow}}}
\put(7800,-3600){}
\multiput(16200,-5400)(7800,0){2}{\multiput(3300,-900)(1800,0){2}{\usebox{\aone}}\put(-300,-300){\usebox{\atwo}}}
\put(15600,-3600){\multiput(4200,-900)(1800,0){2}{\line(0,1){450}}\put(4200,-450){\line(1,0){1800}}}
\put(23400,-3600){}
\multiput(600,-9000)(7800,0){2}{\multiput(3300,-900)(1800,0){2}{\usebox{\aone}}\put(-300,-300){\usebox{\atwo}}}
\put(0,-7200){\put(4600,-1400){\usebox{\toe}}}
\put(7800,-7200){\put(3400,-1400){\usebox{\tow}}}
\put(16200,-9000){\put(-300,-900){\usebox{\aone}}\put(1500,-300){\usebox{\GreyCircle}}\put(3300,-300){\usebox{\atwo}}}
\put(24000,-9000){\put(-300,-300){\usebox{\atwo}}\put(1500,-300){\usebox{\GreyCircle}}\put(3600,0){\circle{600}}\put(5100,-900){\usebox{\aone}}}
\multiput(600,-12600)(7800,0){3}{\put(0,0){\usebox{\dynkinf}}}
\put(600,-12600){\put(-300,-300){\usebox{\atwo}}\put(1500,-300){\usebox{\GreyCircle}}\put(3300,-300){\usebox{\atwo}}}
\put(8400,-12600){\put(-300,-300){\usebox{\atwo}}\multiput(3300,-900)(1800,0){2}{\usebox{\aprime}}}
\put(16200,-12600){\multiput(-300,-300)(3600,0){2}{\usebox{\GreyCircle}}\put(5100,-900){\usebox{\aone}}}
\end{picture}
\end{center}
\end{tablef}

\clearpage

\bigskip
\begin{tablef}\label{tr4ss}
Rank 4 spherical systems, strongly solvable cases.
\begin{center}
\begin{picture}(30000,36000)(0,-36000)%\multiput(0,0)(0,-36000){2}{\line(1,0){30000}}
\multiput(0,-1800)(0,-3600){9}{\multiput(600,0)(7800,0){4}{\put(0,0){\usebox{\dynkinf}}\multiput(-300,-900)(1800,0){4}{\usebox{\aone}}}}
\multiput(600,-900)(1800,0){4}{\line(0,1){450}}\put(600,-450){\line(1,0){5400}}
\put(7800,0){\multiput(600,-900)(1800,0){3}{\line(0,1){450}}\put(600,-450){\line(1,0){3600}}}
\put(15600,0){\multiput(600,-900)(1800,0){2}{\line(0,1){450}}\put(6000,-900){\line(0,1){450}}\put(600,-450){\line(1,0){5400}}}
\put(23400,0){\multiput(600,-900)(1800,0){2}{\line(0,1){450}}\put(6000,-900){\line(0,1){450}}\put(600,-450){\line(1,0){5400}}\put(4200,-1400){\put(400,0){\usebox{\toe}}}}
\put(0,-3600){\multiput(600,-900)(1800,0){2}{\line(0,1){450}}\put(6000,-900){\line(0,1){450}}\put(600,-450){\line(1,0){5400}}\put(4200,-1400){\put(-800,0){\usebox{\tow}}}}
\put(7800,-3600){\multiput(0,0)(3600,0){2}{\multiput(600,-900)(1800,0){2}{\line(0,1){450}}\put(600,-450){\line(1,0){1800}}}}
\put(15600,-3600){\multiput(600,-900)(3600,0){2}{\line(0,1){450}}\put(6000,-900){\line(0,1){450}}\put(600,-450){\line(1,0){5400}}}
\put(23400,-3600){\multiput(600,-900)(3600,0){2}{\line(0,1){450}}\put(6000,-900){\line(0,1){450}}\put(600,-450){\line(1,0){5400}}\put(2400,-1400){\put(-800,0){\usebox{\tow}}}}
\put(0,-7200){\multiput(600,-900)(3600,0){2}{\line(0,1){600}}\multiput(2400,-900)(3600,0){2}{\line(0,1){300}}\put(600,-300){\line(1,0){3600}}\put(2400,-600){\line(1,0){1700}}\put(4300,-600){\line(1,0){1700}}}
\put(7800,-7200){\multiput(600,-900)(3600,0){2}{\line(0,1){600}}\multiput(2400,-900)(3600,0){2}{\line(0,1){300}}\put(600,-300){\line(1,0){3600}}\put(2400,-600){\line(1,0){1700}}\put(4300,-600){\line(1,0){1700}}\put(2800,-1400){\usebox{\toe}}\put(5200,-1400){\usebox{\tow}}}
\put(15600,-7200){\multiput(600,-900)(3600,0){2}{\line(0,1){600}}\multiput(2400,-900)(3600,0){2}{\line(0,1){300}}\put(600,-300){\line(1,0){3600}}\put(2400,-600){\line(1,0){1700}}\put(4300,-600){\line(1,0){1700}}\put(-1800,0){\put(2800,-1400){\usebox{\toe}}\put(5200,-1400){\usebox{\tow}}}}
\put(23400,-7200){\multiput(600,-900)(5400,0){2}{\line(0,1){600}}\multiput(2400,-900)(1800,0){2}{\line(0,1){300}}\put(600,-300){\line(1,0){5400}}\put(2400,-600){\line(1,0){1800}}}
\put(0,-10800){\multiput(2400,-900)(1800,0){3}{\line(0,1){450}}\put(2400,-450){\line(1,0){3600}}}
\put(7800,-10800){\multiput(600,-900)(1800,0){2}{\line(0,1){450}}\put(600,-450){\line(1,0){1800}}}
\put(15600,-10800){\multiput(600,-900)(1800,0){2}{\line(0,1){450}}\put(600,-450){\line(1,0){1800}}\put(4200,-1400){\put(400,0){\usebox{\toe}}}}
\put(23400,-10800){\multiput(600,-900)(1800,0){2}{\line(0,1){450}}\put(600,-450){\line(1,0){1800}}\put(4200,-1400){\put(-800,0){\usebox{\tow}}}}
\put(0,-14400){\multiput(600,-900)(3600,0){2}{\line(0,1){450}}\put(600,-450){\line(1,0){3600}}}
\put(7800,-14400){\multiput(600,-900)(3600,0){2}{\line(0,1){450}}\put(600,-450){\line(1,0){3600}}\put(1600,-1400){\usebox{\tow}}}
\put(15600,-14400){\multiput(600,-900)(3600,0){2}{\line(0,1){450}}\put(600,-450){\line(1,0){3600}}\put(1000,-1400){\usebox{\toe}}\put(3400,-1400){\usebox{\tow}}}
\put(23400,-14400){\multiput(600,-900)(5400,0){2}{\line(0,1){450}}\put(600,-450){\line(1,0){5400}}}
\put(0,-18000){\multiput(600,-900)(5400,0){2}{\line(0,1){450}}\put(600,-450){\line(1,0){5400}}\put(4200,-1400){\put(400,0){\usebox{\toe}}}}
\put(7800,-18000){\multiput(600,-900)(5400,0){2}{\line(0,1){450}}\put(600,-450){\line(1,0){5400}}\put(4200,-1400){\put(-800,0){\usebox{\tow}}}}
\put(15600,-18000){\multiput(600,-900)(5400,0){2}{\line(0,1){450}}\put(600,-450){\line(1,0){5400}}\put(1600,-1400){\usebox{\tow}}}
\put(23400,-18000){\multiput(600,-900)(5400,0){2}{\line(0,1){450}}\put(600,-450){\line(1,0){5400}}\put(1600,-1400){\usebox{\tow}}\put(4200,-1400){\put(400,0){\usebox{\toe}}}}
\put(0,-21600){\multiput(600,-900)(5400,0){2}{\line(0,1){450}}\put(600,-450){\line(1,0){5400}}\put(1600,-1400){\usebox{\tow}}\put(4200,-1400){\put(-800,0){\usebox{\tow}}}}
\put(7800,-21600){\multiput(2400,-900)(1800,0){2}{\line(0,1){450}}\put(2400,-450){\line(1,0){1800}}}
\put(15600,-21600){\multiput(2400,-900)(3600,0){2}{\line(0,1){450}}\put(2400,-450){\line(1,0){3600}}}
\put(23400,-21600){\multiput(2400,-900)(3600,0){2}{\line(0,1){450}}\put(2400,-450){\line(1,0){3600}}\put(4200,-1400){\put(400,0){\usebox{\toe}}}}
\put(0,-25200){\multiput(2400,-900)(3600,0){2}{\line(0,1){450}}\put(2400,-450){\line(1,0){3600}}\put(4200,-1400){\put(-800,0){\usebox{\tow}}}}
\put(7800,-25200){\multiput(2400,-900)(3600,0){2}{\line(0,1){450}}\put(2400,-450){\line(1,0){3600}}\put(2800,-1400){\usebox{\toe}}\put(5200,-1400){\usebox{\tow}}}
\put(15600,-25200){\multiput(4200,-900)(1800,0){2}{\line(0,1){450}}\put(4200,-450){\line(1,0){1800}}}
\put(23400,-25200){\multiput(4200,-900)(1800,0){2}{\line(0,1){450}}\put(4200,-450){\line(1,0){1800}}\put(1600,-1400){\usebox{\tow}}}
\put(0,-28800){}
\put(7800,-28800){\put(4200,-1400){\put(400,0){\usebox{\toe}}}}
\put(15600,-28800){\put(4200,-1400){\put(-800,0){\usebox{\tow}}}}
\put(23400,-28800){\put(1600,-1400){\usebox{\tow}}}
\multiput(600,-34200)(7800,0){2}{\put(0,0){\usebox{\dynkinf}}\multiput(-300,-900)(1800,0){4}{\usebox{\aone}}}
\put(0,-32400){\put(1600,-1400){\usebox{\tow}}\put(4200,-1400){\put(400,0){\usebox{\toe}}}}
\put(7800,-32400){\put(1600,-1400){\usebox{\tow}}\put(4200,-1400){\put(-800,0){\usebox{\tow}}}}
\end{picture}
\end{center}
\end{tablef}

\clearpage

\bigskip
\begin{tablef}\label{tr4o}
Rank 4 spherical systems, other cases.
\begin{center}
\begin{picture}(30000,32400)(0,-32400)%\multiput(0,0)(0,-32400){2}{\line(1,0){30000}}
\multiput(0,-1800)(0,-3600){3}{\multiput(600,0)(7800,0){4}{\put(0,0){\usebox{\dynkinf}}\multiput(-300,-900)(1800,0){4}{\usebox{\aone}}}}
\multiput(600,-900)(1800,0){2}{\line(0,1){450}}\put(6000,-900){\line(0,1){450}}\put(600,-450){\line(1,0){5400}}\multiput(600,-2700)(3600,0){2}{\line(0,-1){450}}\put(600,-3150){\line(1,0){3600}}\put(2800,-1400){\usebox{\toe}}\put(5200,-1400){\usebox{\tow}}
\put(7800,0){\multiput(600,-900)(3600,0){2}{\line(0,1){600}}\multiput(2400,-900)(3600,0){2}{\line(0,1){300}}\put(600,-300){\line(1,0){3600}}\put(2400,-600){\line(1,0){1700}}\put(4300,-600){\line(1,0){1700}}\multiput(600,-2700)(5400,0){2}{\line(0,-1){450}}\put(600,-3150){\line(1,0){5400}}\put(1600,-1400){\usebox{\tow}}\put(4200,-1400){\put(400,0){\usebox{\toe}}}}
\put(15600,0){\multiput(600,-900)(1800,0){2}{\line(0,1){450}}\put(600,-450){\line(1,0){1800}}\multiput(600,-2700)(3600,0){2}{\line(0,-1){450}}\put(600,-3150){\line(1,0){3600}}\multiput(2800,-1400)(1800,0){2}{\usebox{\toe}}}
\put(23400,0){\multiput(600,-900)(5400,0){2}{\line(0,1){450}}\put(600,-450){\line(1,0){5400}}\multiput(600,-2700)(3600,0){2}{\line(0,-1){450}}\put(600,-3150){\line(1,0){3600}}\multiput(3400,-1400)(1800,0){2}{\usebox{\tow}}}
\put(0,-3600){\multiput(600,-900)(5400,0){2}{\line(0,1){450}}\put(600,-450){\line(1,0){5400}}\multiput(600,-2700)(3600,0){2}{\line(0,-1){450}}\put(600,-3150){\line(1,0){3600}}\multiput(1600,-1400)(1800,0){3}{\usebox{\tow}}}
\put(7800,-3600){\multiput(600,-900)(3600,0){2}{\line(0,1){600}}\multiput(2400,-900)(3600,0){2}{\line(0,1){300}}\put(600,-300){\line(1,0){3600}}\put(2400,-600){\line(1,0){1700}}\put(4300,-600){\line(1,0){1700}}\multiput(0,0)(-1800,0){2}{\put(2750,-1400){\usebox{\toe}}\put(5250,-1400){\usebox{\tow}}}}
\put(15600,-3600){\multiput(600,-900)(5400,0){2}{\line(0,1){450}}\put(600,-450){\line(1,0){5400}}\multiput(2400,-2700)(3600,0){2}{\line(0,-1){450}}\put(2400,-3150){\line(1,0){3600}}\put(1000,-1400){\usebox{\toe}}}
\put(23400,-3600){\multiput(600,-900)(5400,0){2}{\line(0,1){450}}\put(600,-450){\line(1,0){5400}}\multiput(2400,-2700)(3600,0){2}{\line(0,-1){450}}\put(2400,-3150){\line(1,0){3600}}\multiput(1000,-1400)(3600,0){2}{\usebox{\toe}}}
\put(0,-7200){\multiput(600,-900)(5400,0){2}{\line(0,1){450}}\put(600,-450){\line(1,0){5400}}\multiput(2400,-2700)(3600,0){2}{\line(0,-1){450}}\put(2400,-3150){\line(1,0){3600}}\put(1000,-1400){\usebox{\toe}}\put(3400,-1400){\usebox{\tow}}}
\put(7800,-7200){\multiput(600,-900)(3600,0){2}{\line(0,1){450}}\put(600,-450){\line(1,0){3600}}\put(950,-1400){\usebox{\toe}}\multiput(3450,-1400)(-1800,0){2}{\usebox{\tow}}}
\put(15600,-7200){\multiput(2400,-900)(3600,0){2}{\line(0,1){450}}\put(2400,-450){\line(1,0){3600}}\multiput(2750,-1400)(1800,0){2}{\usebox{\toe}}\put(5250,-1400){\usebox{\tow}}}
\put(23400,-7200){\multiput(2400,-900)(3600,0){2}{\line(0,1){450}}\put(2400,-450){\line(1,0){3600}}\put(2750,-1400){\usebox{\toe}}\multiput(5250,-1400)(-1800,0){2}{\usebox{\tow}}}
\multiput(0,-12600)(0,-3600){2}{\multiput(600,0)(7800,0){4}{\put(0,0){\usebox{\dynkinf}}\multiput(-300,-900)(3600,0){2}{\usebox{\aone}}\put(5100,-900){\usebox{\aone}}\put(1500,-300){\usebox{\GreyCircle}}}}
\put(0,-10800){\multiput(600,-900)(3600,0){2}{\line(0,1){450}}\put(6000,-900){\line(0,1){450}}\put(600,-450){\line(1,0){5400}}}
\put(7800,-10800){\multiput(600,-900)(5400,0){2}{\line(0,1){450}}\put(600,-450){\line(1,0){5400}}\multiput(600,-2700)(3600,0){2}{\line(0,-1){450}}\put(600,-3150){\line(1,0){3600}}\put(1000,-1400){\usebox{\toe}}\put(5200,-1400){\usebox{\tow}}}
\put(15600,-10800){\multiput(600,-900)(5400,0){2}{\line(0,1){450}}\put(600,-450){\line(1,0){5400}}}
\put(23400,-10800){\multiput(600,-900)(5400,0){2}{\line(0,1){450}}\put(600,-450){\line(1,0){5400}}\put(1000,-1400){\usebox{\toe}}}
\put(0,-14400){\multiput(600,-900)(3600,0){2}{\line(0,1){450}}\put(600,-450){\line(1,0){3600}}}
\put(7800,-14400){\multiput(600,-900)(3600,0){2}{\line(0,1){600}}\put(600,-300){\line(1,0){3600}}
%\put(5700,-900){\multiput(0,0)(-25,25){4}{\line(-1,0){25}}\multiput(-100,100)(-50,25){4}{\line(-1,0){50}}\multiput(-300,200)(-75,25){4}{\line(-1,0){75}}\put(-600,300){\line(-1,0){800}}\put(-1600,300){\line(-1,0){200}}\multiput(-1800,300)(-75,-25){6}{\line(-1,0){75}}\multiput(-2250,150)(-50,-25){6}{\line(-1,0){50}}\multiput(-2550,0)(-25,-25){6}{\line(-1,0){25}}\multiput(-2700,-150)(75,25){6}{\line(1,0){75}}\multiput(-2700,-150)(25,75){6}{\line(0,1){75}}}
\put(5200,-1400){\usebox{\tow}}}
\put(15600,-14400){\multiput(4200,-900)(1800,0){2}{\line(0,1){450}}\put(4200,-450){\line(1,0){1800}}}
\put(23400,-14400){}
\put(600,-19800){\put(0,0){\usebox{\dynkinf}}\multiput(-300,-900)(3600,0){2}{\usebox{\aone}}\put(5100,-900){\usebox{\aone}}\put(1500,-300){\usebox{\GreyCircle}}}
\put(0,-18000){
%\put(5700,-900){\multiput(0,0)(-25,25){4}{\line(-1,0){25}}\multiput(-100,100)(-50,25){4}{\line(-1,0){50}}\multiput(-300,200)(-75,25){4}{\line(-1,0){75}}\put(-600,300){\line(-1,0){1200}}\multiput(-1800,300)(-75,-25){6}{\line(-1,0){75}}\multiput(-2250,150)(-50,-25){6}{\line(-1,0){50}}\multiput(-2550,0)(-25,-25){6}{\line(-1,0){25}}\multiput(-2700,-150)(75,25){6}{\line(1,0){75}}\multiput(-2700,-150)(25,75){6}{\line(0,1){75}}}
\put(5200,-1400){\usebox{\tow}}}
\multiput(8400,-19800)(7800,0){3}{\put(0,0){\usebox{\dynkinf}}\put(-300,-300){\usebox{\atwo}}\put(1500,-300){\usebox{\GreyCircle}}\multiput(3300,-900)(1800,0){2}{\usebox{\aone}}}
\put(7800,-18000){\multiput(4200,-900)(1800,0){2}{\line(0,1){450}}\put(4200,-450){\line(1,0){1800}}}
\put(15600,-18000){}
\put(23400,-18000){\put(5200,-1400){\usebox{\tow}}}
\multiput(600,-23400)(7800,0){4}{\put(0,0){\usebox{\dynkinf}}\put(-300,-300){\usebox{\atwo}}\put(1500,-300){\usebox{\GreyCircle}}\multiput(3300,-900)(1800,0){2}{\usebox{\aone}}}
\put(0,-21600){\put(4600,-1400){\usebox{\toe}}}
\put(7800,-21600){\put(4550,-1400){\usebox{\toe}}\put(5250,-1400){\usebox{\tow}}}
\put(15600,-21600){\put(3400,-1400){\usebox{\tow}}}
\put(23400,-21600){\multiput(3400,-1400)(1800,0){2}{\usebox{\tow}}}
\multiput(600,-27000)(7800,0){4}{\put(0,0){\usebox{\dynkinf}}\multiput(-300,-900)(1800,0){2}{\usebox{\aone}}\multiput(3300,-900)(1800,0){2}{\usebox{\aprime}}}
\put(0,-25200){\multiput(600,-900)(1800,0){2}{\line(0,1){450}}\put(600,-450){\line(1,0){1800}}}
\put(7800,-25200){}
\put(15600,-25200){\put(2800,-1400){\usebox{\toe}}}
\put(23400,-25200){\put(1600,-1400){\usebox{\tow}}}
\put(600,-30600){\put(0,0){\usebox{\dynkinf}}\multiput(-300,-900)(1800,0){4}{\usebox{\aprime}}}
\end{picture}
\end{center}
\end{tablef}